\documentclass[10pt,a4paper]{article}
\usepackage{amsopn}

\usepackage{graphics}
\usepackage{graphicx}
\usepackage{epsfig}

\usepackage{amssymb}
\usepackage{amsthm}
\usepackage{amsmath}
\usepackage{amsfonts}
\usepackage{lipsum}
\usepackage[left=2.5cm,right=2.5cm,bottom=3.5cm,top=3.5cm]{geometry}
\usepackage{empheq}
\usepackage{epstopdf}
\usepackage{changepage}
\usepackage{rotating}
\usepackage{adjustbox}
\usepackage{graphbox}
\usepackage{color}
\usepackage{xcolor}
\usepackage{dsfont}
\usepackage{dutchcal}
\usepackage{hyperref} 
\usepackage{float}
\usepackage{resizegather}
\usepackage{multicol}
\usepackage{booktabs} 
\usepackage{colortbl} 
\usepackage{multirow} 
\usepackage{longtable,tabu} 
\usepackage{hhline}   
\usepackage{anyfontsize} 
\usepackage{everysel}
\usepackage{psfrag}
\usepackage{titlesec}
\usepackage{bbm}
\usepackage{bm}
\usepackage{pgf,tikz}

\newtheorem{remark}{Remark}[section]

\newtheorem*{weakproblem}{Weak problem}

\definecolor{nverde}{RGB}{0,61,0} 
\definecolor{cr1}{RGB}{200,0,0}
\definecolor{cr2}{RGB}{0,0,200}
\definecolor{cr12}{RGB}{100,0,100}

\newcommand{\xx}{\mathbf{x}}
\DeclareMathOperator{\dive}{\nabla \cdot}
\DeclareMathOperator{\gra}{\nabla}
\DeclareMathOperator{\grae}{\nabla}

\DeclareMathOperator{\dV}{\mathrm{dV}}
\newcommand{\g}{\mathbf{g}}
\newcommand{\q}{\mathbf{q}}
\newcommand{\Q}{\mathbf{Q}}
\newcommand{\Flux}{\mathbf{\mathcal{F}}}
\newcommand{\WW}{\mathbf{W}}
\newcommand{\ww}{\mathbf{w}}
\newcommand{\Wur}{\boldsymbol{\mathcal{W}}}
\newcommand{\btau}{\boldsymbol{\tau}}
\newcommand{\bTau}{\mathbf{T}}
\newcommand{\press}{p}
\newcommand{\Press}{P}
\newcommand{\ken}{k}
\newcommand{\Ken}{K}
\newcommand{\rhoken}{\rho k}
\newcommand{\RhoKen}{\rho K}

\newcommand{\vel}{\mathbf{u}}
\newcommand{\wvel}{\mathbf{w}_{\vel}}
\newcommand{\wE}{{w}_{E}}

\newcommand{\wrho}{{w}_{\rho}}
\newcommand{\Vel}{\mathbf{U}}
\newcommand{\Wvel}{\mathbf{W}_{\vel}}
\newcommand{\tWvel}{\widetilde{\mathbf{W}}_{\!\vel}}
\newcommand{\ttWvel}{\widetilde{\widetilde{\mathbf{W}}}_{\!\vel}}
\newcommand{\WE}{{W}_{\!E}}
\newcommand{\WEk}{{W}_{\!E, k}}
\newcommand{\tWE}{\widetilde{{W}}_{\!E}}
\newcommand{\tWEk}{\widetilde{{W}}_{\!E, k}}
\newcommand{\tWEi}{\widetilde{{W}}_{\!E, i}}
\newcommand{\WEi}{{W}_{\!E, i}}
\newcommand{\Wrho}{{W}_{\!\rho}}

\hypersetup{
	pdftitle={A semi-implicit hybrid finite volume / finite element scheme for all Mach number flows on staggered unstructured meshes},    
	pdfauthor={S. Busto, L. Rio-Martin, M.E. Vazquez-Cendon, M. Dumbser},     
	pdfcreator={S. Busto, L. Rio-Martin, M.E. Vazquez-Cendon, M. Dumbser},   
	colorlinks=true,       
}

\allowdisplaybreaks

\begin{document}

\begin{center}
	\textbf{ \Large{A semi-implicit hybrid finite volume / finite element scheme \\ for all Mach number flows on staggered unstructured meshes} }
	
	\vspace{0.5cm}
	{S. Busto\footnote{saray.busto@uvigo.es}, L. R\'io-Mart\'in\footnote{laura.delrio@unitn.it}, M.E. V\'azquez-Cend\'on\footnote{elena.vazquez.cendon@usc.es},  M. Dumbser\footnote{michael.dumbser@unitn.it}}
	
	\vspace{0.2cm}
	{\small
		\textit{$^{(1,2,4)}$ Department of Civil, Environmental and Mechanical Engineering, University of Trento, Via Mesiano 77, 38123 Trento, Italy }
		
		\textit{$^{(3)}$ Department of Applied Mathematics, University of Santiago de Compostela, 15782 Santiago de Compostela, Spain}
	}
\end{center}

\hrule
\vspace{0.4cm}

\begin{center}
	\textbf{Abstract}
\end{center} 

\vspace{0.1cm}
In this paper a new hybrid semi-implicit finite volume / finite element (FV/FE) scheme is presented for the numerical solution of the compressible Euler and Navier-Stokes equations at all Mach numbers on unstructured staggered meshes in two and three space dimensions. The chosen grid arrangement consists of a primal simplex mesh composed of triangles or tetrahedra, and an edge-based / face-based staggered dual mesh. The governing equations are discretized in conservation form. The nonlinear convective terms of the equations, as well as the viscous stress tensor and the heat flux, are discretized on the dual mesh at the aid of an explicit local ADER finite volume scheme, while the implicit pressure terms are discretized at the aid of a continuous $\mathbb{P}^{1}$ finite element method on the nodes of the primal mesh. In the zero Mach number limit, the new scheme automatically reduces to the hybrid FV/FE approach forwarded in \cite{BFTVC17} for the incompressible Navier-Stokes equations. As such, the method is asymptotically consistent with the incompressible limit of the governing equations and can therefore be applied to flows at all Mach numbers. Due to the chosen semi-implicit discretization, the CFL restriction on the time step is only based on the magnitude of the flow velocity and not on the sound speed, hence the method is computationally efficient at low Mach numbers. In the chosen discretization, the only unknown is the scalar pressure field at the new time step. Furthermore, the resulting pressure system is symmetric and positive definite and can therefore be very efficiently solved with a matrix-free conjugate gradient method. 

In order to assess the capabilities of the new scheme, we show computational results for a large set of benchmark problems that range from the quasi incompressible low Mach number regime to compressible flows with shock waves. 

\vspace{0.2cm}

\vspace{0.4cm}

\hrule

\section{Introduction}\label{sec:intro}

Since their first formulation more than 200 years ago, the Euler and Navier-Stokes equations describing the flow of inviscid and viscous fluids have always been a big challenge, both from the theoretical as well as from the numerical point of view. The Euler equations can be directly derived from first principles by considering the conservation of mass, momentum, and total energy. Their extension to the Navier-Stokes equations is then achieved at the aid of appropriate assumptions for the viscous stress tensor and the heat flux. In the most general case, the fluid is assumed to be \textit{compressible}, but different flow regimes can be identified at the aid of the dimensionless \textit{Mach number} $M=\left \|\mathbf{v} \right \|/c$, where $\mathbf{v}$ and $c$ are the fluid velocity and the sound speed, respectively. For $M\to 0$ the behaviour of the fluid becomes the one of an \textit{incompressible} medium, with the well-known condition $\nabla \cdot \mathbf{v} = 0$, which states that the velocity field must become divergence-free when the flow becomes incompressible in the limit $M\to0$. This asymptotic limit was rigorously studied for the first time by Klainerman and Majda in \cite{KlaMaj,KlaMaj82}. 
The asymptotic analysis shows that in the incompressible limit and without compression from the boundary, the pressure can be decomposed into two different contributions: a spatially constant part of the pressure satisfying the equation of state and some fluctuations of the pressure around that constant, governed by the well-known \textit{elliptic} pressure Poisson equation. This change of behaviour for $M \to 0$ is important since the original governing equations are hyperbolic-parabolic they can even exhibit shock waves for high Mach numbers. Because of this changing behaviour of the equations according to the Mach number, it is notoriously difficult to construct suitable numerical schemes which can be simultaneously applied to compressible high Mach number flows with shock waves and also to incompressible or nearly incompressible low Mach number flows. 

Typically, the \textit{incompressible} Euler and Navier-Stokes equations are solved via  \textit{semi-implicit pressure-based} schemes of the finite difference type on staggered grids, see e.g. \cite{HW65,Cho67,Cho68,PS72,Pat80,Kan86,BCG89,HN81,CasulliVOF}, or at the aid of \textit{continuous finite elements} \cite{TH73,BH82,HMM86,For81,Ver91,HR82,HR88}. Instead, for the simulation of \textit{compressible} flows at higher Mach numbers and with shock waves, \textit{explicit density-based} finite volume schemes of the Godunov-type on collocated grids are usually more popular, see \cite{LW60,God59,Roe81,OS82,HLL83,EMPRS9,Munz94,TSS94,LV02,Toro}. 

A first attempt to generalize semi-implicit methods to the more general case of compressible flows was made by Casulli and Greenspan in \cite{CG84}, but the proposed scheme was not conservative and therefore could not be used for the treatment of high Mach number flows with shock waves. Semi-implicit schemes that explicitly make use of the low Mach number asymptotics of the governing partial differential equations can be found in \cite{MeisterMach,munzMPV,Klein95,Klein2001}, while the first conservative staggered semi-implicit pressure-based scheme for compressible flows was introduced by Park and Munz in \cite{PM05}. The scheme \cite{PM05} can be considered as one of the first \textit{all Mach number} flow solvers ever proposed in the literature. The particular splitting of explicit convective terms and implicit pressure terms used in \cite{PM05} was later studied in more detail in \cite{TV12} in order to construct a novel flux-vector splitting method. Since the pioneering work of Park and Munz, the development of all Mach number flow solvers, i.e., of numerical schemes that work at the same time for high Mach number flows with shock waves and in the incompressible limit of the equations, has become a very active research field with many relevant contributions, see e.g. \cite{CordierDegond,DeT11,DC16,RussoAllMach,DLDV18,AbateIolloPuppo,ABIR19,BDLTV2020,BDT2021,BP2021} and references therein. 
For special low Mach number corrections to explicit density-based finite volume schemes, the reader is referred to 
\cite{Thornber,AdamsLowMach}. 

On unstructured simplex meshes, classical continuous finite element methods can nowadays be considered as standard for the numerical solution of the incompressible Navier-Stokes equations. Instead, the construction of \textit{discontinuous} Galerkin finite element schemes for the solution of the compressible and incompressible Navier-Stokes equations on unstructured meshes is still the topic of ongoing research.  
For an overview of high order DG schemes for the compressible and incompressible Navier-Stokes equations, see for example     
\cite{BR97,BO99,BO99b,CS98,CS01,Bassi2007,Ferrer2011,NPC11,RC12,RCV13,CAB13,KKO13}, but this list does not pretend to be complete. 
Concerning high order semi-implicit discontinuous Galerkin methods on \textit{collocated grids} we refer to \cite{Dol08,DF04,DFH07}, while a new family of semi-implicit \textit{staggered} discontinuous Galerkin schemes for the discretization of the incompressible and compressible Navier-Stokes equations was recently forwarded in  \cite{TD15,TD16,TD17,FambriDumbser,AMRDGSI,BTBD20}. 

To round-up this brief literature review, we would also like to point the reader to a very recent and completely different approach for the solution of the Navier-Stokes equations, which consists in embedding the Navier-Stokes equations in a more general first order hyperbolic system with stiff relaxation source terms that is able to describe continuum mechanics as a whole, from nonlinear elastic solids over visco-plastic solids to Newtonian and non-Newtonian fluids, and from which for small enough relaxation times the Navier-Stokes equations are retrieved in the limit of a much more general model that contains continuum mechanics as a whole, see \cite{PeshRom2014,DPRZ16,BCDGP20,SIGPR,GPRNonNewtonian}. This universal model is based on the pioneering work of Godunov and Romenski on symmetric hyperbolic and thermodynamically compatible systems and on nonlinear hyperelasticity, see e.g. \cite{God1961,GodunovRomenski72,God1972MHD,Rom1998,Godunov:2003a} and references therein.     
For an alternative hyperbolic relaxation approach for the discretization of the Navier-Stokes equations, the reader is referred to \cite{NishikawaNS}.

The numerical methods previously discussed were all of a specific type, say finite volume, finite difference, or finite element schemes. More recently in a series of papers a new class of hybrid finite volume / continuous finite element methods on \textit{staggered} unstructured meshes in 2D and 3D has been proposed in \cite{BFSV14,BFTVC17,BSRVC19,Hybrid1} for the solution of the incompressible Navier-Stokes equations and for the low Mach number limit of the weakly compressible equations. In these hybrid schemes, the nonlinear convective part of the equations was solved at the aid of an explicit finite volume scheme on an edge-based staggered dual mesh, see \cite{BTVC16,BFTVC17} for a more detailed analysis, and the pressure equation was solved with a continuous finite element method on the primal grid. The advantage of this hybrid approach is that for each part of the governing PDE system the most appropriate numerical method could be used, since it is well-known that explicit finite volume methods are more suitable for the discretization of nonlinear hyperbolic PDE systems, while the clear strength of continuous finite element methods lies in the discretization of elliptic problems. 

It is therefore the aim of the present paper to provide a novel pressure-based semi-implicit hybrid finite element / finite volume method on staggered unstructured meshes that can solve the compressible Euler and Navier-Stokes equations in a wide range of Mach numbers, which is a very substantial generalization compared to the incompressible and weakly-compressible flow solvers presented in \cite{BFSV14,BFTVC17,Hybrid1}. Following the seminal ideas outlined in \cite{PM05,TV12,DC16}, the nonlinear convective part of the equations will be discretized via an explicit finite volume scheme, while the resulting pressure equation, which is more complex than the simple pressure Poisson equation that typically results from the discretization of the incompressible Navier-Stokes equations, is discretized on the primal mesh at the aid of classical continuous finite elements. 
The semi-implicit discretization allows choosing a time step that is not limited by the sound speed, but only by the velocity magnitude. The new hybrid scheme of this paper is designed to work simultaneously for incompressible and low Mach number flows, as well as for compressible flows including shock waves. For $M \to 0$ the scheme reduces to the hybrid FV/FE method for the incompressible Navier-Stokes equations forwarded in \cite{BFTVC17}. As such, the proposed method is an asymptotic preserving (AP) all Mach number flow solver. 

The rest of the paper is organized as follows: in Section \ref{sec:goveq} we first introduce the governing equations considered in this paper; next, in Section \ref{sec:numdisc} we present their discretization via the new semi-implicit hybrid finite-volume / finite-element scheme on staggered meshes. In Section \ref{sec:numericalresults} we present numerical results for a wide range of Mach numbers, from almost incompressible flows to supersonic flows with shock waves. The conclusions and an outlook to further work are given in Section \ref{sec:conclusions}. 

\section{Governing partial differential equations} \label{sec:goveq}
Let us denote by $\rho$ the density, $\vel=(u,v,w)$ is the velocity vector and $E$ is the specific total energy then, the compressible Navier-Stokes equations given in conservative form read
\begin{gather}
\frac{\partial \rho}{\partial t} + \dive \left( \rho\vel\right)=0,\label{eq:mass1}\\
\frac{\partial \rho\vel}{\partial t} + \dive \left( \rho\vel\otimes \vel\right)  + \grae \press - \dive \btau = \rho \g,\label{eq:momentum1}\\
\frac{\partial \rho E}{\partial t} + \dive \left[ \vel\left(\rho E + \press \right) \right]   - \dive \left(\btau\vel\right) +\dive \q = \rho \g \cdot \vel,\label{eq:totenerg1}
\end{gather}
where $\g$ is the gravity vector, $\btau$ is the tensor of the viscous stresses, 
\begin{equation}
\btau = \mu \left(\gra \mathbf{u} + \gra \vel^{T} \right) -\frac{2}{3} \mu \left( \dive \vel  \right) \mathbf{I}, \label{eq:stresstensor}
\end{equation}
and $\q$ denotes the heat flux,
\begin{equation}
\q = -\lambda \gra \theta. \label{eq:heatflux}
\end{equation}
Here, $\theta$ is the temperature and $\lambda$ denotes the thermal conductivity. 
In this paper we use the simple ideal gas equation of state (EOS) to close the system: 
\begin{equation}
\press = \rho R \theta,\label{eq:stateequation}
\end{equation}
where $R = c_p - c_v$ is the specific gas constant, with $c_{p}$ being the heat capacity at constant pressure, while $c_{v}$ denotes the heat capacity at constant volume. 
Accounting for \eqref{eq:stateequation}, the relation between the total energy, the kinetic energy $k$ and the specific internal energy $e$ reads 
\begin{equation}
\rho E = \rho e +\rho \ken = \frac{1}{\gamma-1}\press + \frac{1}{2}\rho \left| \mathbf{u} \right|^{2} \label{eq:E.as.pk}
\end{equation}
with $\gamma=\frac{c_{\press}}{c_{v}}$ being the ratio of specific heats.  Introducing the enthalpy, 
\begin{equation}
h = e+\frac{\press}{\rho} = \frac{\gamma}{\gamma-1} \frac{\press}{\rho}
\end{equation}
we get
\begin{equation}
\frac{\partial \rho E}{\partial t} + \dive \left( \rho k\vel \right)  + \dive \left( \rho h \vel\right)  - \dive \left(\btau\vel\right) +\dive \q = \rho \g \cdot \vel.
\end{equation}

The former system, \eqref{eq:mass1}-\eqref{eq:totenerg1}, can be rewritten in terms of the conservative variables, $\ww=\left(\wrho, \wvel, \wE \right)^{T}=\left(\rho, \rho\vel, \rho E \right)^{T}$ as
\begin{gather}
\frac{\partial \wrho}{\partial t} + \dive \left( \wvel\right)=0,\label{eq:mass2}\\
\frac{\partial \wvel}{\partial t} + \dive \left(\frac{1}{\rho} \wvel\otimes \wvel\right)  + \grae \press - \dive \btau = \rho \g,\label{eq:momentum2}\\
\frac{\partial \wE}{\partial t} + \dive \left[ \frac{1}{\rho} \wvel\left(\wE + \press \right) \right]   - \dive \left(\frac{1}{\rho} \btau\wvel\right) +\dive \q =  \g \cdot \wvel.\label{eq:totenerg2}
\end{gather}

\section{Numerical method} \label{sec:numdisc}
Discretization of system \eqref{eq:mass2}-\eqref{eq:totenerg2} is performed extending the hybrid finite volume / finite element method proposed in  \cite{BFSV14,BBFSTVC17,BFTVC17,BBDFSVC20}. We start by considering a semi-discrete scheme where only time discretization is applied leading to
\begin{eqnarray}
& &\frac{1}{\Delta t} \left(\Wrho^{n+1}-\Wrho^{n}\right) + \dive \left( \Wvel^{n}\right)=0,\label{eq:mass_td}\\
& &\frac{1}{\Delta t}\left( \Wvel^{n+1} - \Wvel^{n} \right)  + \dive \left(\frac{1}{\rho^{n}} \Wvel^{n}\otimes \Wvel^{n}\right)  + \grae \Press^{n+1} - \dive \bTau^{n} = \rho^{n} \g,\label{eq:momentum_td}\\
& &\frac{1}{\Delta t}\left(  \WE^{n+1} - \WE^{n} \right)  + \dive \left( K^{n} \Wvel^{n} \right)  + \dive \left( H^{n+1} \Wvel^{n+1}\right)   - \dive \left(\frac{1}{\rho^{n}}\bTau^{n}\Wvel^{n}\right) +\dive \Q^{n} =  \g \cdot \Wvel^{n},\quad\label{eq:totenerg_td}
\end{eqnarray}
with $\WW^{n},\, \Press^{n}$ approximations of the solution, $\ww\left( \xx,t^{n}\right) ,\, \press\left( \xx,t^{n}\right) $, at time $t^{n}\in\mathbb{R}^{+}$ and $\xx\in\mathbb{R}^{d}$ the spatial coordinate. We now introduce the following notation for an intermediate approximation of the linear momentum
\begin{equation}
\tWvel = \Wvel^{n} - \Delta t \left( \dive \left(\frac{1}{\rho^{n}} \Wvel^{n}\otimes \Wvel^{n}\right) - \dive \bTau^{n} - \rho^{n} \g \right) \label{eq:wtilde}
\end{equation}
and define
\begin{equation}
\ttWvel := \tWvel -\Delta t \grae \Press^{n}, \qquad \delta\Press^{n+1} := \Press^{n+1}-\Press^{n} \label{eq:ttwvel}
\end{equation}
so
\begin{gather}
 \tWvel = \ttWvel + \Delta t \grae \Press^{n}, \\
 \Wvel^{n+1} = \tWvel - \Delta t \grae \Press^{n+1} = \ttWvel - \Delta t \grae \delta\Press^{n+1} . \label{eq:wn1_dp}
\end{gather}
In such a way, we have derived a pressure-correction formulation in which the computation of the pressure and the linear momentum are ``decoupled".
Similarly, we can define an intermediate auxiliary variable for the computation of the total energy,
\begin{equation}
\tWE = \WE^{n} -\Delta t\left( \dive \left(  K^{n}\Wvel^{n}\right)   - \dive \left(\frac{1}{\rho^{n}}\bTau^{n}\Wvel^{n}\right) +\dive \Q^{n} -  \g \cdot \Wvel^{n}\right).\label{eq:wenerg_td}
\end{equation}
Later, $\WE^{n+1}$ would be recovered from
\begin{equation}
\WE^{n+1} = \tWE - \Delta t \dive \left( H^{n+1}\Wvel^{n+1} \right). \label{eq:wenergt_td}
\end{equation}
On the other hand, equation \eqref{eq:totenerg_td} can be rewritten in terms of the pressure and the kinetic energy by using relation \eqref{eq:E.as.pk} and the ideal gas equation of state as follows:
\begin{equation}
	\frac{\Press^{n+1}}{\gamma-1}+\left(\RhoKen\right)^{n+1} = \frac{\Press^{n}}{\gamma-1}-\frac{\Press^{n}}{\gamma-1}+\tWE - \Delta t \dive \left( H^{n+1} \Wvel^{n+1}\right) 
\end{equation}
Substitution of \eqref{eq:wn1_dp} yields
\begin{gather}
\frac{\Press^{n+1}}{\gamma-1}-\frac{\Press^{n}}{\gamma-1} =-\left(\RhoKen\right)^{n+1} +\tWE -\frac{\Press^{n}}{\gamma-1} - \Delta t \dive \left( H^{n+1} \ttWvel\right)  +\Delta t^{2} \dive\left(  H^{n+1} \grae \delta\Press^{n+1}\right).
\end{gather}
Hence,
\begin{equation}
\frac{1}{\gamma-1} \delta \Press^{n+1}-\Delta t^{2} \dive\left(  H^{n+1} \grae \delta\Press^{n+1}\right) = \tWE-\frac{\Press^{n}}{\gamma-1}-\left(\RhoKen\right)^{n+1} - \Delta t \dive \left( H^{n+1} \ttWvel\right). \label{eq:pressure_dt}
\end{equation}
A Picard procedure is applied to deal with the crossed $t^{n+1}$ terms, i.e., $\Press^{n+1}$ in equation \eqref{eq:wn1_dp}, $ \left(\RhoKen\right)^{n+1}=\frac{1}{2\rho^{n+1}}\left|\Wvel^{n+1}\right|^{2}$ and $H^{n+1}$ in \eqref{eq:pressure_dt}, and $\Wvel^{n+1}$ in \eqref{eq:wenergt_td}, since we do not want to solve a highly nonlinear system. The final system of equations to be discretized in space reads
\begin{eqnarray}
\hspace*{-20pt}& & \Wrho^{n+1} = \Wrho^{n} - \Delta t \dive \left( \Wvel^{n}\right),\label{eq:mass_td_p}\\
\hspace*{-20pt}& &\ttWvel = \Wvel^{n} - \Delta t \left( \dive \left(\frac{1}{\rho^{n}} \Wvel^{n}\otimes \Wvel^{n}\right)  +\grae \Press^{n} - \dive \bTau^{n} - \rho^{n} \g \right), \label{eq:momentumtt_dt_p}\\
\hspace*{-20pt}& &
\tWE = \WE^{n} -\Delta t\left( \dive \left(  K^{n}\Wvel^{n}\right)   - \dive \left(\frac{1}{\rho^{n}}\bTau^{n}\Wvel^{n}\right) +\dive \Q^{n} -  \g \cdot \Wvel^{n}\right),
\label{eq:totenerg_td_p}\\
\hspace*{-20pt}& &\frac{1}{\gamma-1} \delta \Press^{n+1,k+1}-\!\Delta t^{2} \dive\!\left(  H^{n+1,k} \grae \delta\Press^{n+1,k+1}\right) = \! \tWE-\!\frac{\Press^{n}}{\gamma-1}-\!\left(\RhoKen\right)^{n+1,k}\! - \!\Delta t \dive\! \left( H^{n+1,k} \ttWvel\right)\!, \qquad\label{eq:dpressure_dt_p}\\
\hspace*{-20pt}& &\Wvel^{n+1,k+1} = \ttWvel - \Delta t \grae \delta\Press^{n+1,k+1}, \label{eq:momentum_dt_p}\\
\hspace*{-20pt}& & \Press^{n+1,k+1} = \Press^{n}+ \delta\Press^{n+1,k+1},\label{eq:pressure_dt_p}\\
\hspace*{-20pt}& &\WE^{n+1} = \tWE - \Delta t \dive \left( H^{n+1,k+1}\Wvel^{n+1,k+1} \right), \label{eq:energy_dt_p}
\end{eqnarray}
with $k$ the Picard iteration index, $k=1,\dots,N_{\mathrm{Pic}}$. 

Let us remark that the method is by construction asymptotic preserving in the low Mach number limit. 
In the limit $M\rightarrow 0$, we have $c^{2}\rightarrow \infty$ with $c^{2}=\frac{H}{\gamma-1}$. According to \cite{KlaMaj,KlaMaj82}, the pressure, and for constant density also the enthalpy, tend to a constant. In this limit, we can now divide equation \eqref{eq:dpressure_dt_p} by the enthalpy and neglecting terms of the order $1/c^2$ we obtain the following equation
\begin{equation} 
 \nabla^2 \delta\Press^{n+1} = \frac{1}{\Delta t} \dive \left(\ttWvel\right), 
\end{equation}
which, together with the momentum equation \eqref{eq:momentumtt_dt_p}, corresponds to the pressure correction system obtained for the incompressible Navier-Stokes equations in \cite{BFTVC17}.  

For general equations of state, relation \eqref{eq:pressure_dt} needs to be replaced by 
\begin{equation}
W_E(\Press^{n+1},\Wrho^{n+1})  -\Delta t^{2} \dive\left(  H^{n+1} \grae \Press^{n+1}\right) = \tWE-\left(\RhoKen\right)^{n+1} - \Delta t \dive \left( H^{n+1} \ttWvel\right), \label{eq:pressure:GEOS}
\end{equation}
where the density at the new time $t^{n+1}$ is readily available from eqn. \eqref{eq:mass_td_p}, and thus the only unknown remains the scalar pressure field $\Press^{n+1}$ at the new time, see also \cite{DC16}. When appropriate mass-lumping is used within a finite-element discretization of \eqref{eq:pressure:GEOS}, the resulting mildly nonlinear pressure system can be solved very efficiently at the aid of the (nested) Newton-type methods of Brugnano and Casulli \cite{BrugnanoCasulli,BrugnanoCasulli2,BrugnanoSestini}  and Casulli and Zanolli \cite{CasulliZanolli2010,CasulliZanolli2012}, and for which convergence has been rigorously proven. For finite elements without mass lumping, i.e., for non-diagonal mass matrices, the theorems which are the basis of the convergence proofs in the aforementioned works of Casulli \textit{et al.} do \textit{not} directly apply and still need to be generalized to more general non-diagonal but symmetric positive definite mass matrices. In the rest of this paper, we therefore assume that the simple ideal gas equation of state holds.

Spatial discretization is done by choosing a numerical method adapted to the nature of each equation: finite volumes are applied to approximate the transport-diffusion equations, whereas the Poisson problem is solved using continuous finite elements. The use of staggered grids avoids  the checker-board phenomena, which are typical for many numerical methods on collocated grids. The use of unstructured meshes increases the applicability of the methodology with respect to Cartesian grids since the meshing of complex domains becomes straightforward. The overall algorithm can be divided into four main stages:
\begin{itemize}
	\item Transport-diffusion stage. The equations \eqref{eq:mass_td_p}, \eqref{eq:momentumtt_dt_p}, and \eqref{eq:totenerg_td_p} are solved using explicit finite volumes in the dual mesh. To attain second order in space and time, a local ADER method combined with an ENO reconstruction is considered. Within this stage we get the new density, $\rho^{n+1}$, and the intermediate approximations of the momentum, $\ttWvel$, as well as the total energy density, $\tWE$ at each cell of the dual mesh.
	
	\item Pre-projection stage. The intermediate states for the total energy density and for the linear momentum are transferred from the dual to the primal grid. Next, the auxiliary variables that will be needed within the next stage, as the enthalpy, $H^{n+1,k}$, and the kinetic energy density, $\left(\rho K\right)^{n+1,k}$,
	are also computed. Let us note that the intermediate variables are calculated only once per time step, meanwhile the auxiliary variables are updated at each Picard iteration.
	
	\item Projection stage. A $\mathbb{P}^{1}$ finite element scheme is employed in order to solve the pressure equation \eqref{eq:dpressure_dt_p} implicitly. The resulting $\delta \Press^{n+1,k+1}$ is computed on the vertexes of the primal simplex mesh.
	
	\item Post-projection stage. The pressure correction at the new time $t^{n+1}$ is substituted in \eqref{eq:momentum_dt_p} and \eqref{eq:pressure_dt_p} to update the linear momentum $\Wvel^{n+1,k+1}$ and the pressure $\Press^{n+1,k+1}$. Once the Picard iterations have finished, the total energy, $\WE^{n+1}$, is recovered following \eqref{eq:energy_dt_p}. 
\end{itemize}
In what follows, we will further detail each stage of the algorithm.

\subsection{Staggered unstructured mesh}
In this paper we make use of two overlapping unstructured staggered meshes to discretize the domain $\Omega$. For the sake of simplicity, we focus here on the 2D case introducing the main notation needed. Further details on the construction of three-dimensional face-type staggered meshes can be found in \cite{BFSV14,TD16,BFTVC17,BTBD20,BBDFSVC20}.

\begin{figure}[h]
	\begin{center}
	\includegraphics[width=0.3\linewidth]{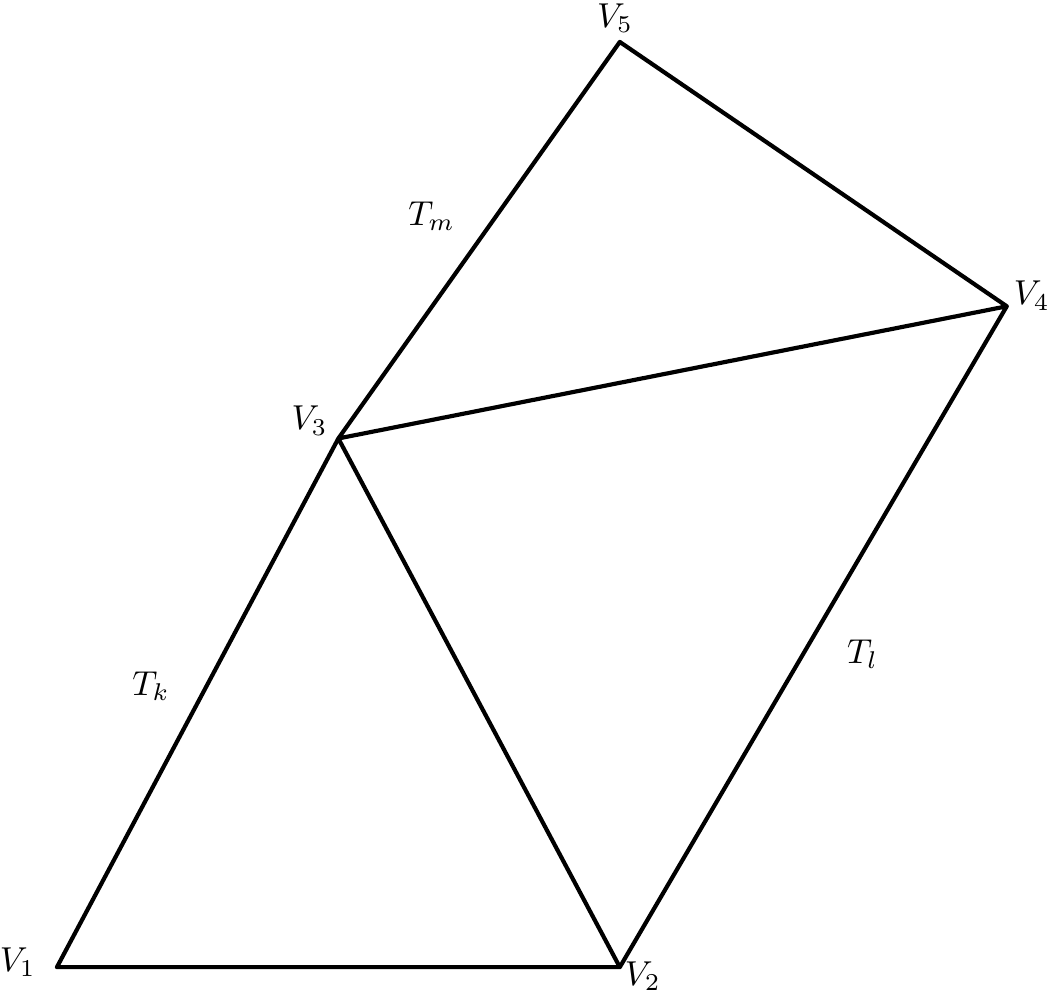}
	\includegraphics[width=0.3\linewidth]{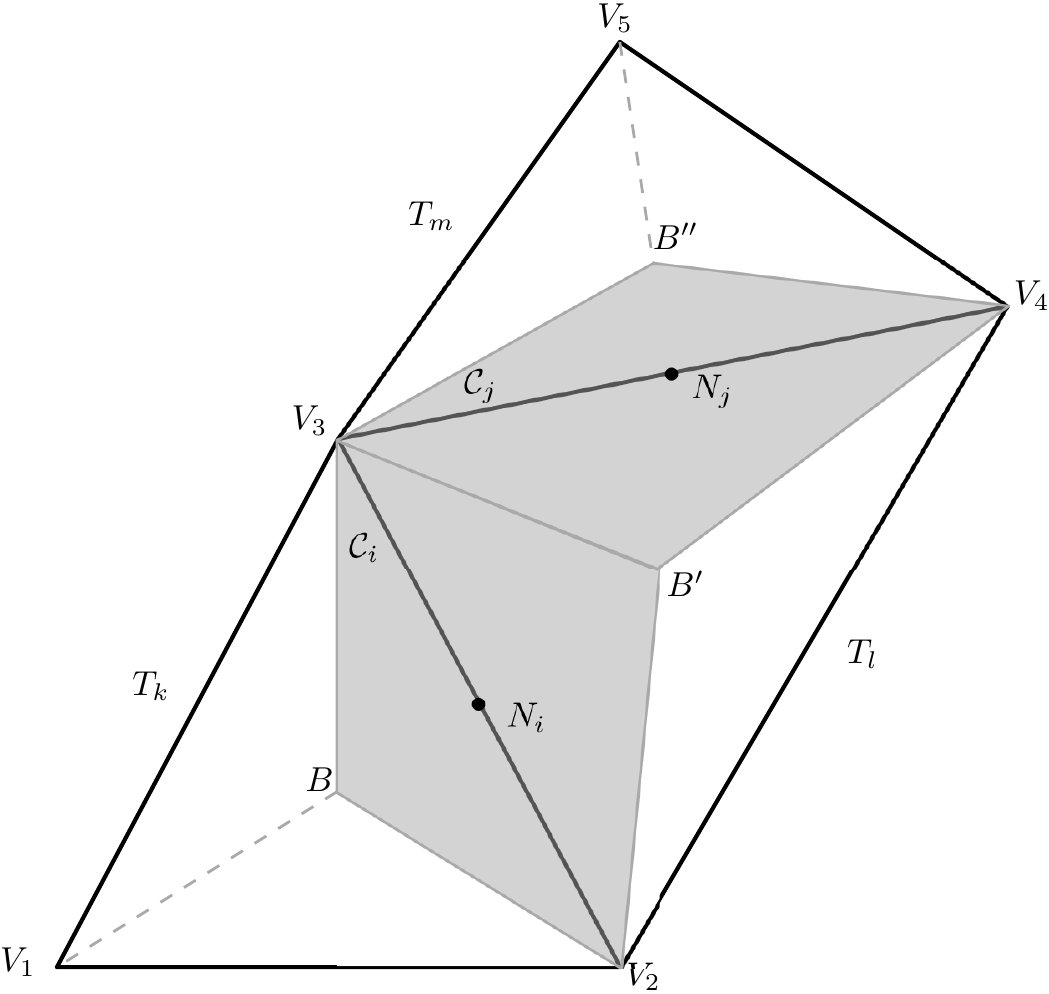}
	\includegraphics[width=0.3\linewidth]{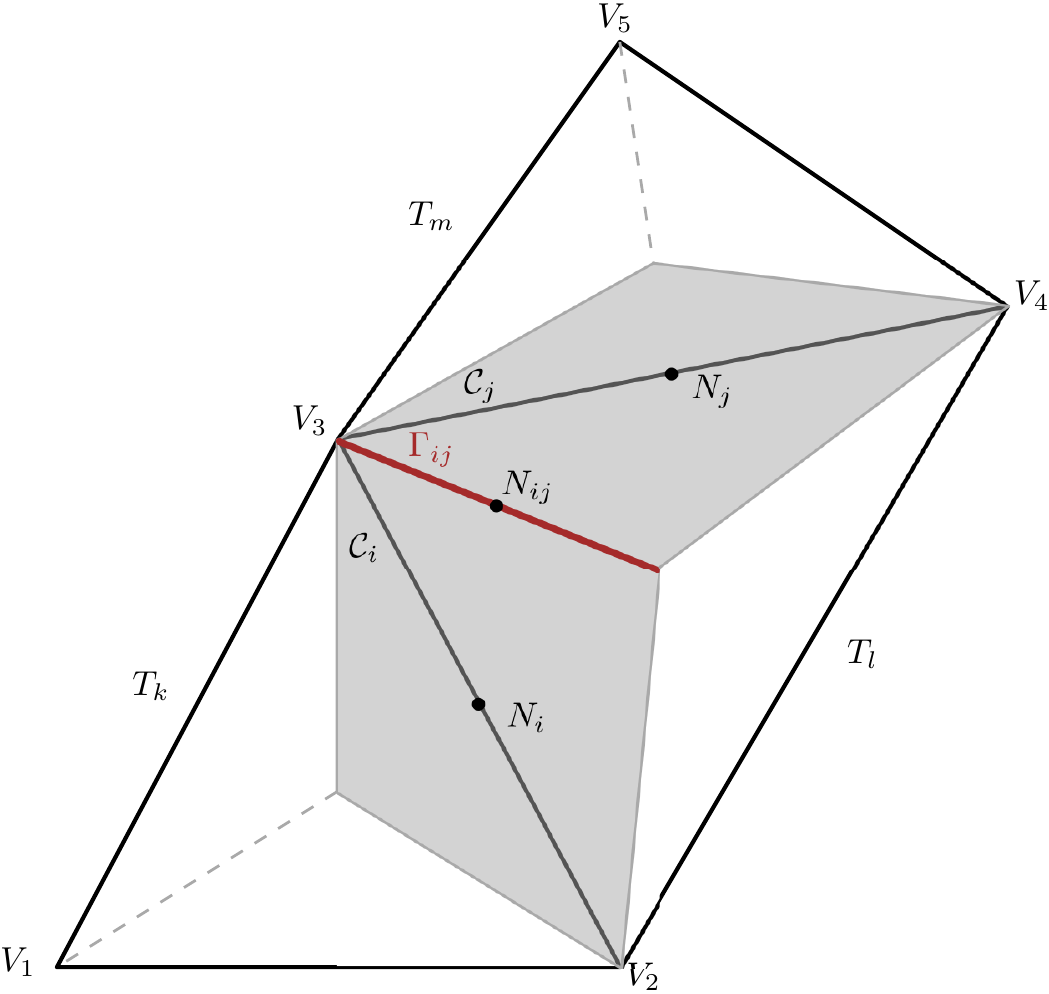}
	\end{center}
	\caption{Construction of face-type dual elements from a 2D primal mesh. Left: primal mesh made of elements $T_{k}$, $T_{l}$, $T_{m}$ and vertex $V_{n},\, n=1,\dots, 5$. Center: dual interior cells $C_{i}$, $C_{j}$ (shadowed in grey), white triangles correspond to boundary cells. Right: boundary face, $\Gamma_{ij}$, between the dual elements $C_{i}$, $C_{j}$.}
	\label{fig:2Dmesh}
\end{figure}

\begin{figure}[h]
	\begin{minipage}{0.6\linewidth}
		\begin{figure}[H]
			\begin{center}
				\begin{picture}(100,90)
				\put(0,0){\makebox(100,90){
				\vspace{-2cm}
				\includegraphics[width=6.2cm]{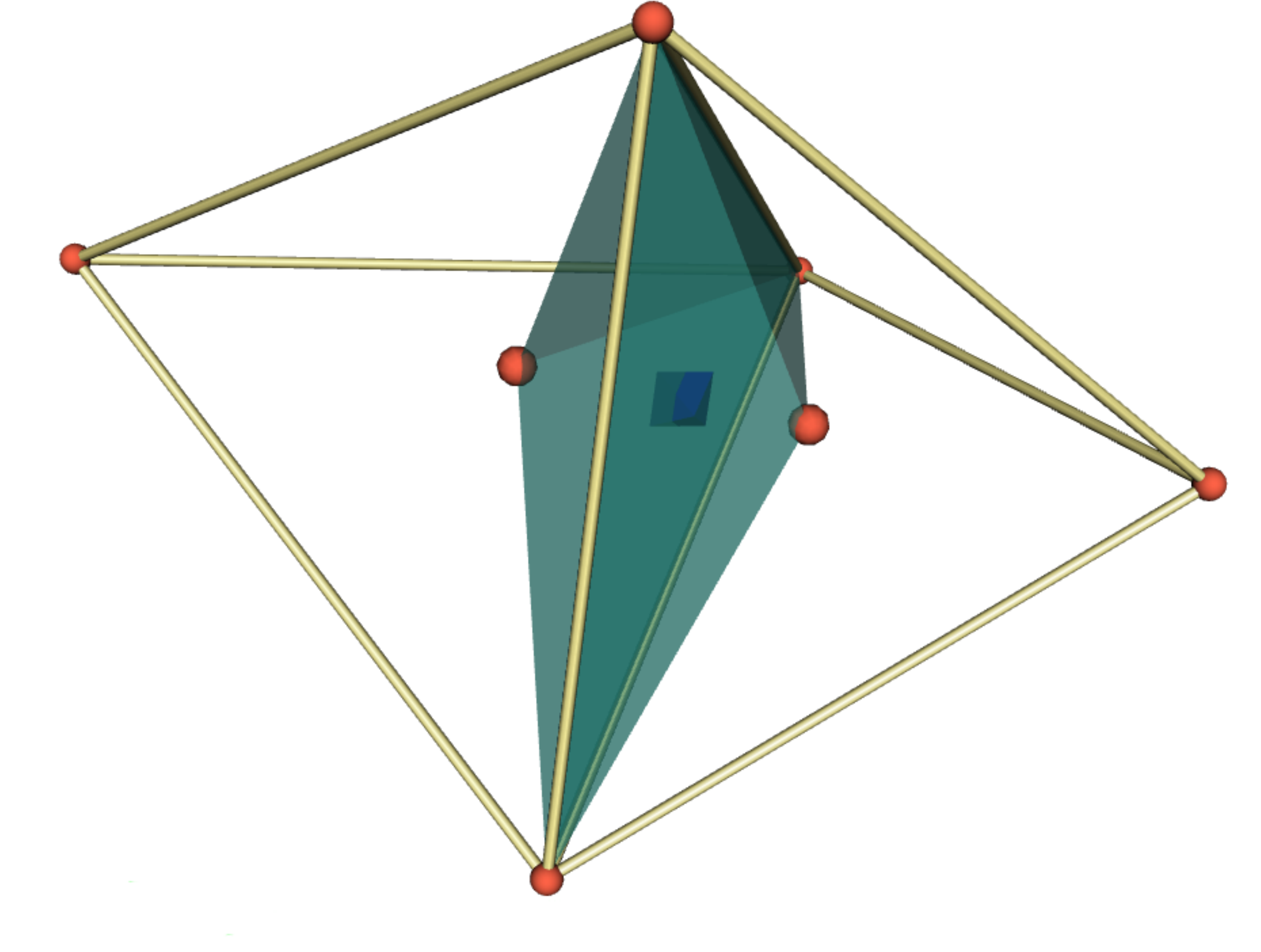}}}
				\put(33,-55){\normalsize{$V_{2}$}}
				\put(75,48){\normalsize{$V_{3}$}}
				\put(47,90){\normalsize{$V_{1}$}}
				\put(-40,52){\normalsize{$V_{4}$}}
				\put(140,12){\normalsize{$V_{4}'$}}
				\put(22,30){\normalsize{$B$}}
				\put(84,20){\normalsize{$B'$}}
				\put(56,33){\normalsize{$N_i$}}
				\end{picture}
				\vspace{1.5cm}
			\end{center}
		\end{figure}
	\end{minipage}
	\begin{minipage}{0.3\linewidth}
		\begin{figure}[H]
			\begin{center}
				\begin{picture}(100,90)
				\put(0,0){\makebox(100,90){
						\vspace{-2cm}
						\includegraphics[width=3.6cm]{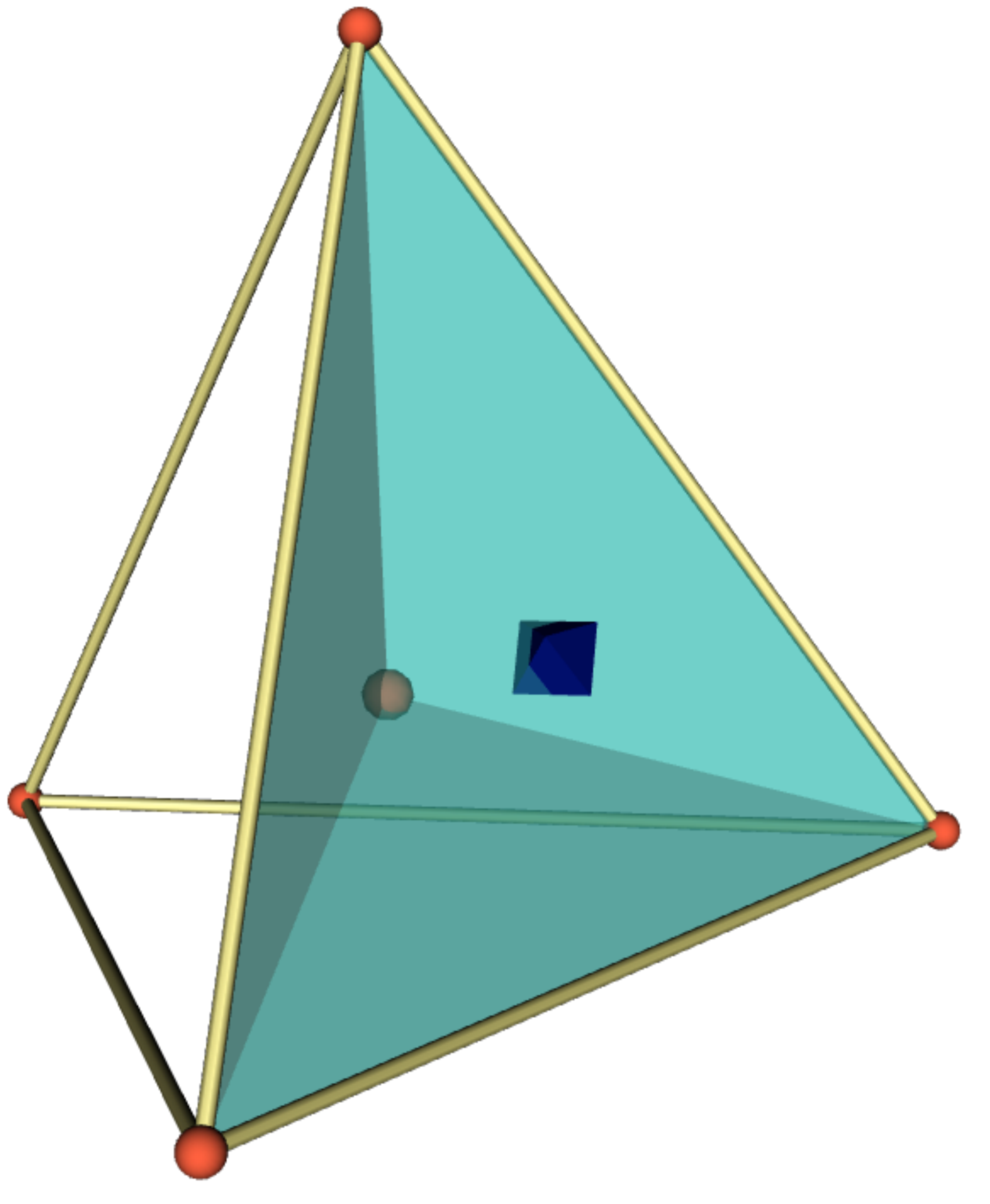}}}
				\put(12,-55){\normalsize{$V_{2}$}}
				\put(104,-13){\normalsize{$V_{3}$}}
				\put(30,87){\normalsize{$V_{1}$}}
				\put(-15,-12){\normalsize{$V_{4}$}}
				\put(30,0){\normalsize{$B$}}
				\put(58,15){\normalsize{$N_i$}}
				\end{picture}
				\vspace{1.5cm}
			\end{center}
		\end{figure}
	\end{minipage}
	\caption{Example of an interior finite volume (left) and of a boundary finite volume (right) which constitute the staggered dual mesh in three space dimensions. 
	}\label{fig:3Dmesh}
\end{figure}

Let us consider a triangular primal mesh $\{ T_k, \, k=1,\dots ,nel\}$ with vertex $\{V_{j}, \, j=1,\dots, nver\}$ (Figure~\ref{fig:2Dmesh} left). We now define the two triangles with basis one interior edge of a primal element and opposite vertex the barycenters, $B$, $B^{\prime}$, of the two primal elements sharing this face. The dual element, $C_{i}$, is then built by merging these two triangles (Figure~\ref{fig:2Dmesh} center). Similarly, a dual boundary element is a triangle which has as basis a primal boundary edge and as opposite vertex the barycenter of the primal element. Let us note that on the dual mesh the nodes $\left\lbrace N_i,\, i=1,\dots, nnod\right\rbrace$ are associated with the edges/faces of the simplex elements on the primal mesh. The remaining notation related to the mesh is as follows:
\begin{itemize}
	\item ${\cal K}_{i}$ is the set of neighbouring nodes of a node $N_i$ 
	consisting of the barycenters of the dual cells sharing a face with $C_{i}$.
	\item $\Gamma_{i}$ is the boundary of a cell $C_{i}$ and $\boldsymbol{\widetilde{\eta}}_{i}$ its outward unit normal.
	\item $\Gamma_{ij}$ is the edge between cells $C_{i}$ and $C_{j}$. $N_{ij}$ is the barycenter of $\Gamma_{ij}$ (Figure \ref{fig:2Dmesh} right). Note that 
	$\Gamma_i =\displaystyle \bigcup_{N_j \in {\cal K}_{i}} \Gamma_{ij}$.
	\item  $\left|C_{i}\right|$ is the area of $C_{i}$.
	\item $\widetilde{\boldsymbol{\eta}}_{ij} $ is the outward unit normal vector to $\Gamma_{ij}$. 	We define $\boldsymbol{\eta_{ij}}:=\boldsymbol{\widetilde{\eta}_{ij}} ||\boldsymbol{\eta_{ij}} ||$, where,
	$||\boldsymbol{\eta_{ij}} ||$  represents the length of $\Gamma_{ij}$. 
\end{itemize}
Sketches of the 2D and 3D staggered meshes are depicted in Figures \ref{fig:2Dmesh} and \ref{fig:3Dmesh}, respectively.

\subsection{Transport-diffusion stage}\label{sec:transport}
To solve the transport dominated equations, we use a finite volume method on the dual grid. As result, we will obtain the value of the averaged new density, $\Wrho^{n+1}=\rho^{n+1}$, on each dual cell, as well as intermediate approximations for the cell averaged linear momentum, $\ttWvel$, and total energy density, $\tWE$. 
We start integrating equations \eqref{eq:mass_td_p}-\eqref{eq:totenerg_td_p} on each dual cell $C_{i}$ and applying Gauss theorem which yields
\begin{gather}
\rho^{n+1}_{i}= \rho^{n}_{i} - \frac{\Delta t}{\left|C_i\right|} \int_{\Gamma_i} \Flux^{\Wrho}\left(\mathbf{W}^n\right)\boldsymbol{\widetilde{\eta}}_{i}\,dS, \label{eq:discr_mass}
\\
 \widetilde{\widetilde{\mathbf{W}}}_{\mathbf{u},\, i} = \mathbf{W}_{\mathbf{u},\, i}^{n} -\frac{\Delta t}{\left|C_i\right|} \left( \int_{\Gamma_i}\Flux^{\mathbf{W}_{\mathbf{u}}}\left(\mathbf{W}^n\right)\boldsymbol{\widetilde{\eta}}_{i}\,dS + \int_{C_i} \grae \Press^{n} dV -\int_{\Gamma_i} \bTau^n \,\boldsymbol{\widetilde{\eta}}_{i}\, dS- \int_{C_i} \rho^{n}\g dV\right) ,\label{eq:discr_momentum}
\\
\tWEi = \WEi^{n} -\frac{\Delta t}{\left|C_i\right|}\left(  \int_{\Gamma_i} \Flux^{\WE}\left(\mathbf{W}^n\right) \cdot \boldsymbol{\widetilde{\eta}}_{i}\,dS  - \int_{\Gamma_i} \left(\frac{1}{\rho^{n}}\bTau^{n}\Wvel^{n}\right)\cdot\boldsymbol{\widetilde{\eta}}_{i}\,dS +\int_{\Gamma_i} \Q^{n}\cdot\boldsymbol{\widetilde{\eta}}_{i}\,dS - \int_{C_i} \g \cdot \Wvel^{n} dV\right).
\label{eq:discr_energy}
\end{gather}
where 
\begin{equation}
\Flux^{\Wrho} (\WW^{n}) := \Wvel^{n}, \qquad \Flux^{\Wvel} (\WW^{n}) :=\frac{1}{\rho^{n}}\Wvel^{n}\otimes\Wvel^{n} ,\qquad \Flux^{\WE} (\WW^{n}) := \Ken^{n} \Wvel^{n}
\end{equation}
are the convective fluxes of mass, momentum, and total energy, respectively.

\subsubsection{Convective numerical fluxes}\label{sec:advectionterm}
The global normal flux through the boundary of a dual cell is denoted by
\begin{equation}
{\mathcal{Z} }(\WW^{n},  {\boldsymbol{\tilde{\eta}}}_{i} ): ={
	\Flux} (\WW^{n})  {\boldsymbol{\tilde{\eta}}}_{i},\quad\mathrm{ with }\quad
\Flux (\WW^{n})= \left(\Flux^{\Wrho} (\mathbf{W}^{n}) , \Flux^{\Wvel} (\WW^{n}), \Flux^{\WE} (\WW^{n}) \right)^{T}.
\end{equation}
Let us recall that the flux contribution in the energy equation accounts only for the kinetic energy density contribution, $\rhoken$, instead of the total energy density, $\rho E$. Within the flux computation the value of $K^{n}$ is recovered from the linear momentum and density fields, $K^{n}=\frac{1}{2\rho^{n}} \left|\Wvel^{n}\right|^{2}$.
The integral of the flux term on $\Gamma_{i}$ can be split into the sum of the integral on the cell faces, $\Gamma_{ij}$,
\begin{equation}
\int_{ \Gamma_{i} } \mathbf{\mathcal{F}}(
\mathbf{W}^{n})  \boldsymbol{ \widetilde{\eta}}_{i} \,\mathrm{dS}
=\displaystyle \sum_{N_j \in {\cal K}_{i}} \displaystyle
\int_{\Gamma_{ij}} {
	\mathcal{Z} }(\mathbf{W}^{n}, {\boldsymbol{\tilde{\eta}}}_{ij} )  \,\mathrm{dS},
\end{equation}
and approximated using an upwind scheme to get a stable discretization. In particular, we consider a modified Rusanov flux function, \cite{Rus62,BBDFSVC20},
\begin{gather}
\boldsymbol{\phi} \left( \WW_{i}^{n}\, ,\WW_{j}^{n}\, , \boldsymbol{\eta}_{ij} \right)
=
\left( 
\phi_{\rho} \left( \WW_{i}^{n}\, ,\WW_{j}^{n}\, , \boldsymbol{\eta}_{ij} \right),
\boldsymbol{\phi}_{\vel} \left( \WW_{i}^{n}\, ,\WW_{j}^{n}\, , \boldsymbol{\eta}_{ij} \right),
\phi_{E}  \left( \WW_{i}^{n}\, ,\WW_{j}^{n}\, , \boldsymbol{\eta}_{ij} \right)
\right)^{T} \nonumber\\
= \displaystyle\frac{1}{2} ({\cal Z}
(\WW_{i}^n\, ,\boldsymbol{\eta}_{ij})+{\cal Z}(\WW_{j}^n \, , \boldsymbol{\eta}_{ij}))-\frac{1}{2}
\alpha^{n}_{RS,\, ij} \left( \Wur_{j}^n-\Wur_{i}^n\right) \label{eq:flux}
\end{gather}
with
\begin{equation}
	\Wur := \left(\Wrho, \Wvel, \rho K \right)^{T}
\end{equation}
the modified conservative variables vector,
\begin{equation}\alpha^{n}_{RS,\, ij}=\alpha_{RS}(\mathbf{W}_{i}^n,\mathbf{W}_{j}^n,\boldsymbol{\eta}_{ij}):=\max \left\lbrace \left|\mathbf{U}_i^{n}\cdot\boldsymbol{\eta}_{ij}\right|,\left|\mathbf{U}_j^{n}\cdot\boldsymbol{\eta}_{ij}\right|\right\rbrace + c_{\alpha} \left\|\boldsymbol{\eta}_{ij} \right\|  \label{eq:alphars_decoupled1_comp}\end{equation}
the maximum signal speed on the edge and $c_{\alpha}\in\mathbb{R}_0^+$ an artificial viscosity coefficient that may be activated on particular tests to increase the stability properties of the final scheme when large variations of the density and energy fields are encountered in the presence of small velocities. Substitution in \eqref{eq:discr_mass}-\eqref{eq:discr_energy} gives  
\begin{gather}
\rho^{n+1}_{i}= \rho^{n}_{i} - \frac{\Delta t}{\left|C_i\right|}  \sum_{N_j \in {\cal K}_{i}} \phi_{\rho} \left( \WW_{i}^{n}\, ,\WW_{j}^{n}\, , \boldsymbol{\eta}_{ij} \right), \label{eq:discrf_mass}
\\
\widetilde{\widetilde{\mathbf{W}}}_{\mathbf{u},\, i} = \mathbf{W}_{\mathbf{u},\, i}^{n} -\frac{\Delta t}{\left|C_i\right|} \left( \sum_{N_j \in {\cal K}_{i}} \boldsymbol{\phi}_{\vel} \left( \WW_{i}^{n}\, ,\WW_{j}^{n}\, , \boldsymbol{\eta}_{ij} \right) + \int_{C_i} \grae \Press^{n} dV -\int_{\Gamma_i} \bTau^n \,\boldsymbol{\widetilde{\eta}}_{i}\, dS- \int_{C_i} \rho^{n}\g dV\right) ,\label{eq:discrf_momentum}
\\
\tWEi = \WEi^{n} -\frac{\Delta t}{\left|C_i\right|}\left(  \sum_{N_j \in {\cal K}_{i}} \phi_{E} \left( \WW_{i}^{n}\, ,\WW_{j}^{n}\, , \boldsymbol{\eta}_{ij} \right)  - \int_{\Gamma_i} \left(\frac{1}{\rho^{n}}\bTau^{n}\Wvel^{n}\right)\cdot\boldsymbol{\widetilde{\eta}}_{i}\,dS +\int_{\Gamma_i} \Q^{n}\cdot\boldsymbol{\widetilde{\eta}}_{i}\,dS - \int_{C_i} \g \cdot \Wvel^{n} dV\right).
\label{eq:discrf_energy}
\end{gather}

The scheme proposed above would result in a first order scheme in space and time.
To increase the order of accuracy attaining second order in space and time, a local ADER methodology (LADER) is employed, see \cite{BFTVC17,Bus18,BBDFSVC20}. The reader is referred to \cite{TMN01,Toro,Mill99,BTVC16} for further details on the original ADER methodology and to \cite{DET08,BD14,BCDGP20} for an alternative variant of ADER schemes that allow to avoid the cumbersome Cauchy-Kovalevskaya procedure thanks to the use of a general space-time finite element predictor.  
In what follows, we briefly recall the main steps to be performed in the LADER algorithm:
\begin{description}
	\item[Step 1.] Piecewise polynomial reconstruction in the neighbourhood of each boundary edge of the cell. Considering an scalar conservative variable, $W$, and one of the cell boundaries, $\Gamma_{ij}$, the related reconstruction polynomials read	
	\begin{equation}P^{i}_{ij}(N)=W_i+(N-N_{i})\left( \grae\, W\right)^{i}_{ij},\quad P^{j}_{ij}(N)=W_j+(N-N_{j})\left( \grae\, W\right)^{j}_{ij}.\end{equation}
	To circumvent Godunov's theorem and to develop a second order scheme avoiding spurious oscillations, we introduce a non linearity via the use of a nonlinear Essentially Non-Oscillatory (ENO) reconstruction. Accordingly, the gradients are computed as 
	\begin{gather*}
	\left( \grae\, W\right)^{i}_{ij}= \left\lbrace
	\begin{array}{lc}
	\left(\grae\, W \right)_{T_{ijL}}, & \textrm{if }\left| \left(\grae\, W \right)_{ T_{ijL}}\cdot \left(N_{ij}-N_{i}\right)\right| \leq \left|\left(\grae\, W \right)_{ T_{ij}}\cdot \left(N_{ij}-N_{i}\right) \right|,\\[8pt]
	\left(\grae\, W \right)_{ T_{ij}}, & \mathrm{otherwise};
	\end{array}
	\right.
	\end{gather*}
	\begin{gather*}
	\left( \grae\, W \right)^{j}_{ij}= \left\lbrace
	\begin{array}{lc}
	\left(\grae\, W\right)_{ T_{ijR}}, & \textrm{if }\left|\left( \grae\, W \right)_{ T_{ijR}}\cdot \left(N_{ij}-N_{j}\right)\right| \leq \left|\left(\grae\, W \right)_{ T_{ij}}\cdot \left(N_{ij}-N_{j}\right)\right|,\\[8pt]
	\left(\grae\, W\right)_{ T_{ij}},  & \mathrm{otherwise},
	\end{array}
	\right.
	\end{gather*}
	with $T_{ij}$, $T_{ijL}$ and $T_{ijR}$ the centered and upwind primal elements to the face $\Gamma_{ij}$ where the gradients are computed using a Galerkin approach (Crouzeix-Raviart finite elements). An alternative to the ENO-based reconstruction is the use of  classical slope limiters like the Barth and Jespersen limiter \cite{BJ89}, or the minmod limiter of Roe \cite{Roe85,Toro}. Also, \textit{a posteriori} limiting strategies like the MOOD approach, \cite{CDL11}, could be used and will be part of future research.
	
	\item[Step 2.] Calculation of the necessary boundary--extrapolated data in $\mathbf{x}_{N_{ij}}$ of each edge/face of the FV mesh,	
	\begin{gather} W_{i\, N_{ij}} =p^{i}_{ij}(N_{ij})
	= W_{i}+(N_{ij}-N_{i})   \left( \grae\, W\right)^{i}_{ij},\\
	W_{j\, N_{ij}}=p^{j}_{ij}(N_{ij})=
	W_{j}+(N_{ij}-N_{j})   \left( \grae\, W\right)^{j}_{ij}.  \end{gather}
	
	\item[Step 3.] Use of the mid-point rule to get a second order of accuracy approximation in time. A temporal Taylor series expansion in combination with the Cauchy-Kovalevskaya procedure, based on the mass, momentum, and energy equations \eqref{eq:mass2}, \eqref{eq:momentum2}, \eqref{eq:dpressure_dt_p}, are
	employed in order to approximate the conservative variables at the time  $t^n + \frac{\Delta t}{2}$:
	\begin{align}\overline{\WW}_{i\, N_{ij}}=\rho_{i\, N_{ij}} + {\WW}_{ N_{ij}}^{\ast},\qquad \overline{\WW}_{j\, N_{ij}}= \rho_{j\, N_{ij}} + {\WW}_{N_{ij}}^{\ast}\label{eq:Wevolij}\end{align}
	where	
	\begin{eqnarray}{\WW}_{\rho\, N_{ij}}^{\ast}&:=&-\frac{\Delta t}{2\mathcal{L}_{ij}}\left[ \mathcal{Z}^{\rho}(\mathbf{W}_{i\, N_{ij}},\boldsymbol{\eta}_{ij}  )+\mathcal{Z}^{\rho}(\mathbf{W}_{j\, N_{ij}},\boldsymbol{\eta}_{ij}  )\right],\label{eq:densityevol}\\
	{\WW}_{\mathbf{u}\, N_{ij}}^{\ast}&:= &-\frac{\Delta t}{2\mathcal{L}_{ij}}\left[ \mathcal{Z}^{\vel}(\mathbf{W}_{i\, N_{ij}},\boldsymbol{\eta}_{ij}  )+\mathcal{Z}^{\vel}(\mathbf{W}_{j\, N_{ij}},\boldsymbol{\eta}_{ij}  )\right]  +
	\frac{\Delta t}{2\mathcal{L}_{ij}} \bTau^{\ast}\boldsymbol{ \eta }_{ij} 
	\notag\\&&
	-\frac{ \Delta t}{2}  \left(\grae\, \Press\right)_{T_{ij}} +\frac{ \Delta t}{4}   \left(\rho_{i\, N_{ij}} + \rho_{j\, N_{ij}}\right) \g, \\
	{W}_{E\, N_{ij}}^{\ast}&:= &-\frac{\Delta t}{2\mathcal{L}_{ij}}\left[ \mathcal{Z}^{E}(\mathbf{W}_{i\, N_{ij}},\boldsymbol{\eta}_{ij}  )+\mathcal{Z}^{E}(\mathbf{W}_{j\, N_{ij}},\boldsymbol{\eta}_{ij}  )\right]
	+\frac{\Delta t}{2\mathcal{L}_{ij}} \left(\bTau\Vel\right)^{\ast} \cdot\boldsymbol{\eta}_{ij}
	 \notag\\&&	
	+\frac{ \Delta t}{4} \g \cdot \left(\WW_{\vel\, i\, N_{ij}} + \WW_{\vel, j\, N_{ij}}\right), 
	\end{eqnarray}
	with
	\begin{gather}
	\bTau^{\ast}  = \mu \left[ \left( \gra\, \mathbf{W}_{\mathbf{u}}+\gra\, \mathbf{W}_{\mathbf{u}}^{T}\right)_{i\, N_{ij}}  +  \left( \gra\, \mathbf{W}_{\mathbf{u}}+\gra\, \mathbf{W}_{\mathbf{u}}^{T}\right)_{j\,N_{ij}} 
	 -\frac{2}{3} \left( \dive \mathbf{W}_{\mathbf{u},\, i \, N_{ij}} \mathbf{I} +  \dive  \mathbf{W}_{\mathbf{u},\, j \, N_{ij}} \mathbf{I}\right) \right],\\
	\left(\bTau\Vel\right)^{\ast}  = \mu \left[ \left( \gra\, \mathbf{W}_{\mathbf{u}}+\gra\, \mathbf{W}_{\mathbf{u}}^{T}\right)_{i\, N_{ij}} \Vel_{i\, N_{ij}} +  \left( \gra\, \mathbf{W}_{\mathbf{u}}+\gra\, \mathbf{W}_{\mathbf{u}}^{T}\right)_{j\,N_{ij}} \Vel_{j\, N_{ij}} \phantom{\frac{2}{3}\!\!\!}\right.\notag\\ \left.
	 -\frac{2}{3} \left( \dive \mathbf{W}_{\mathbf{u},\, i \, N_{ij}} \Vel_{i\, N_{ij}}+  \dive  \mathbf{W}_{\mathbf{u},\, j \, N_{ij}} \Vel_{j\, N_{ij}}\right) \right].
	\end{gather}
		
	\item[Step 4.]  Calculation of the numerical flux for the convective terms using \eqref{eq:flux},
	\begin{gather}
	\boldsymbol{\phi} \left(\overline{\mathbf{W}}_{i\, N_{ij}}^{n},\overline{\mathbf{W}}_{j\, N_{ij}}^{n},\boldsymbol{\eta}_{ij}\right)
	= \displaystyle\frac{1}{2} \left[ {\cal Z}
	(\overline{\mathbf{W}}_{i\, N_{ij}}^n,\boldsymbol{\eta}_{ij})+{\cal Z}(\overline{\mathbf{W}}_{j\, N_{ij}}^n,\boldsymbol{\eta}_{ij})\right]
	-\frac{1}{2}
	\overline{\alpha}^{n}_{RS,\, ij} \left(\overline{\Wur}_{j\, N_{ij}}^n -\overline{\Wur}_{ i \, N_{ij}}^n\right).\label{eq:flux_lmnt_lader}
	\end{gather}
\end{description}

\subsubsection{Viscous term}\label{sec:viscousterm}
Gauss' theorem allows to rewrite the volume integrals of the viscous stress tensor and of the heat flux into 
surface integrals over the boundary, $\Gamma_{i}$,  which can further be rewritten as a sum of integrals over the 
individual cell faces, $\Gamma_{ij}$. This leads to 
\begin{gather}
\int_{C_i}\dive \bTau^n dV = \sum_{N_j\in\mathcal{K}_i}\int_{\Gamma_{ij}} \bTau^n\boldsymbol{\widetilde{\eta}}_{ij}  \, \mathrm{dS} 
= \sum_{N_j\in\mathcal{K}_i}\int_{\Gamma_{ij}}\mu\left[
\gra\, \Vel^{n}+\left( \gra\,\Vel^{n}\right)^T  - \frac{2}{3} \dive  \Vel^{n} I\right]\boldsymbol{\widetilde{\eta}}_{ij}  \,\mathrm{dS},\label{eq:visc_mom}\\
\int_{C_i}\dive \left( \frac{1}{\rho^{n}}\bTau^n\Wvel^{n}\right)  dV = \sum_{N_j\in\mathcal{K}_i}\int_{\Gamma_{ij}} \frac{1}{\rho^{n}}\bTau^n\Wvel^{n}\cdot\boldsymbol{\widetilde{\eta}}_{ij}  \, \mathrm{dS} 
= \sum_{N_j\in\mathcal{K}_i}\int_{\Gamma_{ij}}\mu\left[\left( 
\gra\, \Vel^{n}+\left( \gra\,\Vel^{n}\right)^T  - \frac{2}{3} \dive  \Vel^{n} I\right) \Vel^{n}\right] \cdot\boldsymbol{\widetilde{\eta}}_{ij} \, \mathrm{dS},\label{eq:visc_toten}
\end{gather}
where the gradients are computed on the primal element containing the face, $T_{ij}$, following the Galerkin approach already introduced in Section \ref{sec:advectionterm} within the LADER reconstruction.
The corresponding numerical diffusion functions read
\begin{gather}
\varphi_{\vel}\left( \Vel_{i}^n,\Vel_{j}^n, \boldsymbol{\eta}_{ij} \right) \approx \int_{\Gamma_{ij}}   \mu \left[ \gra\, \Vel^n +\left( \gra\, \Vel^n\right)^{T} - \frac{2}{3} \dive \Vel^n \right] {\boldsymbol{\widetilde{\eta}}}_{ij}\mathrm{dS}, \label{eq:numviscfluxfunU}\\
\varphi_{E}\left( \mathbf{U}_{i}^n,\mathbf{U}_{j}^n, \boldsymbol{\eta}_{ij} \right) \approx \int_{\Gamma_{ij}} \mu\left[\left( 
\gra\, \Vel^{n}+\left( \gra\,\Vel^{n}\right)^T  - \frac{2}{3} \dive  \Vel^{n} I\right) \Vel^{n}\right] \cdot\boldsymbol{\widetilde{\eta}}_{ij}\mathrm{dS} . \label{eq:numviscfluxfunE}
\end{gather}
Similar to what has been done for the advection term, a Taylor series expansion combined with the Cauchy-Kovalevskaya procedure could be applied to get a second order accurate approximation of the viscous terms in space and time. The main difference with respect to the flux terms computation is that we now can neglect the presence of the flux term on the reconstruction of the linear momentum field. As it has been shown in \cite{BFTVC17} for the scalar advection-diffusion-reaction equation, the specific way of computing the gradients makes the mixed contribution of advection and diffusion terms be completely included in the time evolution of the flux to the half time level. We should notice that the evolution of viscous terms may lead to a more restrictive CFL stability condition, so the time step would be smaller than when applying LADER only to convective terms, see \cite{BTVC16}. Once the evolution to the half time level of the linear momentum, $\overline{\overline{\WW}}_{\vel,\,i}^{n}$, is computed, it is divided by the reconstructed density to approximate the evolved velocity, $\overline{\overline{ \Vel}}_{i}^{n}=\left( \overline{\rho}_{i}^{n} \right)^{-1}\overline{\overline{ \WW}}_{\vel,\,i}^{n}$, to be inserted in \eqref{eq:numviscfluxfunU}-\eqref{eq:numviscfluxfunE}.

\subsubsection{Pressure term}
To account for the pressure contribution at the previous time step, we transform the integral of its gradient on the dual cell into the sum of the normal projection on each face
\begin{equation}
\int_{C_i} \grae \Press^{n} dV = \sum_{N_j \in {\cal K}_{i}} \Press^{n}_{ij}\boldsymbol{ {\eta}_{ij}}.\label{eq:pressterm}
\end{equation} 
The value of the pressure at each dual face, $\Press^{n}_{ij}$, is obtained as the average of the pressure at its vertexes. Regarding the vertex corresponding to the barycenter of the primal element, the pressure is approximated again as the averaged value on the nodes of the primal element. To get second order in space and time, the pressure in \eqref{eq:pressterm} is replaced by its half in time reconstructed value. The LADER methodology is applied like for the viscous term computation to get also the value of $\overline{\overline{\WW}}_{E,\,i}^{n}$ from which the evolved pressure $\overline{\overline{\Press}}^{n}_{ij}$ can be recovered using relation \eqref{eq:E.as.pk}.

\subsubsection{Gravity term}
The gravity source terms in \eqref{eq:momentumtt_dt_p}-\eqref{eq:totenerg_td_p} are integrated on each spatial control volume $C_{i}$ by taking the averaged density and linear momentum fields on the cell:
\begin{equation}
\int_{C_i} \rho^{n}\g  dV = \left|C_i\right| \rho^{n}_{i} \g,\quad
\int_{C_i} \g \cdot \Wvel^{n} dV = \left|C_i\right| \g\cdot\WW^{n}_{\vel\,i}.
\end{equation}

\subsubsection{Heat flux term}
Let us assume that the averaged cell temperature $\Theta_{i}^{n}$ is known. Then, it can be used to approximate the temperature gradients on each primal cell as already done for the gradients of conservative variables in the flux and viscous terms. Finally, these values are employed to approximate the integral of the heat flux term after applying Gauss theorem, 
\begin{equation}
\int_{C_i} \mathrm{div} \mathbf{Q}^{n} dV = 
\sum_{N_j \in {\cal K}_{i}}\mathbf{Q}_{T_{ij}}^{n}  \cdot {\boldsymbol{{\eta}}}_{ij}=
-\sum_{N_j \in {\cal K}_{i}}\lambda \left(\gra\, \Theta^{n} \right)_{T_{ij}} \cdot {\boldsymbol{{\eta}}}_{ij}
\end{equation}
From the implementation point of view, the temperature, $\Theta_{i}^{n}$ can be computed at the previous time step using the averaged values of the pressure and density at the dual cells: 
\begin{equation}
\Theta_{i}^{n} = \frac{\Press_{i}^{n}}{\rho_{i}^{n}R}. 
\end{equation}

\subsection{Pre-projection stage}\label{sec:prepro}
Some of the terms in the pressure system require for the preprocessing of the involved variables since they need to be transferred from the dual mesh to the primal one, or they do not belong to the original unknowns of the system to be solved. 

Given a scalar variable at the dual cells, $W_{i}$,  its value on the primal element $T_{k}$ is computed as the weighted average of the values on the subelements of a primal element associated to each face, $T_{ki}$, $i\in \mathcal{K}_{k}$, $\mathcal{K}_{k}$ the set of indexes identifying the faces of primal elements,
\begin{equation}
W_{k} = \sum_{i\in \mathcal{K}_{k}} W_{i} \frac{\left|T_{ki}\right|}{\left| T_{k}\right|}. \label{eq:interpolation_dual2primal}
\end{equation}
This approach is used for the computation of the density, $\Wrho^{n+1}$, the intermediate velocity, $\ttWvel$, and the intermediate total energy density, $\tWE$, by primal element.
Then, the first guess for the kinetic energy density,
\begin{equation}
\left(\rho K\right)^{n+1,1}= \frac{1}{2\rho^{n+1}}\left|\Wvel^{n+1,1}\right|^2,
\end{equation}
is obtained. On the other hand, the first guess for the enthalpy is initially computed at the dual mesh,
\begin{equation}
H^{n+1,1}_{i} = \frac{\gamma}{\gamma-1}\frac{P^{n}_{i}}{\rho^{n+1}_{i}}
\end{equation}
and used when its value on the faces of the primal elements is needed. 
Passing to a value per primal element is again done following \eqref{eq:interpolation_dual2primal}.

\subsection{Projection stage}\label{sec:projection}
For the projection stage a classical $\mathbb{P}_1$ continuous finite element method is employed to approximate the pressure correction $\delta\Press^{n+1,k+1} $ in the primal grid nodes in each Picard iteration. To get the weak formulation of the Laplacian problem we  start by multiplying equation \eqref{eq:dpressure_dt_p} by a test function $z\in V_{0}$, 
$V_{0}:=\left\lbrace z\in H^{1}\left(\Omega\right): \int_{\Omega} z dV = 0\right\rbrace$, integrating in $\Omega$ and applying Green's formula:
\begin{gather}
\frac{1}{\gamma-1} \int_{\Omega} \delta \Press^{n+1,k+1} z \dV
- \Delta t^{2} \int_{\Omega} H^{n+1,k} \grae \delta\Press^{n+1,k+1}\cdot \grae z dV = \nonumber\\  \Delta t^{2} \int_{\Gamma} H^{n+1,k} \grae \delta\Press^{n+1,k+1}\cdot \boldsymbol{\eta} z dA
+ \int_{\Omega}\left( \tWE -\left(\RhoKen\right)^{n+1,k} \right) z dV
\nonumber\\ 
- \frac{1}{\gamma-1} \int_{\Omega} \Press^{n} z dV
+ \Delta t \int_{\Omega}   H^{n+1,k} \ttWvel\cdot \grae z dV
- \Delta t \int_{\Gamma}   H^{n+1,k} \ttWvel\cdot \boldsymbol{\eta} z dA, \qquad\label{eq:dpressure_dt_p_ef}
\end{gather}
Next, taking into account \eqref{eq:ttwvel}, we have
\begin{equation}
\ttWvel = \tWvel -\Delta t \grae \Press^{n} = \Wvel^{n+1} + \Delta t \grae \Press^{n+1} = \Wvel^{n+1} + \Delta t \grae \delta \Press^{n+1}.
\end{equation}
Multiplication by the time step, $\Delta t$, the enthalpy, $H^{n+1,k}$, the normal vector, $\boldsymbol{\eta}$ and the test function, $z$, and integration in the boundary of the computational domain, $\Gamma$, gives
\begin{equation}
\Delta t\int_{\Gamma} H^{n+1,k}\ttWvel \cdot \boldsymbol{\eta} z dA = \Delta t\int_{\Gamma}H^{n+1,k} \Wvel^{n+1} \cdot \boldsymbol{\eta} z dA  + \Delta t^{2}\int_{\Gamma}H^{n+1,k}\grae \delta \Press^{n+1} \cdot \boldsymbol{\eta}  z dA .
\end{equation}
Rearranging the above equation, we get
\begin{equation}
\Delta t\int_{\Gamma} H^{n+1,k}\ttWvel \cdot \boldsymbol{\eta} z dA - \Delta t^{2}\int_{\Gamma}H^{n+1,k}\grae \delta \Press^{n+1} \cdot \boldsymbol{\eta}  z dA   = \Delta t\int_{\Gamma} H^{n+1,k} \Wvel^{n+1} \cdot \boldsymbol{\eta} z dA .\label{eq:hdpetaz}
\end{equation}
Substitution of \eqref{eq:hdpetaz} into \eqref{eq:dpressure_dt_p_ef} leads to the following weak problem:
\begin{weakproblem}
	In each Picard iteration $k$ find the pressure correction $\delta\Press^{n+1,k+1} \in V_{0}$ that satisfies 
	\begin{gather}
	\frac{1}{\gamma-1} \int_{\Omega} \delta \Press^{n+1,k+1} z \dV
	- \Delta t^{2} \int_{\Omega} H^{n+1,k} \grae \delta\Press^{n+1,k+1}\cdot \grae z \, dV = \nonumber\\  - \Delta t\int_{\Gamma} H^{n+1,k} \Wvel^{n+1} \cdot \boldsymbol{\eta} z dA 
	+ \int_{\Omega}\left( \tWE -\left(\RhoKen\right)^{n+1,k} \right) z dV
	- \frac{1}{\gamma-1} \int_{\Omega} \Press^{n} z dV
	+ \Delta t \int_{\Omega}   H^{n+1,k} \ttWvel\cdot \grae z dV,
	\label{eq:weak_problem}
	\end{gather} 
	for all $z\in V_{0}$.
\end{weakproblem}
The presence of the enthalpy and the kinetic energy on \eqref{eq:weak_problem} would make the above system highly nonlinear if $H^{n+1}$ and $\left(\RhoKen\right)^{n+1}$ would have been employed instead of $H^{n+1,k}$ and $\left(\RhoKen\right)^{n+1,k}$. To avoid direct resolution of such a complex system a classical approach consists in employing a Picard iteration procedure. Following the ideas introduced in \cite{CZ09,DBTM08,DC16,TD15,TD17} to circumvent the non-linearities arising in semi-implicit and locally implicit schemes for nonlinear PDEs, the enthalpy and the kinetic energy term on the right hand side of \eqref{eq:weak_problem} are discretised at the previous Picard iteration, becoming thus explicit. Consequently, at each Picard iteration, $k=1,\dots,N_{Pic}$, we have got a \textcolor{black}{symmetric and positive definite} system \textcolor{black}{for the pressure correction} that can be efficiently solved using classical numerical algorithms for linear systems like the conjugate gradient method using the solution at time $t^{n}$ as initial guess for $k=1$. \textcolor{black}{Since in the numerical tests of this paper only the ideal gas equation of state has been considered, no mass lumping was applied.} Let us also remark that for ODEs the Picard iteration procedure allows to gain one order in time per iteration so $N_{Pic}=2$ might be enough to keep the accuracy of the developed scheme. Once \eqref{eq:weak_problem} is solved, the solution is replaced into \eqref{eq:pressure_dt_p} to update the pressure and the enthalpy,
\begin{equation}
H^{n+1,k+1}= \frac{\gamma}{\gamma-1}\frac{\Press^{n+1,k+1}}{\rho^{n+1}},
\end{equation}
and the density kinetic energy,  
\begin{equation}
\left( \rho K\right) ^{n+1,k+1}=\frac{1}{2 \rho^{n+1}} \left|\Wvel^{n+1,k+1}\right|^{2},\quad \Wvel^{n+1,k+1} =  \ttWvel - \Delta t \grae \delta\Press^{n+1,k+1},
\end{equation}
can be computed to be used in the next Picard iteration.

\subsection{Post-projection stage}\label{sec:postpro}
Once the pressure $\Press^{n+1}$ at the new time $t^{n+1}$ and the new momentum $\Wvel^{n+1}$, are computed from \eqref{eq:pressure_dt_p} and \eqref{eq:momentum_dt_p}, the total energy density must be updated. 
Integrating equation \eqref{eq:energy_dt_p} on a spatial control volume $T_{k}$ of the primal mesh and applying Gauss theorem yields
\begin{equation}
\WEk^{n+1} = \tWEk - \frac{\Delta t}{\left|T_{k}\right|} \sum_{N_{l}\in\mathcal{K}_{k}} \int_{\Gamma_{kl}}  H^{n+1}\Wvel^{n+1} \cdot \boldsymbol{\tilde{\eta}}_{kl} dS.\label{eq:WEupdatek}
\end{equation}
where $\boldsymbol{\tilde{\eta}}_{kl}$ denotes the unit normal vector to the $l$-face $\Gamma_{kl}$ of the primal element $T_{k}$. To approximate the integral on the face we assume a constant value for the linear momentum given by its averaged value on the dual cell $C_{i}$ containing the face $\Gamma_{kl}$. Regarding the enthalpy, the averaged value of the vertex of the face is employed. The averaged total energy density at each primal element is then interpolated into the dual mesh using a weighted average,
\begin{equation}
\WEi^{n+1} = \frac{1}{\left|C_{i}\right|}\sum_{N_{k}\in\mathcal{K}_{i}} \left|T_{ki}\right|\WEk^{n+1},
\end{equation}
with $\left|T_{ki}\right|$ the area/volume of the intersection between the primal element $T_{k}$ and the dual element $C_{i}$.

\begin{remark}
	An alternative to update the total energy density would consist on computing on the primal grid only the contribution of the last term of \eqref{eq:WEupdatek}. Then, it can be interpolated to the dual grid and added to the intermediate value $\tWE$ originally computed on the dual grid. This approach would reduce dissipation arising from forward-backward interpolation between meshes.
\end{remark}

\subsection{Boundary conditions}\label{sec:bc}
Boundary conditions of the numerical tests to be presented in Section \ref{sec:numericalresults} are constructed as a combination of:
\begin{itemize}
	\item \textit{Periodic} boundary conditions. For the implementation of periodic boundary conditions, we assume that a periodic mesh is provided. The pairs of matching dual boundary elements are combined in order to define a new dual cell which then becomes of the interior type. Concerning the finite element method, the corresponding vertexes are merged, resulting in a reduced number of unknowns for the pressure system. All connections between the elements and nodes need to be updated accordingly. 
	
	\item \textit{Strong Dirichlet} boundary conditions on FV. The exact value at the boundary is imposed as the averaged value on the cell. Let us note that for inviscid flows only the normal component is set. When adiabatic walls are selected, we impose zero heat flux instead of defining the value of the density field. Definition of the exact value of the linear momentum is conveyed to the pressure system where it is used to compute the last term in equation \eqref{eq:weak_problem}. 
	
	\item \textit{Weak Dirichlet} boundary conditions on FV. The value of a variable on the boundary is employed to compute the contribution of the different terms of the corresponding conservative equations on boundary cells. Accordingly, the computation of gradients with the Galerkin approach makes use of the exact value of the variable at the boundary node whereas the numerical flux is computed considering an auxiliary state:
	\begin{equation}
	\WW_{\mathrm{aux}}^{\vel} = 2 \WW_{\mathrm{exact}}^{\vel} - \WW_{i}^{\vel} 
	\end{equation}
	in the viscous case and
	\begin{equation}	
	\WW_{\mathrm{aux}}^{\vel} = \WW_{i}^{\vel} -2 \WW_{i}^{\vel} \cdot \widetilde{\boldsymbol{\eta}} \widetilde{\boldsymbol{\eta}}
	\end{equation}
	for inviscid wall boundary conditions.
	Likewise strong Dirichlet boundary conditions, they are usually combined with Neumann boundary conditions for the pressure field.
	
	\item Neumann boundary conditions on FV. They are generally linked to Dirichlet boundary conditions on the pressure system so the exact value is imposed on boundary vertex. No further computations are needed for the definition of inflow and outflow conditions on the velocity field.
\end{itemize}


\section{Numerical results}\label{sec:numericalresults}
The developed methodology is assessed using classical benchmarks from the incompressible limit to high Mach number flows. For all tests presented in the following, we consider SI units. The time step is determined according to the  condition 
\begin{equation}
\Delta t = \min_{C_{i}}\left\lbrace \Delta t_{i}\right\rbrace, \qquad \Delta t_{i} = \textnormal{CFL} \frac{r_{i}^2}{( \left|\zeta\right|_{\max} + c_\alpha) r_{i} +2\left[  \dfrac{4}{3}\left( \left|\nu\right|_{\max}+\dfrac{\gamma \lambda}{c_{\press} \rho} \right) \right] }
\label{eqn.cfl.condition}
\end{equation}          
where $\left|\zeta\right|_{\max}$, $\left|\nu\right|_{\max}$ and $\lambda$ denote the maximum absolute eigenvalues related to the convective and diffusive terms, respectively, which have been discretized explicitly, \textcolor{black}{and $c_\alpha$ is an artificial viscosity parameter, which in this paper is assumed to be constant in space and time for simplicity. If not explicitly stated otherwise, $c_\alpha=0$ is used as default value. In \eqref{eqn.cfl.condition}, the symbol $r_i$ denotes the incircle diameter of each dual control volume. The default value for the CFL number for all test cases is CFL$=1/d$, with $d$ the number of space dimensions. Besides, the CFL related to the sound speed, $\textnormal{CFL}_{c}$, indicated in the tests is computed as
\begin{equation}
\textnormal{CFL}_{c}= \max_{C_{i}} \left( c_{i} \frac{\Delta t}{r_i}\right), \quad c_{i} = \sqrt{\gamma\frac{\Press_{i}}{\rho_i}}.
\end{equation}
In the rest of this section, gravity is neglected, hence the gravity vector is set to $\mathbf{g}=\mathbf{0}$.} 

\subsection{Taylor-Green vortex and numerical convergence results} \label{sec:TGV}
To study the accuracy of the new method proposed in this paper the Taylor-Green vortex problem is solved in 2D. 
We consider a computational domain $\Omega=\left[0,2\pi\right]\times\left[0,2\pi\right]$ discretized using the meshes described in Table \ref{TGV_mesh} and a final simulation time $t_{\mathrm{end}}=0.1$. 
\begin{table}[t]
	\renewcommand{\arraystretch}{1.2}
	\begin{center}
		\begin{tabular}{cccc}
			\hline 
			Mesh & Elements & Vertices & Dual elements \\\hline
			$M_1$ & $512 $ & $289 $ & $800 $ \\ 
			$M_2$ & $2048 $ & $1089 $ & $3136 $ \\ 
			$M_3$ & $8192 $ & $4225 $ & $12416 $ \\ 
			$M_4$ & $32768 $ & $16641 $ & $49408 $ \\ 
			$M_5$ & $131072 $ & $66049 $ & $197120 $ \\ 
			$M_6$ & $524288 $ & $263169 $ & $787456 $ \\ 
			\hline 
		\end{tabular}
		\caption{Different meshes used for the numerical convergence study based on the Taylor-Green vortex. } \label{TGV_mesh}
	\end{center}
\end{table}
The exact solution for this test case reads 
\begin{equation}
\rho\left(\mathbf{x},t\right) = 1,\quad
\vel \left(\mathbf{x},t\right) = \left( \begin{array}{r} 
  \sin(x)\cos(y) \\ 
 -\cos(x)\sin(y) \end{array} \right), \quad 
\press \left(\mathbf{x},t\right) = \frac{\press_{0}}{\gamma-1} + \frac{1}{4} \left(\cos(2x)+\cos(2y) \right), 
\end{equation}
with $\gamma=1.4$, $\press_{0}=10^{5}$. This leads to a characteristic Mach number $M\approx1.7 \cdot10^{-3}$, so it corresponds to the low Mach number regime. In this particular test case, the time step has not been automatically computed from the CFL condition, instead a fixed value starting from $\Delta t=0.25$ for the coarsest mesh and decreasing linearly according to the mesh size has been employed. The $L^{2}$ errors in space, computed at the final time step, and the order of accuracy attained,
\begin{gather}
L^{2}_{\Omega,\, M_{i}}\left(W\right) = \left\|W_{\mathrm{exact}}-W_{M_i} \right\|_{L^2(\Omega)},\quad
\mathcal{O}_{M_i/M_j} \left(W\right) = \frac{\log\left( L^{2}_{\Omega,\, M_{i}}\left(W\right)/L^{2}_{\Omega,\, M_{j}}\left(W\right)\right) }{\log\left( h_{M_i}/h_{M_j}\right) },\end{gather}
are depicted in Table \ref{TGV_errors}. The sought order of accuracy is reached for the main flow variables with both the first order and LADER schemes, \textcolor{black}{as for the hybrid FV/FE scheme for the incompressible Navier-Stokes equations proposed in \cite{BFTVC17}. At this point, we would like to remark that the proposed scheme is nominally only first order accurate in time due to the employed operator splitting technique. In order to achieve high order also in time, we recommend the use of an IMEX Runge-Kutta scheme, see e.g. \cite{PareschiRusso2000,RussoAllMach,DLDV18,BDLTV2020,BDT2021,BP2021}. }    
\begin{table}[ht]
	\renewcommand{\arraystretch}{1.2}
	\begin{center}
		\begin{tabular}{ccccccccc}
			\hline 
		    Mesh &$L^{2}_{\Omega}\left(\rho\right)$ & $\mathcal{O}\left(\rho\right)$                  
			&$L^{2}_{\Omega}\left(\vel\right)$ & $\mathcal{O}\left(\vel\right)$ & $L^{2}_{\Omega}\left(E\right)$ & $\mathcal{O}\left(E\right)$  & $L^{2}_{\Omega}\left(\press\right)$ & $\mathcal{O}\left(\press\right)$ \\ \hline
			\multicolumn{9}{c}{First order scheme}\\ \hline
			M1 & $1.51E\!-\!01$ & $    $ &$2.23E\!-\!01$ & $    $ & $1.74E\!+\!03$ & $    $ & $2.77E\!-\!01$ & $    $ \\
			M2 & $8.19E\!-\!02$ & $0.88$ &$1.13E\!-\!01$ & $0.98$ & $3.74E\!+\!02$ & $2.22$ & $5.47E\!-\!02$ & $2.34$ \\
			M3 & $4.12E\!-\!02$ & $0.99$ &$5.66E\!-\!02$ & $0.99$ & $8.71E\!+\!01$ & $2.10$ & $2.95E\!-\!02$ & $0.89$ \\
			M4 & $2.06E\!-\!02$ & $1.00$ &$2.84E\!-\!02$ & $1.00$ & $2.13E\!+\!01$ & $2.03$ & $1.55E\!-\!02$ & $0.93$ \\
			M5 & $1.03E\!-\!02$ & $1.00$ &$1.42E\!-\!02$ & $1.00$ & $5.28E\!+\!00$ & $2.01$ & $7.96E\!-\!03$ & $0.96$ \\
			M6 & $5.14E\!-\!03$ & $1.00$ &$7.11E\!-\!03$ & $1.00$ & $1.32E\!+\!00$ & $2.00$ & $4.03E\!-\!03$ & $0.98$ \\ \hline
			\multicolumn{9}{c}{LADER scheme}\\ \hline
			M1 & $4.23E\!-\!03$ & $    $ &$3.06E\!-\!02$ & $    $ & $7.20E\!+\!01$ & $    $ & $1.10E\!-\!01$ & $    $ \\
			M2 & $9.71E\!-\!04$ & $2.12$ &$7.69E\!-\!03$ & $1.99$ & $3.82E\!+\!00$ & $4.24$ & $2.99E\!-\!02$ & $1.87$ \\
			M3 & $2.37E\!-\!04$ & $2.04$ &$1.92E\!-\!03$ & $2.00$ & $2.29E\!-\!01$ & $4.06$ & $7.56E\!-\!03$ & $1.98$ \\
			M4 & $5.84E\!-\!05$ & $2.02$ &$4.78E\!-\!04$ & $2.00$ & $1.50E\!-\!02$ & $3.93$ & $1.90E\!-\!03$ & $1.99$ \\
			M5 & $1.45E\!-\!05$ & $2.01$ &$1.19E\!-\!04$ & $2.00$ & $1.21E\!-\!03$ & $3.64$ & $4.76E\!-\!04$ & $2.00$ \\
			M6 & $3.60E\!-\!06$ & $2.01$ &$2.99E\!-\!05$ & $2.00$ & $1.63E\!-\!04$ & $2.89$ & $1.19E\!-\!04$ & $2.00$ \\
			\hline 
		\end{tabular}
		\caption{Taylor-Green vortex. Spatial $L_{2}$ error norms obtained at time $t_{\mathrm{end}}=0.1$, and convergence rates for the first order and local ADER schemes.} \label{TGV_errors}
	\end{center}
\end{table}

\subsection{Riemann problems}
\begin{figure}[!htbp]
	\centering
	\includegraphics[trim= 5 5 5 5,clip,width=0.325\linewidth]{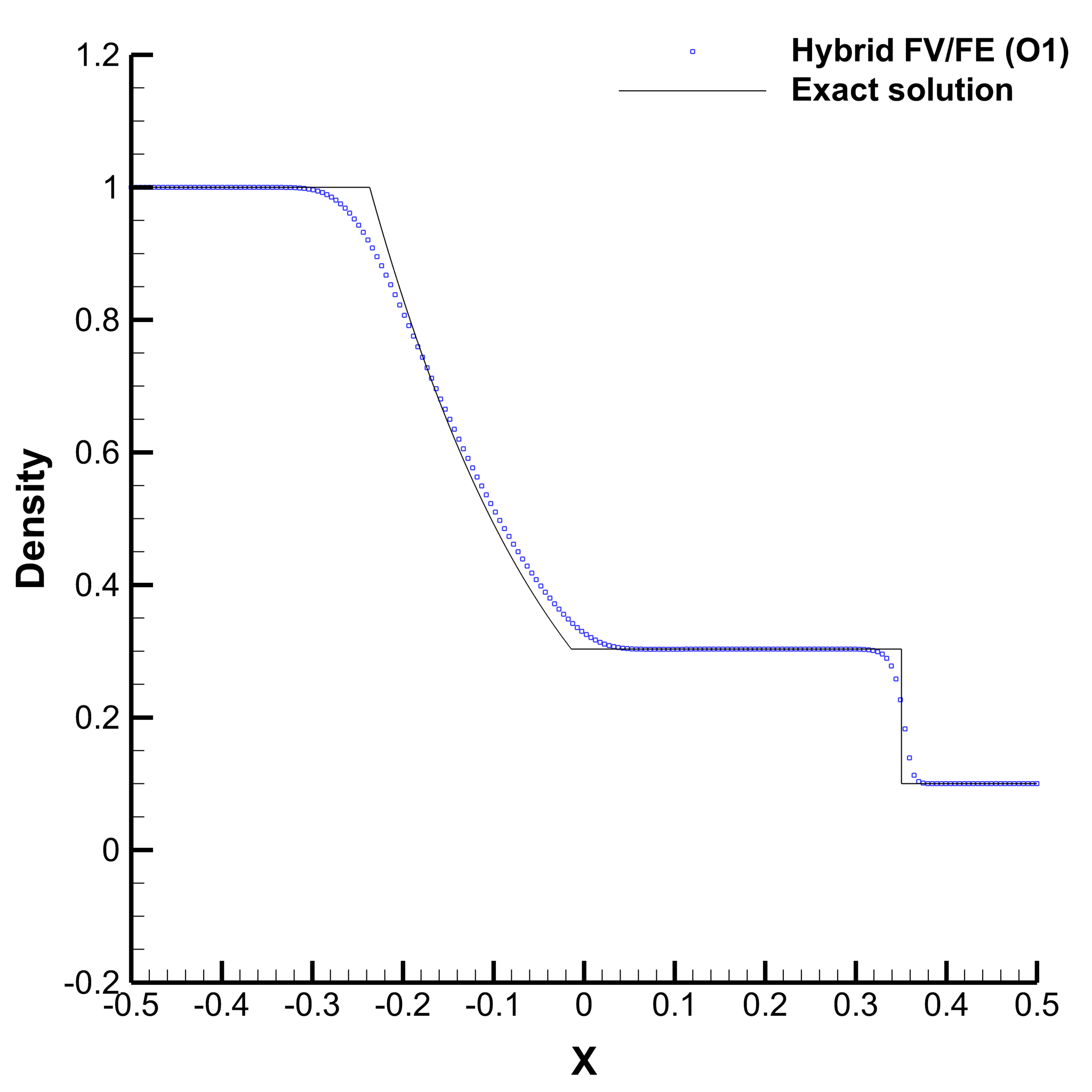}
	\includegraphics[trim= 5 5 5 5,clip,width=0.325\linewidth]{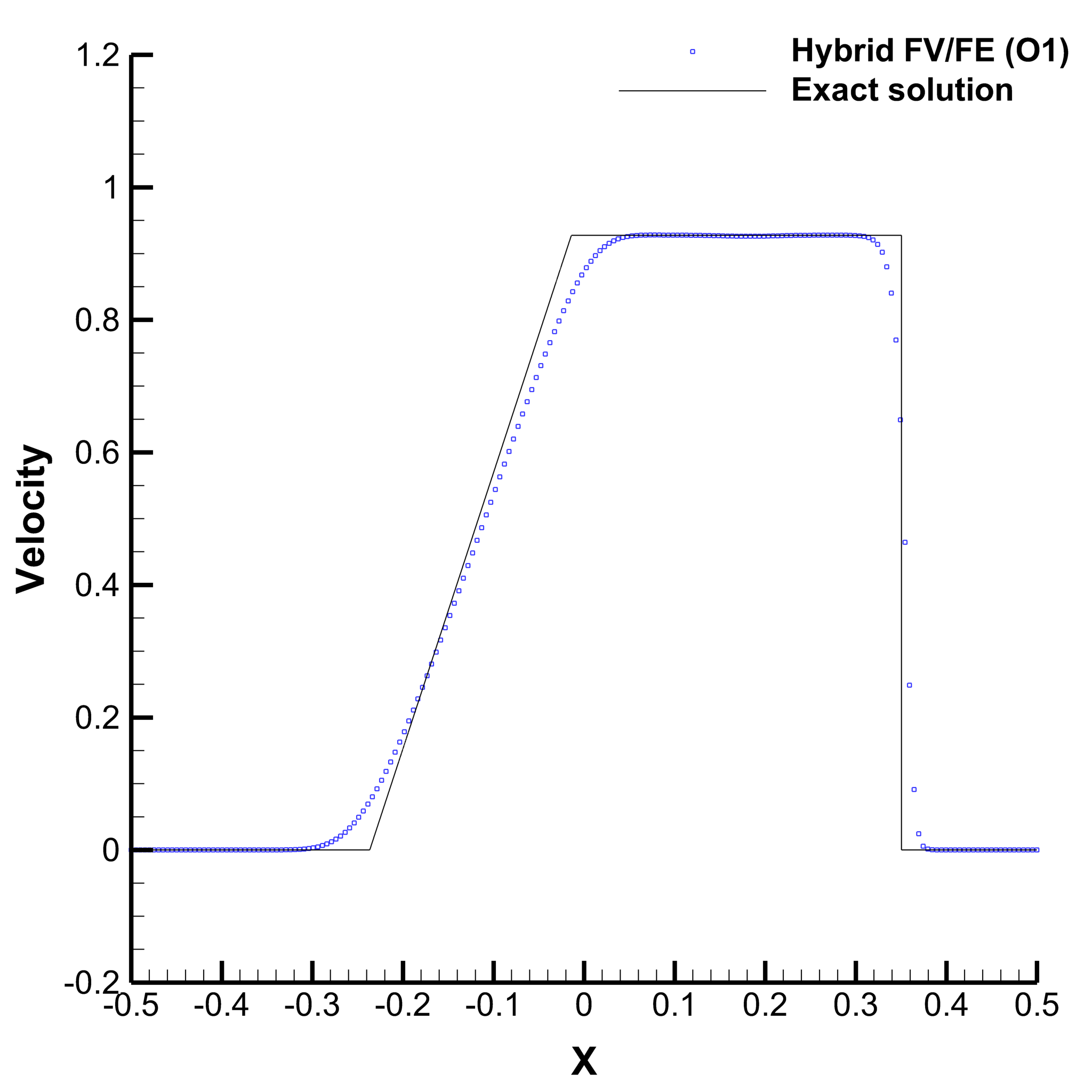}
	\includegraphics[trim= 5 5 5 5,clip,width=0.325\linewidth]{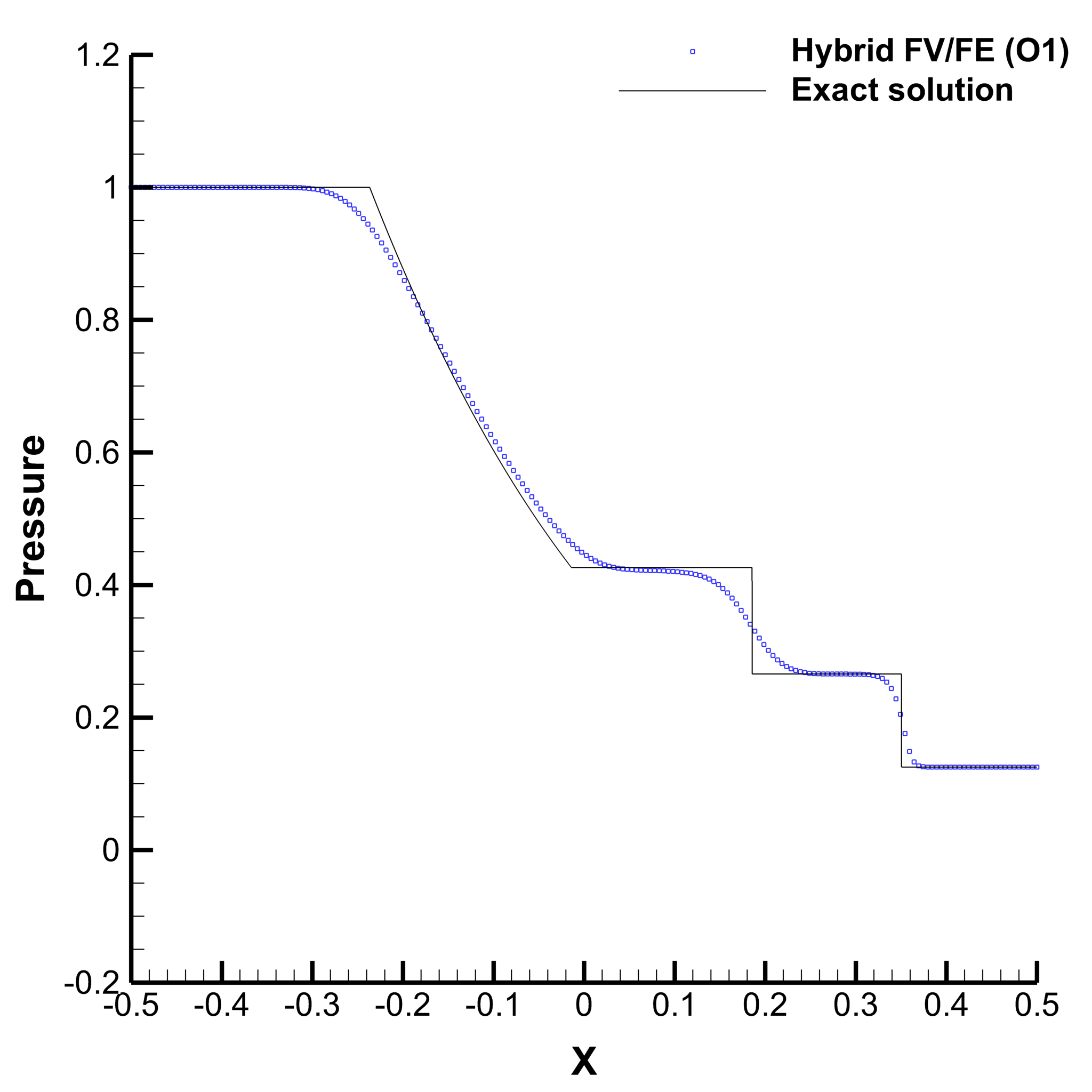}
	
	\caption{Riemann problem 1 (Sod). 1D cut through the numerical results along the line $y=0$ for $\rho$, $u$ and $p$ at $t_{\mathrm{end}}=0.2$ using the first order method  ($\mathrm{CFL}_{c}=3.35$, $c_{\alpha}=1$, $M\approx 0.93$).}
	\label{fig:RP1_o1}
\end{figure}

In this section, we analyse the performance of the proposed methodology for the compressible Euler equations in presence of medium to strong shocks. We consider a two-dimensional computational domain with $x \in [-0.5,0.5]$ and a variable width depending on the number of cells in the horizontal direction so that the final elements have a good aspect ratio and a small number of layers in the $y$-direction to reduce the computational cost of the simulation. The initial condition is defined as 
\begin{equation}
\rho^{0}\left(\mathbf{x}\right) = \left\lbrace \begin{array}{lc}
\rho_{L} & \mathrm{ if } \; x \le  x_{c},\\
\rho_{R} & \mathrm{ if } \; x>  x_{c};
\end{array}\right.\qquad
{u}_{1}^{0} \left(\mathbf{x}\right) =\left\lbrace \begin{array}{lc}
u_{L}  & \mathrm{ if } \; x \le  x_{c},\\
u_{R} & \mathrm{ if } \; x >  x_{c};
\end{array}\right. \qquad
{u}_{2}^{0} \left(\mathbf{x}\right) =0\qquad
\press^{0} \left(\mathbf{x}\right) = \left\lbrace \begin{array}{lc}
\press_{L} & \mathrm{ if } \; x \le  x_{c},\\
\press_{R} & \mathrm{ if } \; x >  x_{c};
\end{array}\right. 
\end{equation}
where $\rho_{L}$, $\rho_{R}$, $\press_{L}$, $\press_{R}$, $u_{L}$, $u_{R}$, $x_{c}$ are summarized in Table \ref{tab:RP_IC} for the diverse tests, selected among those presented in \cite{Toro,TD17}. The final time of each simulation, as well as the number of mesh divisions along the $x$-axis \textcolor{black}{($N_x$)}, have also been reported in Table \ref{tab:RP_IC}. \textcolor{black}{The characteristic mesh spacing is therefore equal to $h=1/N_x$}. All tests have been run with the first order and LADER schemes using Dirichlet boundary conditions in the $x$-direction and periodic boundary conditions in the $y$-direction. 
\begin{table}[ht]
	\renewcommand{\arraystretch}{1.2}
	\begin{center}
		\begin{tabular}{cccccccccc}
			Test &  $\rho_{L}$ &  $\rho_{R}$  &  $u_{L}$ &  $u_{R}$ &  $\press_{L}$ &  $\press_{R}$ & $x_{c}$ & $t_{\mathrm{end}}$ & $N_{x}$ \\ \hline
			RP1 & $ 1 $ & $ 0.125 $ & $ 0 $ & $ 0 $  & $ 1 $ & $ 0.1 $& $ 0 $ &  $ 0.2 $  & $200$\\  
			RP2 & $ 1 $ & $ 1 $  & $ -1 $ & $ 1 $ & $ 0.4 $ & $ 0.4 $& $ 0 $ &  $ 0.15 $  & $300$\\ 
			RP3 & $ 0.445 $ & $ 0.5 $ & $ 0.698 $ & $ 0 $ & $ 3.528 $ & $ 0.571 $ & $ 0 $ &  $ 0.14 $  & $200$\\
			RP4 & $ 5.99924 $ & $ 5.99242 $  & $ 19.5975 $ & $ -6.19633 $ & $ 460.894 $ & $ 46.095 $& $ -0.2 $ &  $ 0.035 $  & $200$\\ 
			RP5 & $ 1 $ & $ 1 $ & $ -19.59745 $ & $ -19.59745 $ & $ 1000.0 $ & $ 0.01 $ & $ 0.3 $ &  $ 0.01 $ & $300$ \\
			RP6 & $ 1 $ & $ 1 $ &  $ 2 $ & $ -2 $ &$ 0.1 $ & $ 0.1 $ & $ 0 $ &  $ 0.8 $ & $200$ \\
		\end{tabular} 
	\end{center}
	\caption{Riemann problems. Initial condition, initial position of the discontinuity, $x_{c}$, final time, $t_{\mathrm{end}}$, and number of mesh cells on $x$-direction, $N_{x}$, for each Riemann problem.}
	\label{tab:RP_IC}
\end{table}
The first test analysed, RP1, is the classical Sod problem presented for the first time in \cite{Sod78}. Figure \ref{fig:RP1_o1} shows a good agreement between the numerical and the exact solution for the shock, the contact, and the rarefaction waves. 
\begin{figure}[!htbp]
	\centering
	\includegraphics[trim= 5 5 5 5,clip,width=0.325\linewidth]{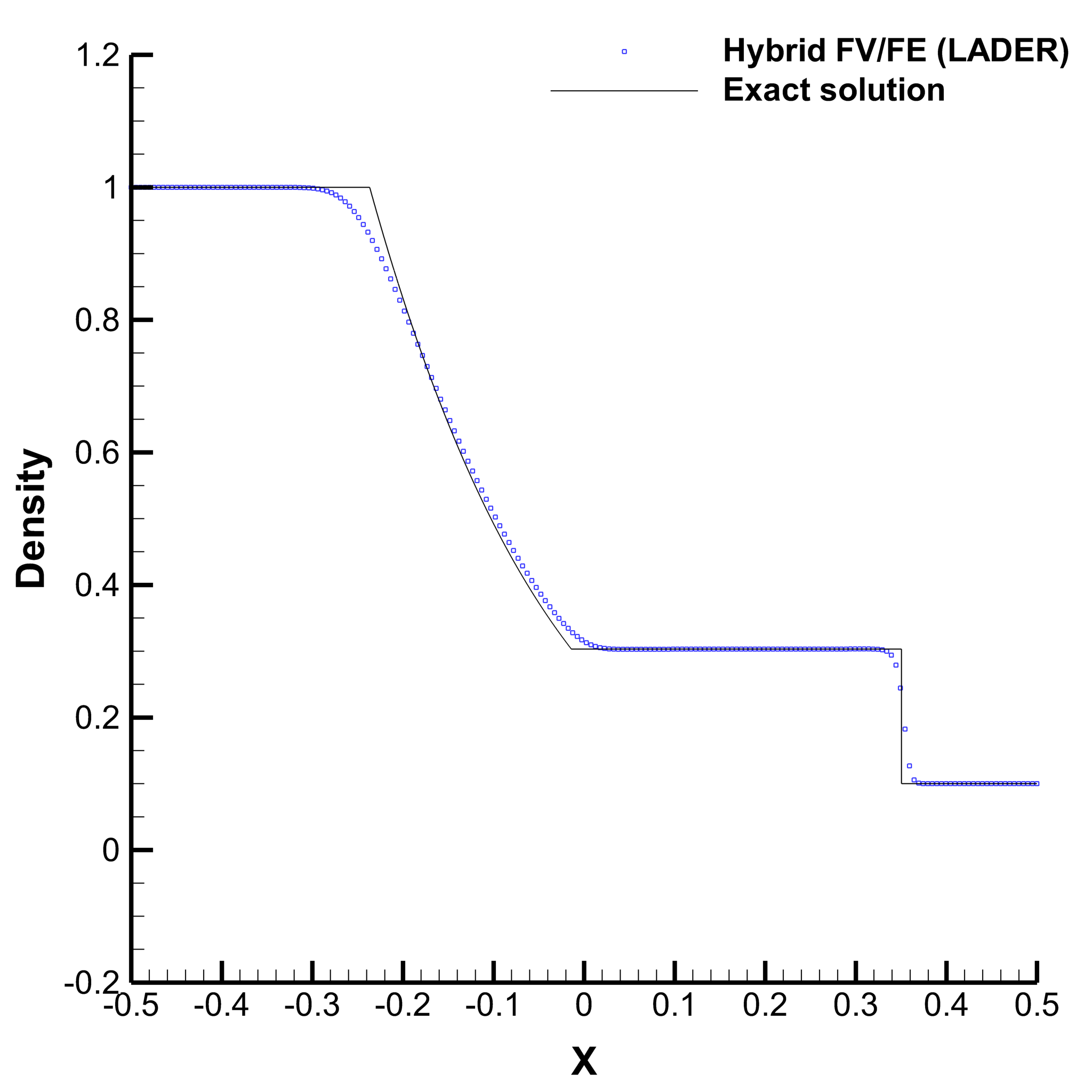}
	\includegraphics[trim= 5 5 5 5,clip,width=0.325\linewidth]{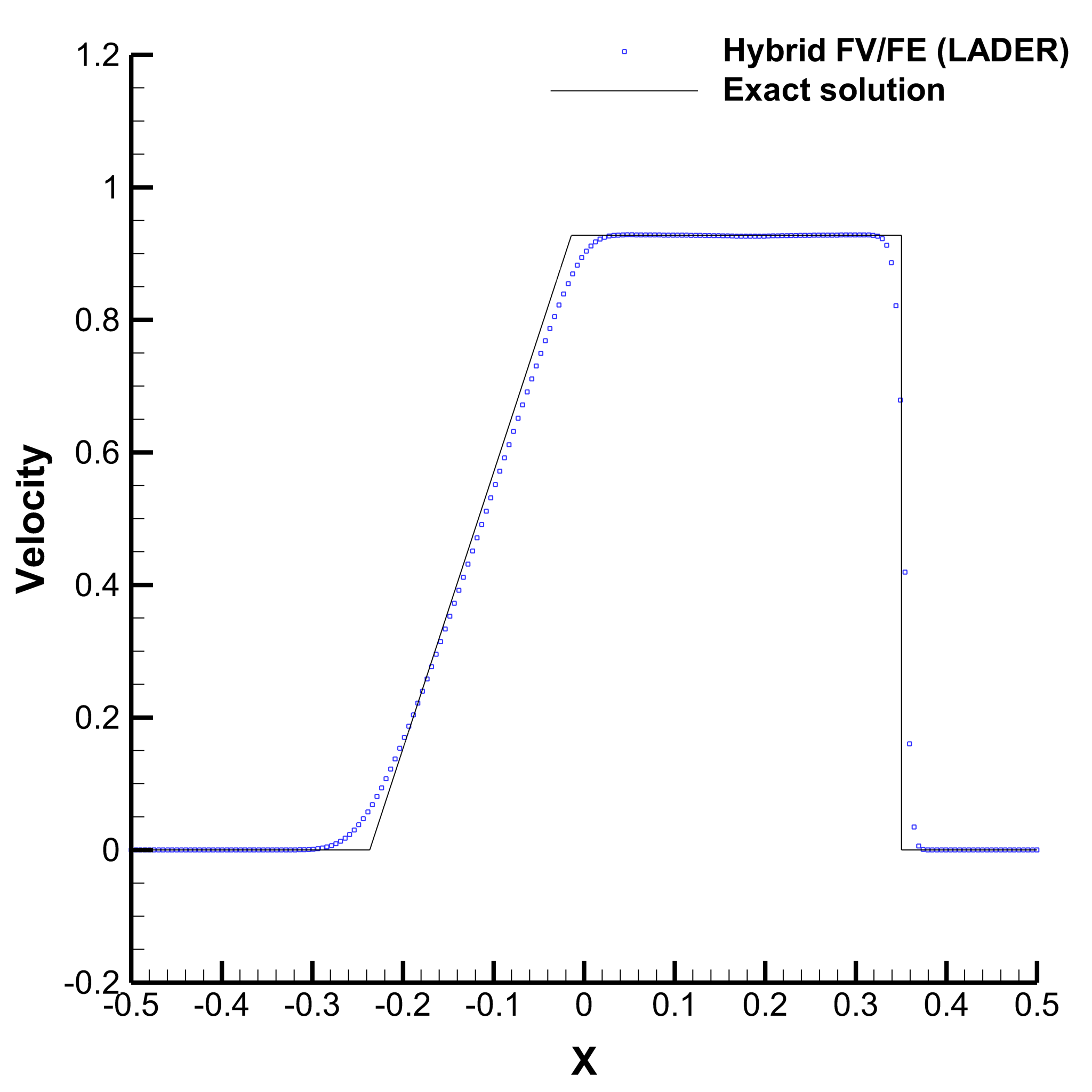}
	\includegraphics[trim= 5 5 5 5,clip,width=0.325\linewidth]{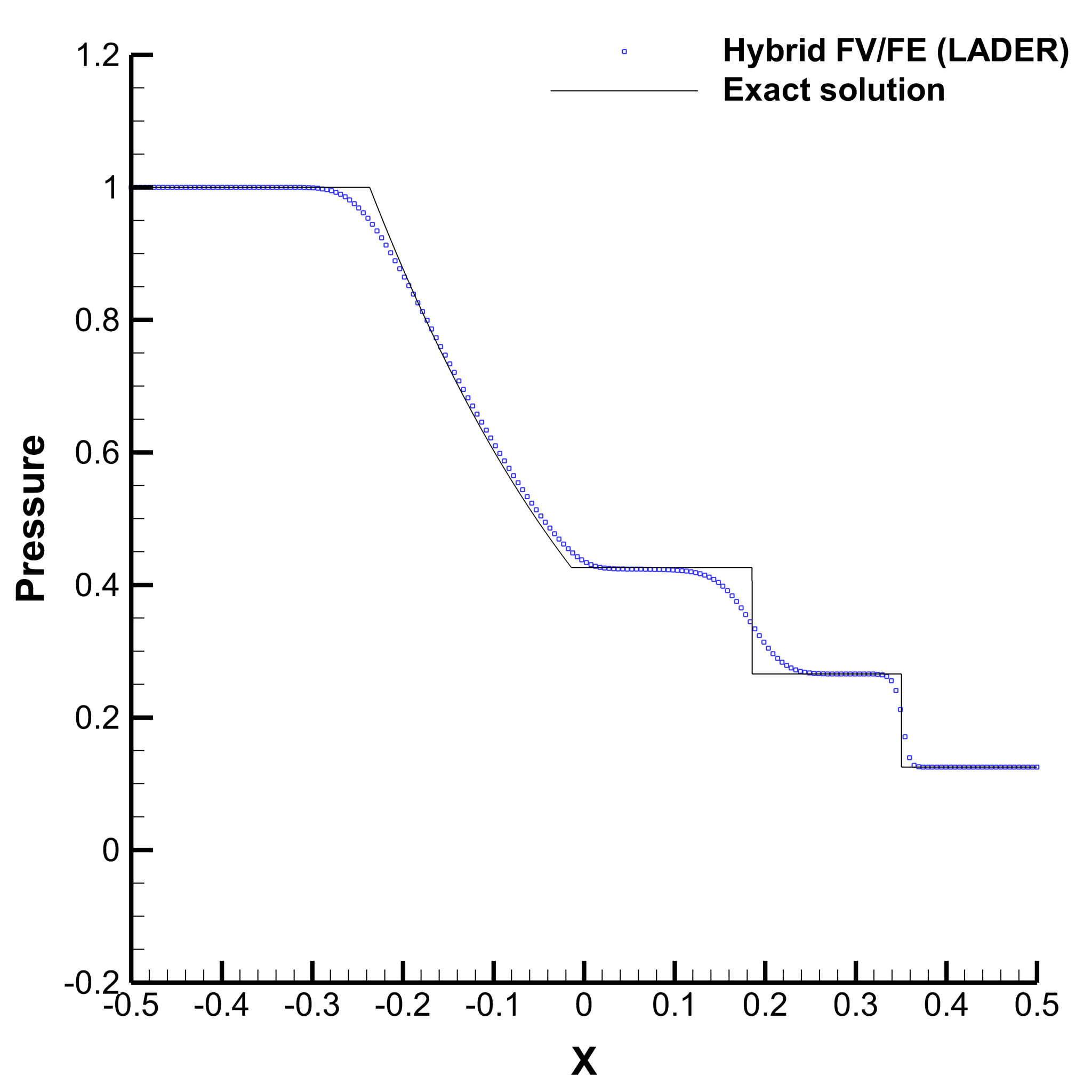}
	
	\caption{Riemann problem 1 (Sod). 1D cut through the numerical results along the line $y=0$ for $\rho$, $u$ and $p$ at $t_{\mathrm{end}}=0.2$ using the LADER-ENO method  ($\mathrm{CFL}_{c}=3.35$, $c_{\alpha}=1$, $M\approx 0.93$).}
	\label{fig:RP1_lader}
\end{figure}

\begin{figure}[!htbp]
	\centering
	\includegraphics[trim= 5 5 5 5,clip,width=0.325\linewidth]{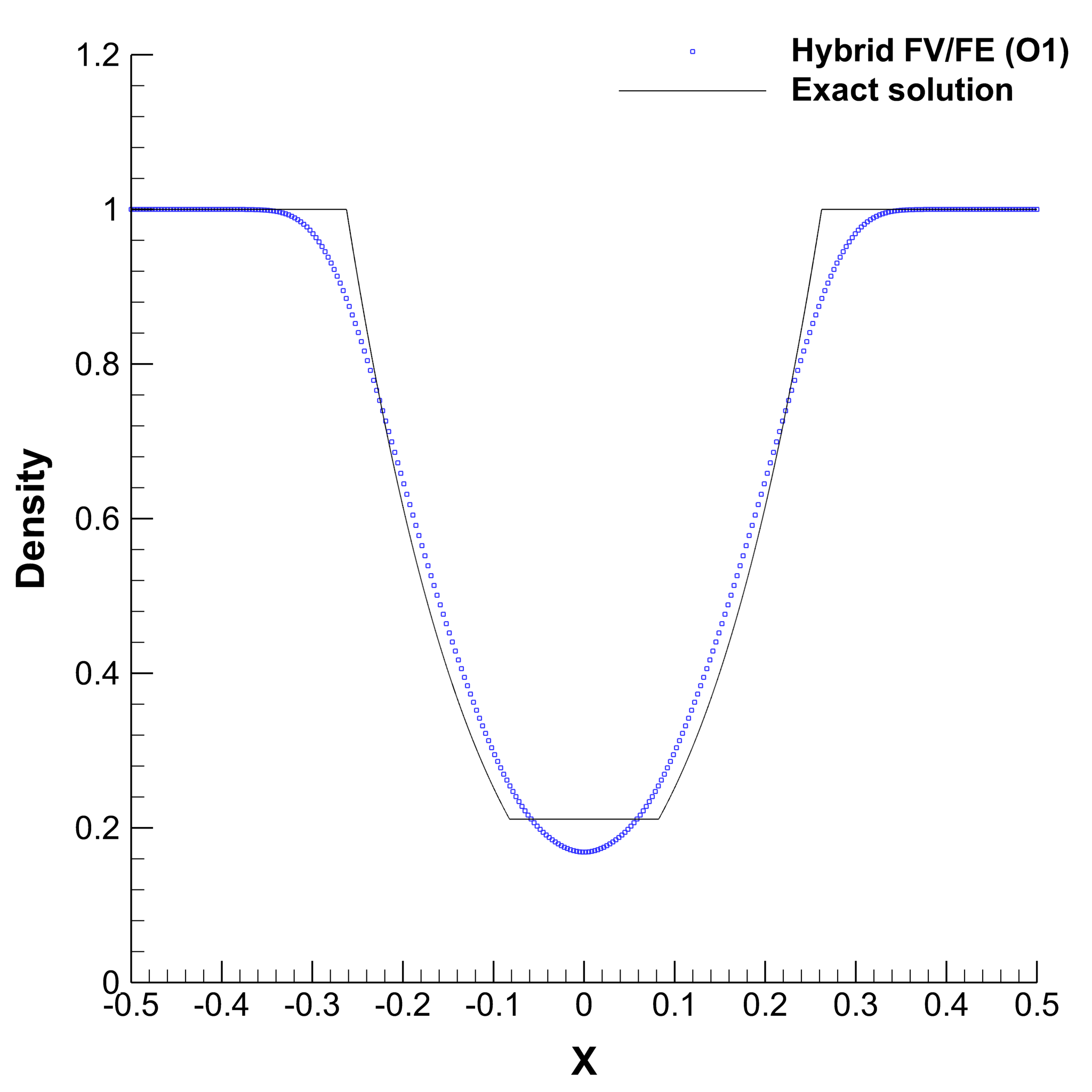}
	\includegraphics[trim= 5 5 5 5,clip,width=0.325\linewidth]{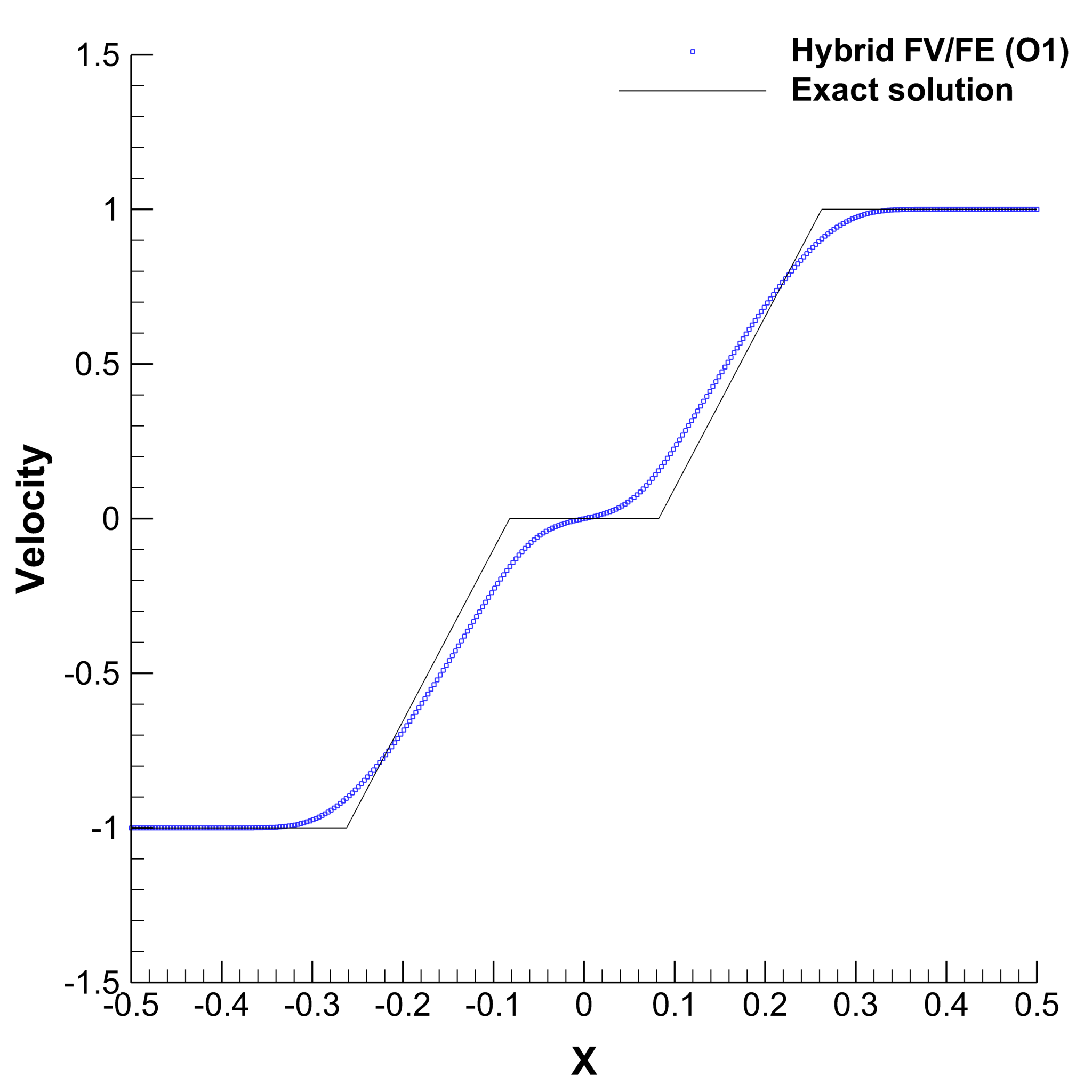}
	\includegraphics[trim= 5 5 5 5,clip,width=0.325\linewidth]{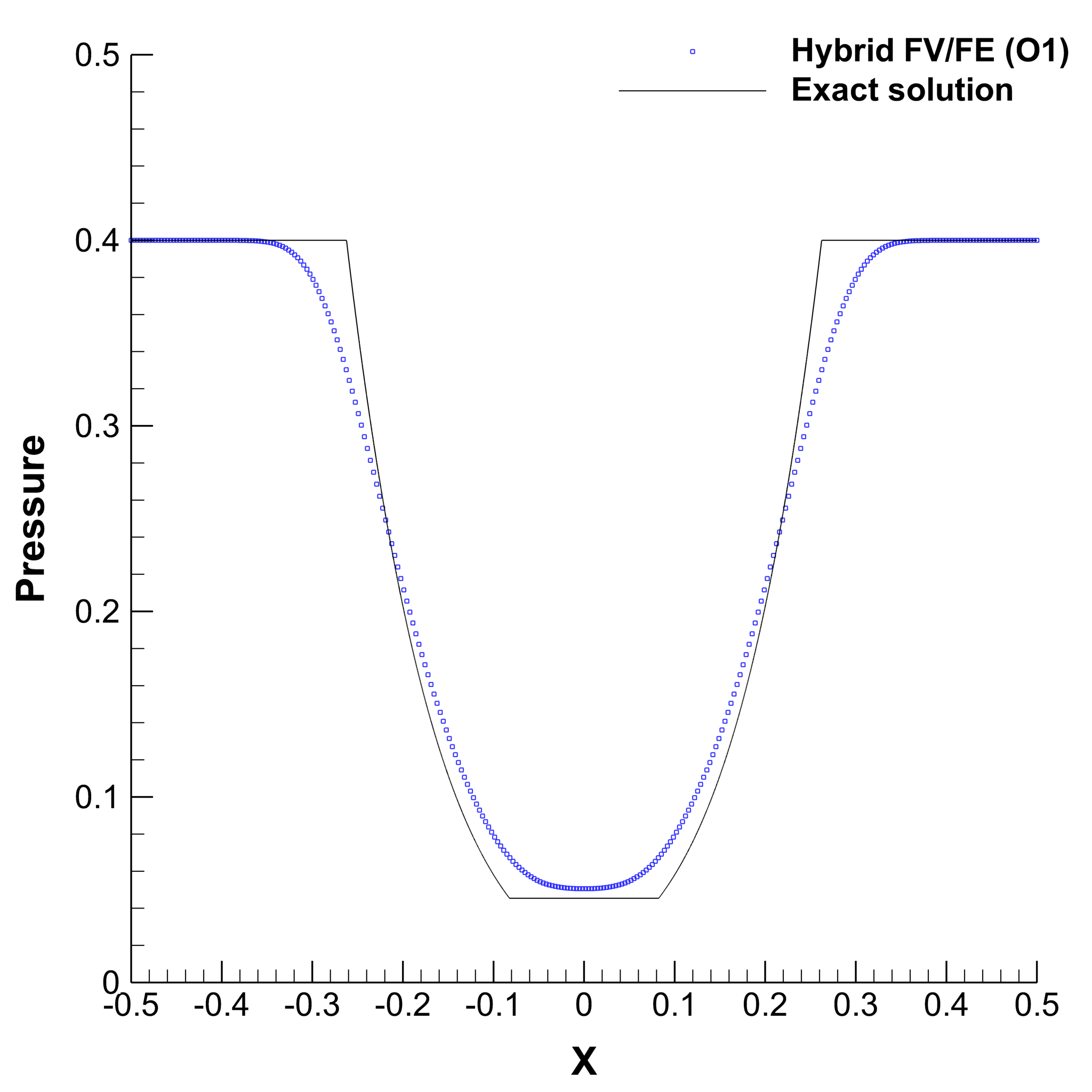}
	
	\caption{Riemann problem 2 (Double rarefaction). 1D cut through the numerical results along the line $y=0$ for $\rho$, $u$ and $p$ at $t_{\mathrm{end}}=0.15$ using the first order method  ($\mathrm{CFL}_{c}=0.4$, $c_{\alpha}=2$, $M\approx 1.37$).}
	\label{fig:RP2_o1}
\end{figure}

\begin{figure}[!htbp]
	\centering
	\includegraphics[trim= 5 5 5 5,clip,width=0.325\linewidth]{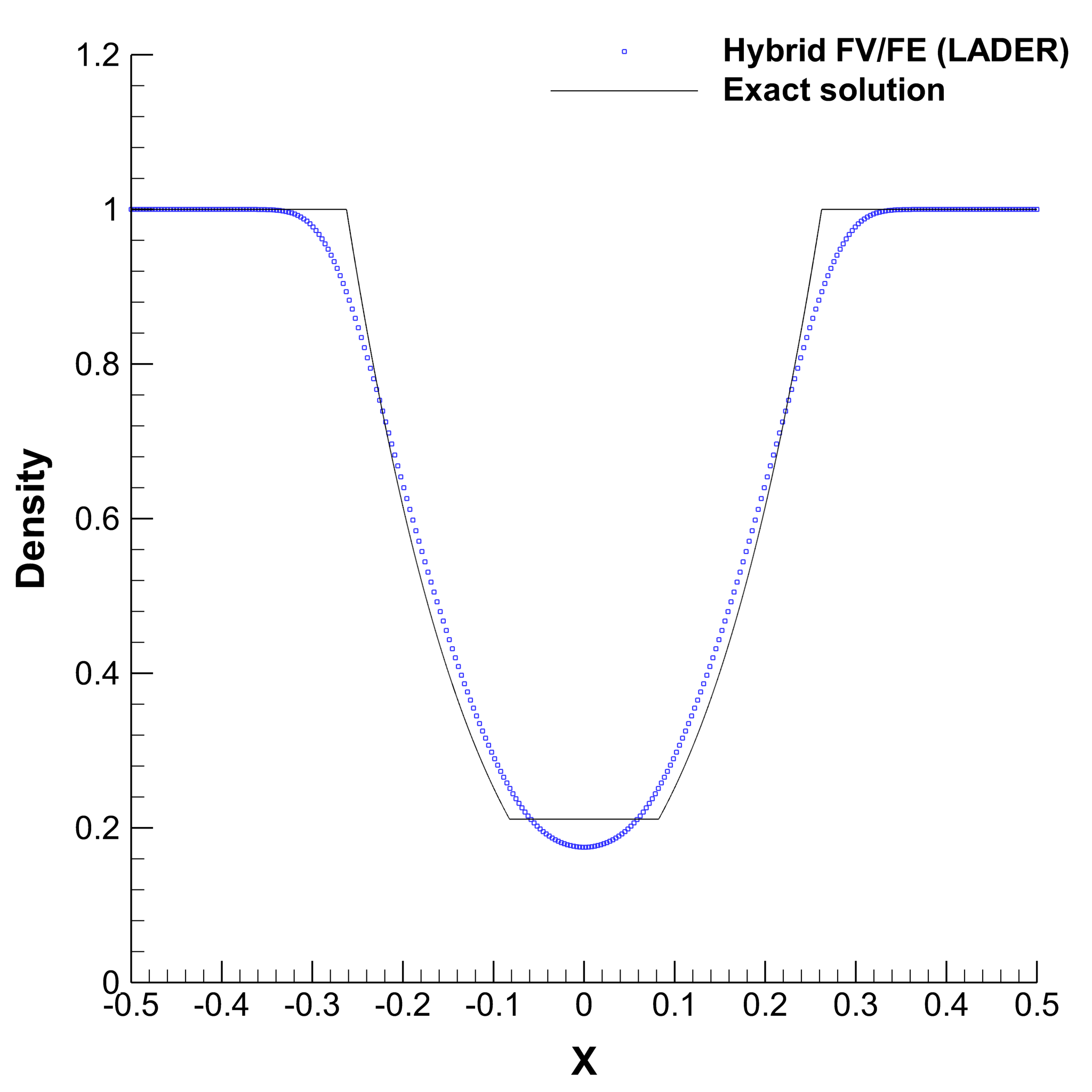}
	\includegraphics[trim= 5 5 5 5,clip,width=0.325\linewidth]{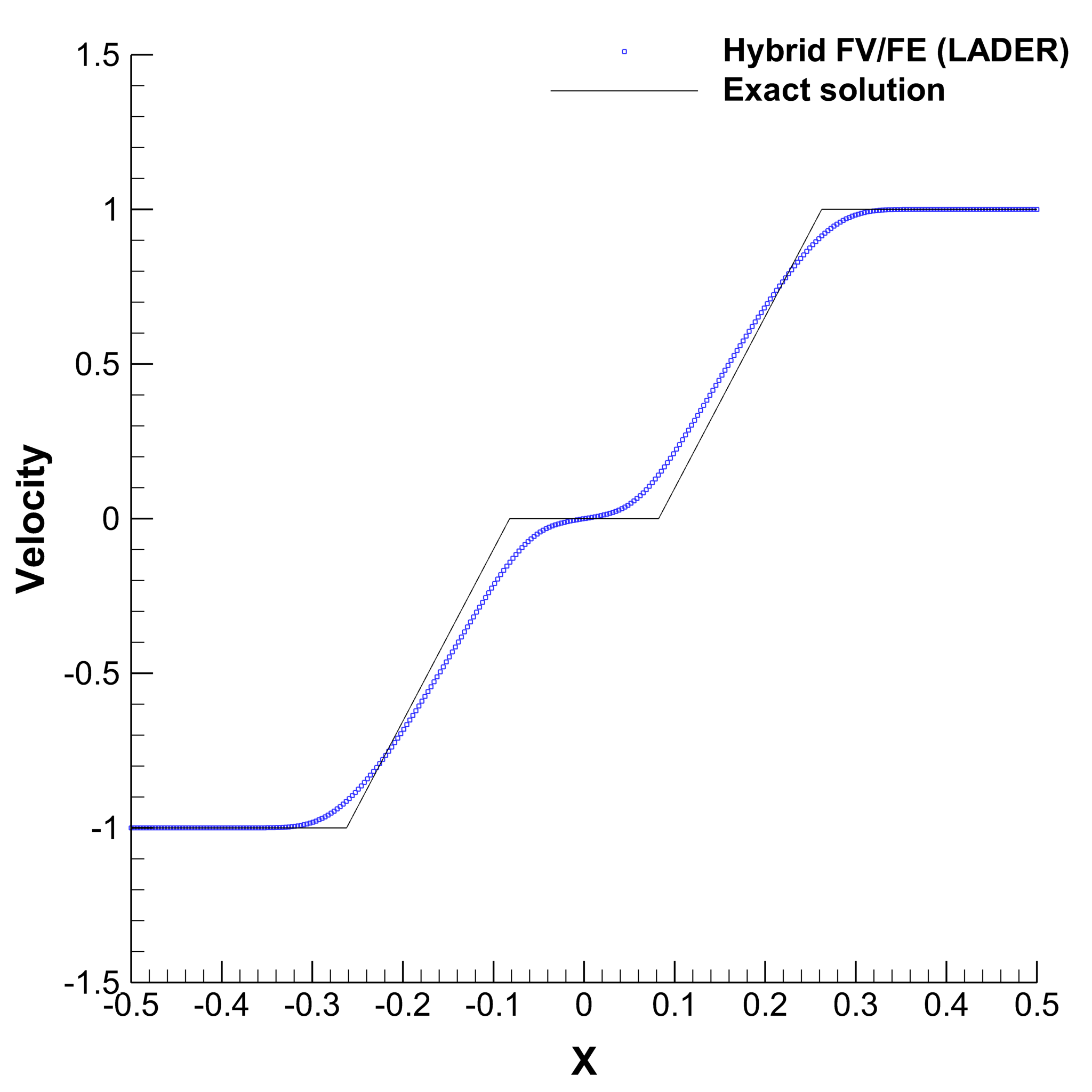}
	\includegraphics[trim= 5 5 5 5,clip,width=0.325\linewidth]{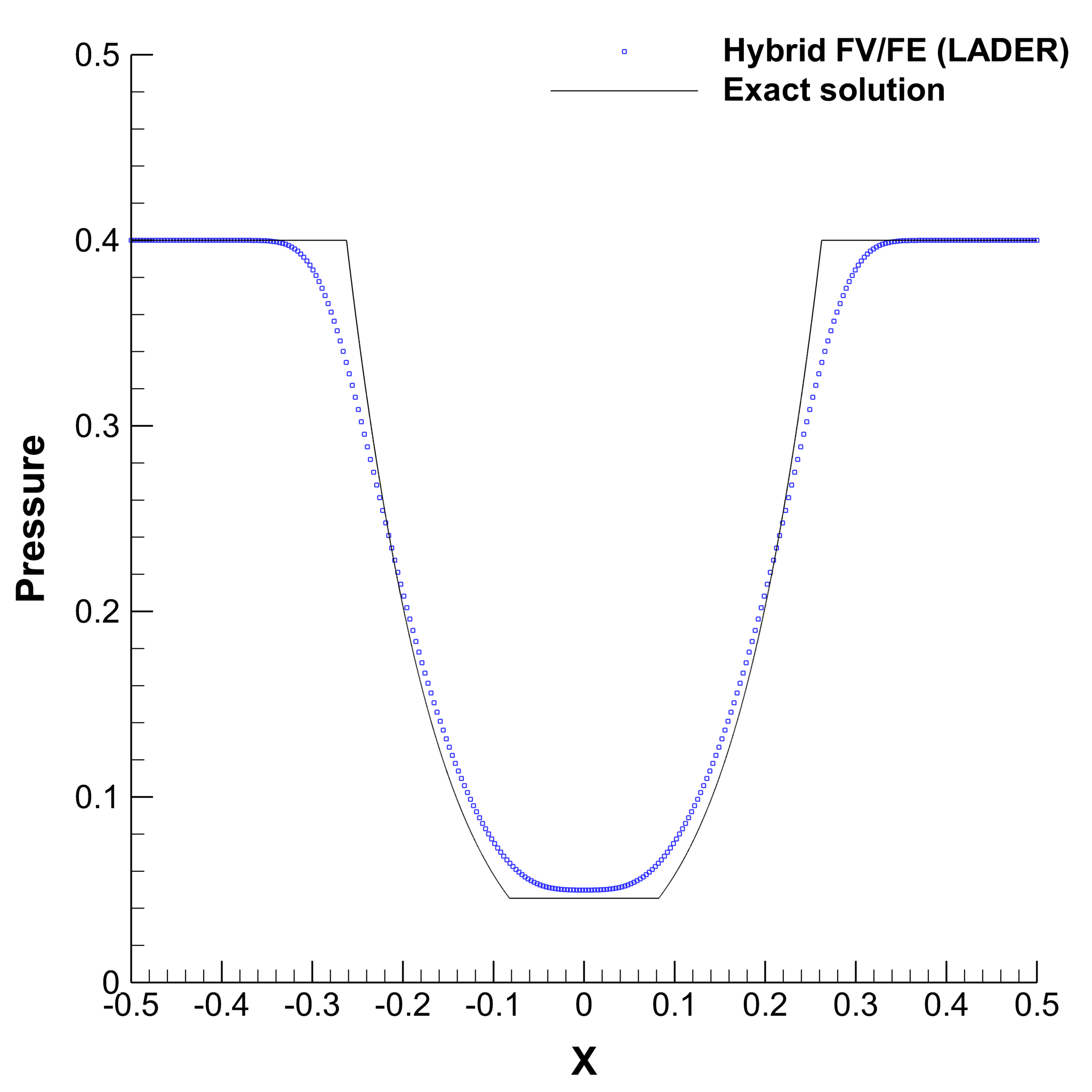}
	
	\caption{Riemann problem 2 (Double rarefaction). 1D cut through the numerical results along the line $y=0$ for $\rho$, $u$ and $p$ at $t_{\mathrm{end}}=0.15$ using the LADER-ENO method  ($\mathrm{CFL}_{c}=0.4$, $c_{\alpha}=2$, $M\approx 1.37$).}
	\label{fig:RP2_lader}
\end{figure}

RP2 corresponds to a double rarefaction problem. Overall the shape of the exact solution is captured even if a finer mesh would be useful to better approximate the constant contact discontinuity between the two rarefactions, Figures \ref{fig:RP2_o1}-\ref{fig:RP2_lader}. 

\begin{figure}[!htbp]
	\centering
	\includegraphics[trim= 5 5 5 5,clip,width=0.325\linewidth]{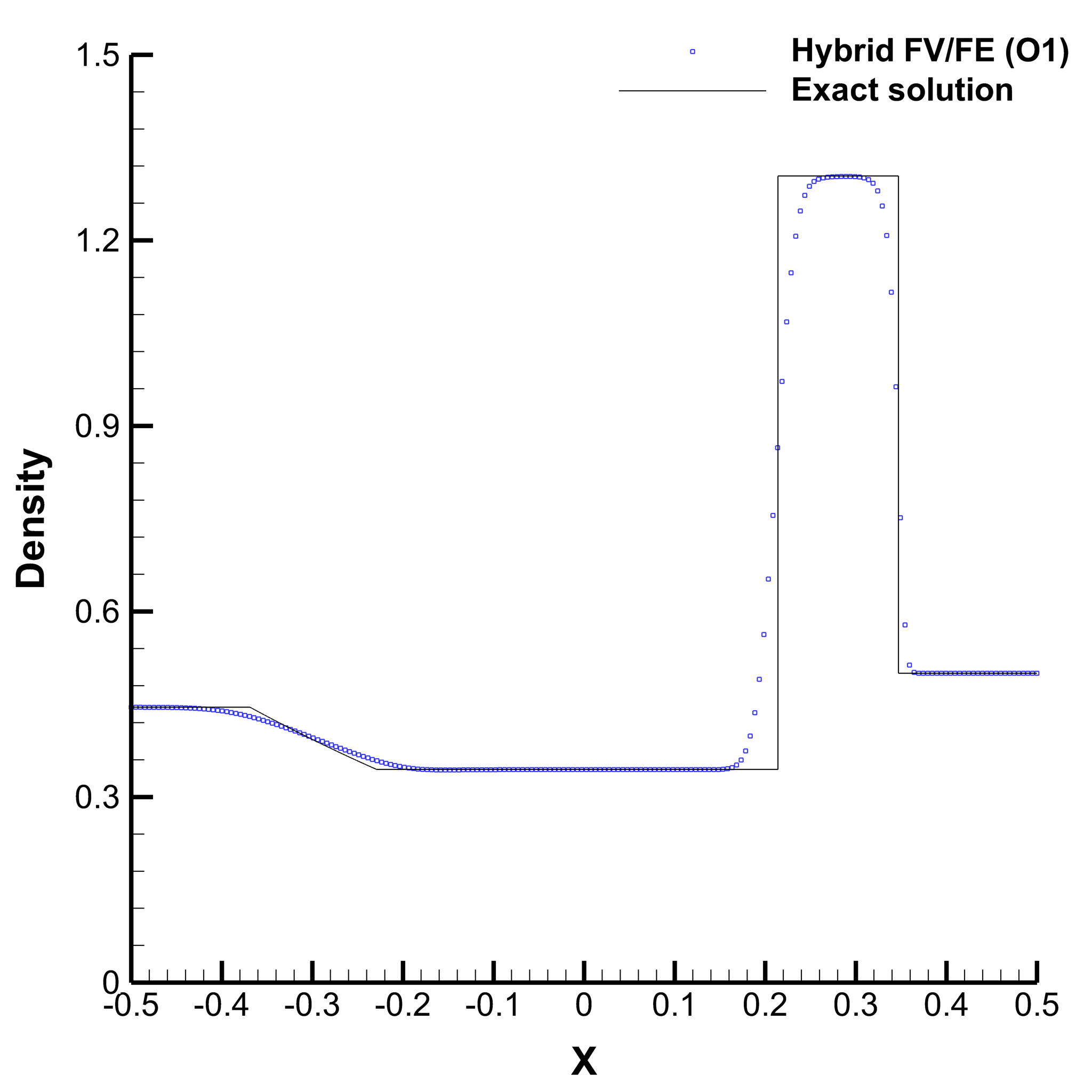}
	\includegraphics[trim= 5 5 5 5,clip,width=0.325\linewidth]{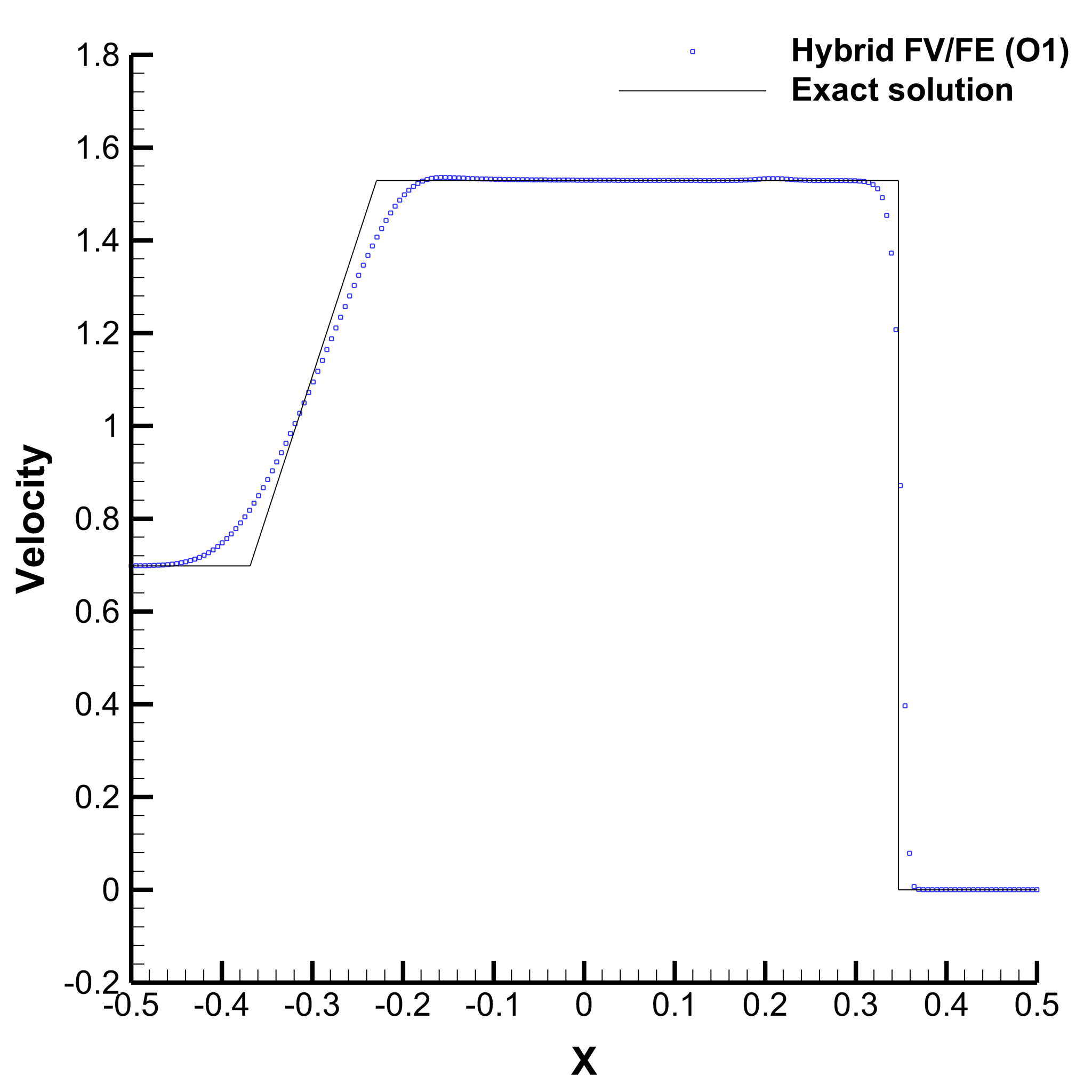}
	\includegraphics[trim= 5 5 5 5,clip,width=0.325\linewidth]{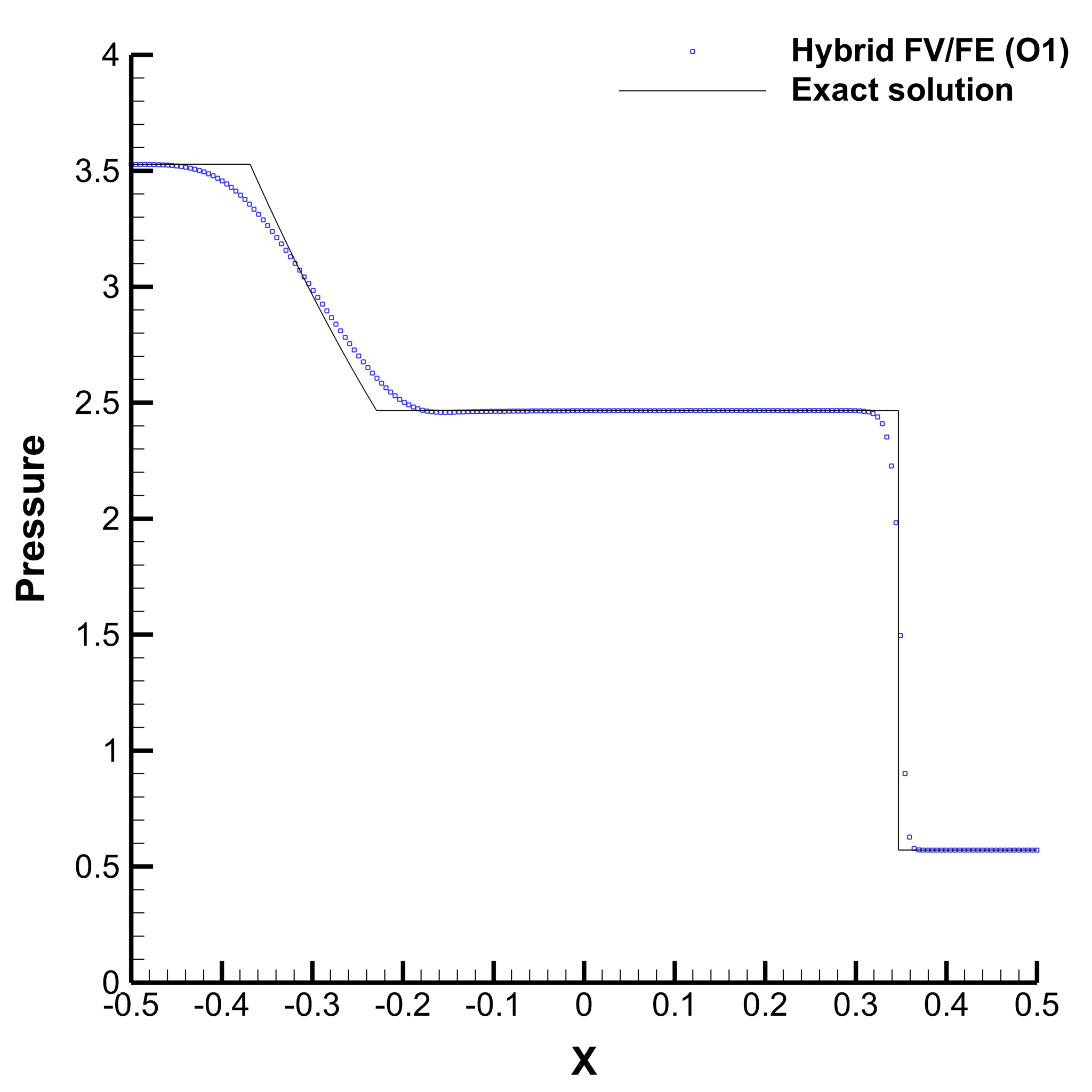}
	
	\caption{Riemann problem 3 (Lax). 1D cut through the numerical results along the line $y=0$ for $\rho$, $u$ and $p$ at $t_{\mathrm{end}}=0.14$ using the first order method on mesh M1  ($\mathrm{CFL}_{c}=2.78$, $M\approx 0.94$).} 
	\label{fig:RP3_o1}
\end{figure}

\begin{figure}[!htbp]
	\centering
	\includegraphics[trim= 5 5 5 5,clip,width=0.325\linewidth]{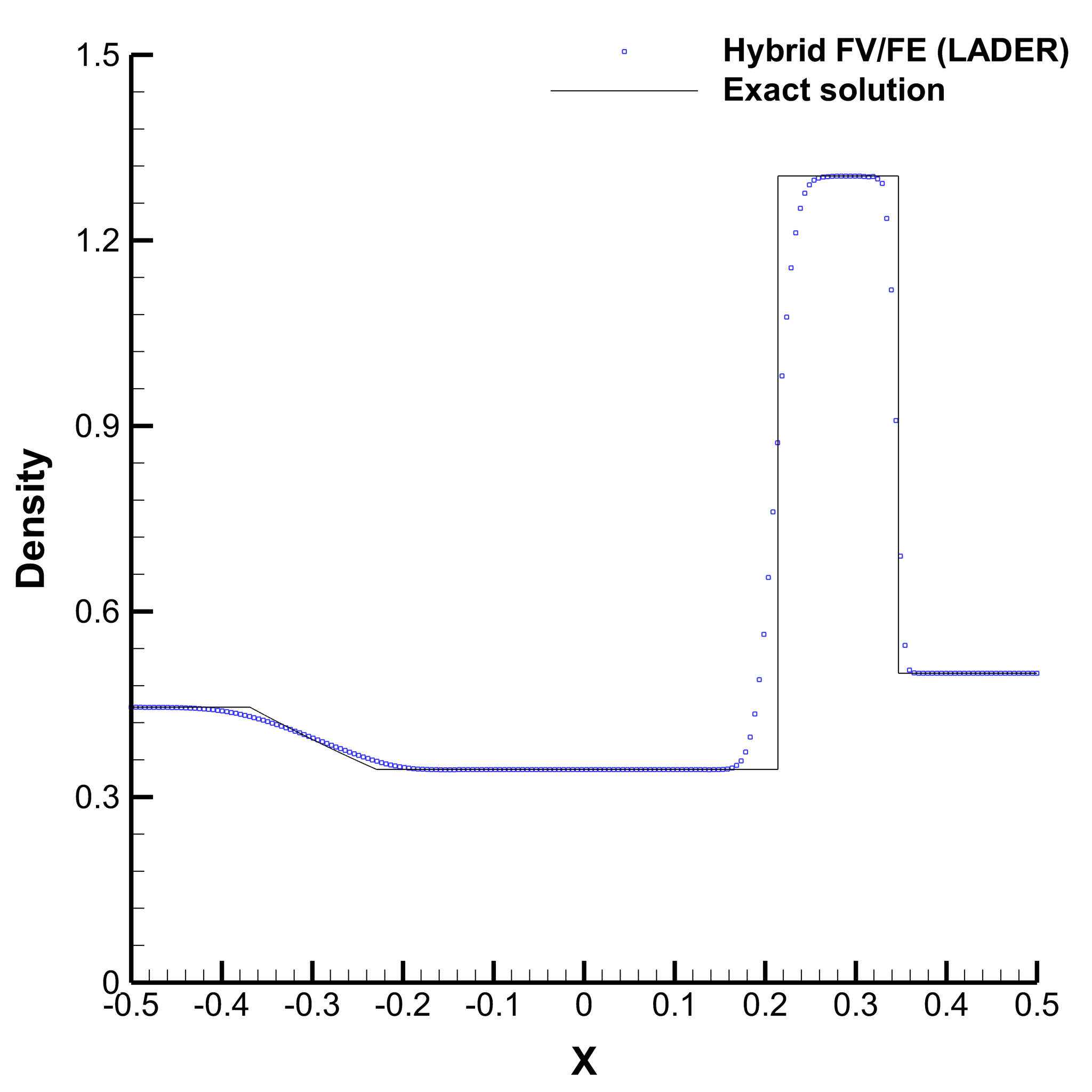}
	\includegraphics[trim= 5 5 5 5,clip,width=0.325\linewidth]{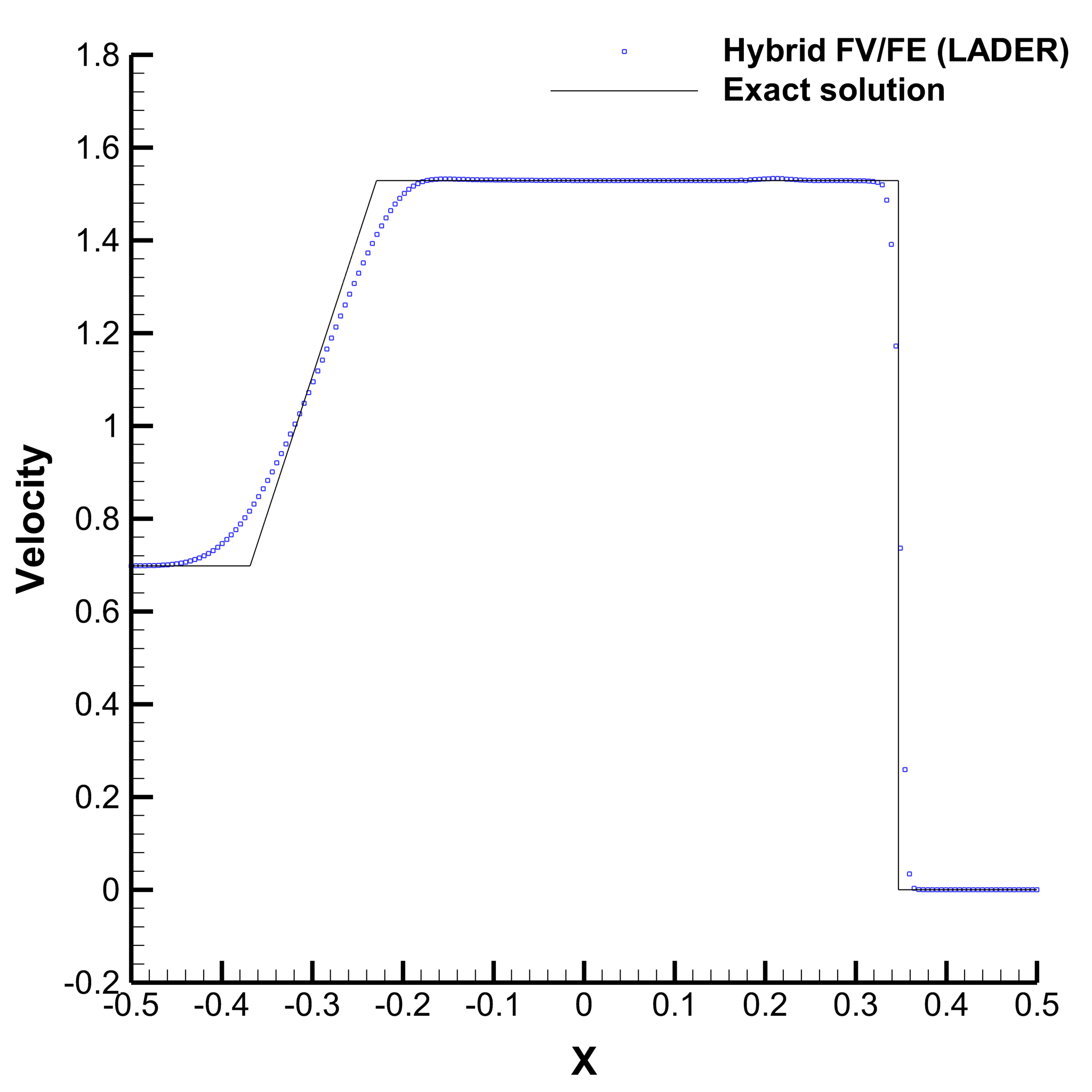}
	\includegraphics[trim= 5 5 5 5,clip,width=0.325\linewidth]{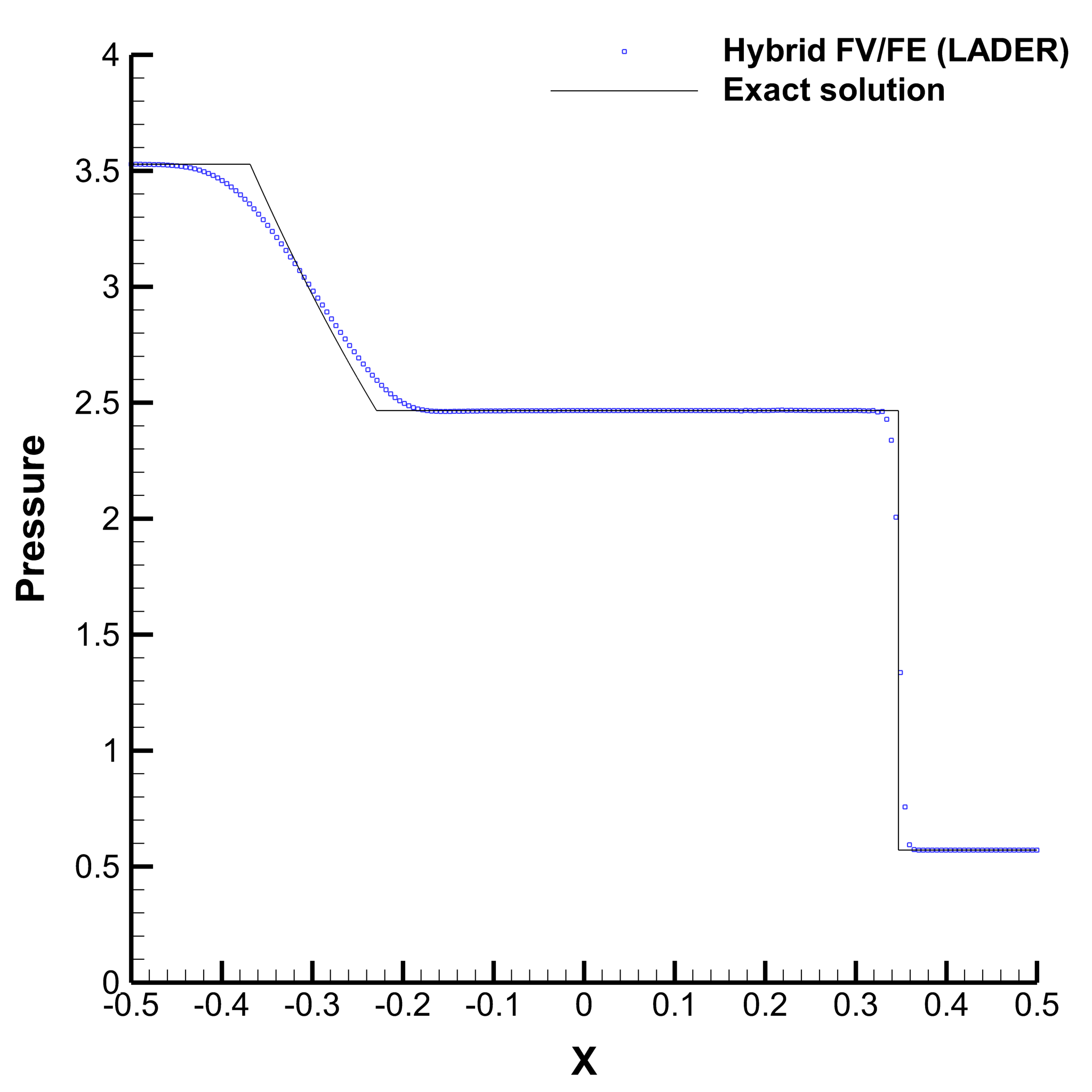}
	
	\caption{Riemann problem 3 (Lax). 1D cut through the numerical results along the line $y=0$ for $\rho$, $u$ and $p$ at $t_{\mathrm{end}}=0.14$ using the LADER-ENO method on mesh M1 ($\mathrm{CFL}_{c}=2.78$, $M\approx 0.94$).} 
	\label{fig:RP3_lader}
\end{figure}

The third test, RP3, corresponds to the Lax shock tube and is used to assess the ability of the method to deal with simple waves. The obtained results, presented in Figures \ref{fig:RP3_o1}-\ref{fig:RP3_lader}, match pretty well the exact reference solution.

\begin{figure}[!htbp]
	\centering
	\includegraphics[trim= 5 5 5 5,clip,width=0.325\linewidth]{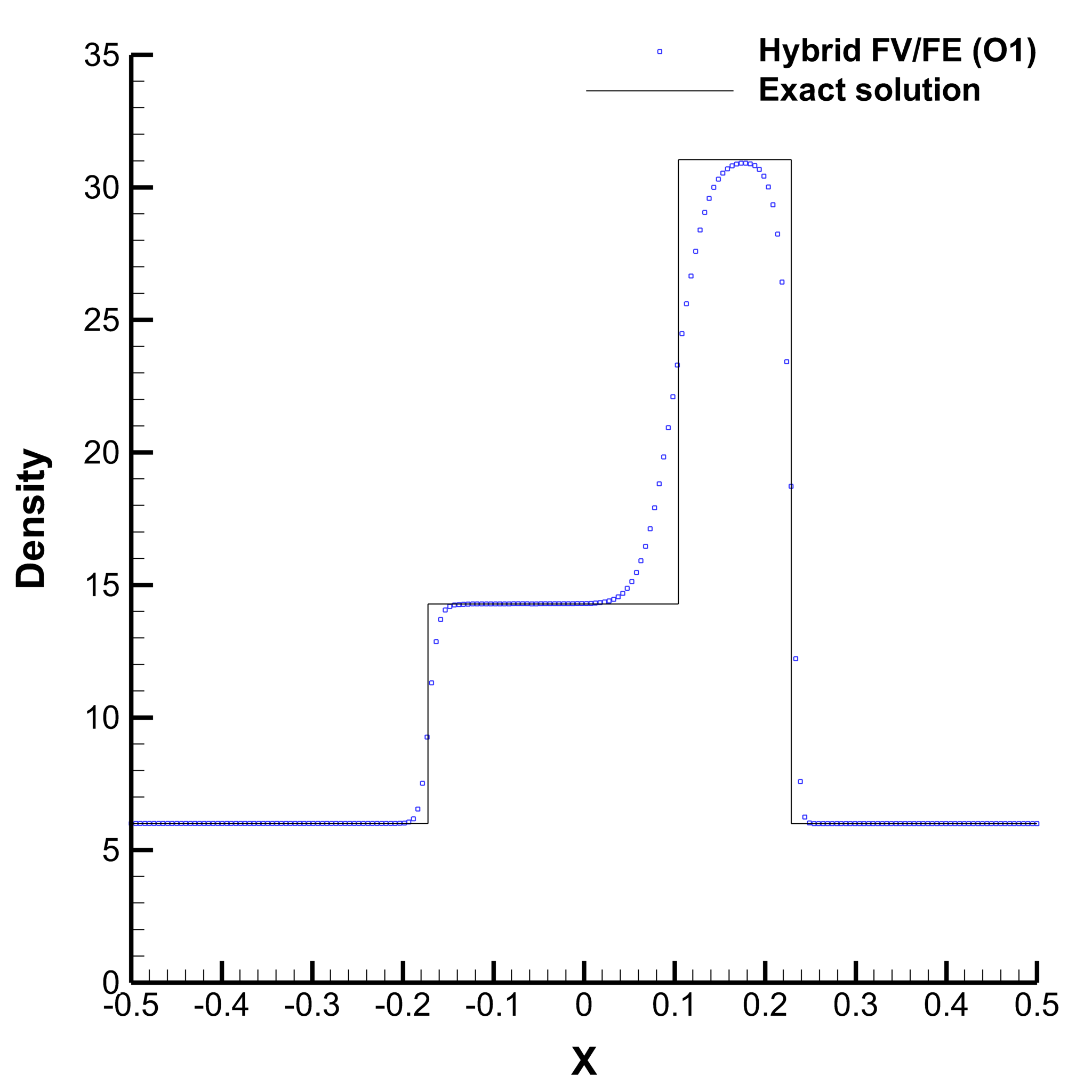}
	\includegraphics[trim= 5 5 5 5,clip,width=0.325\linewidth]{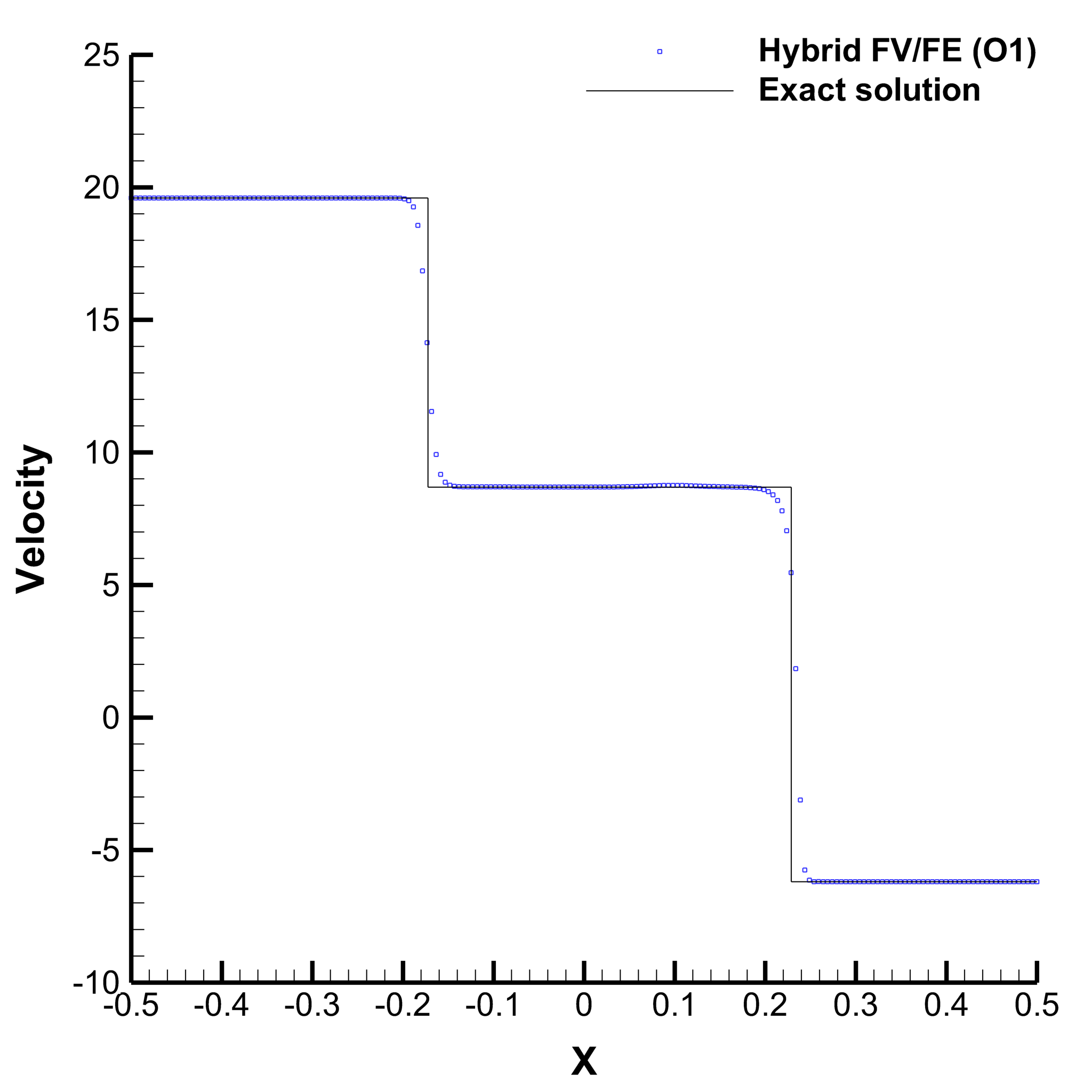}
	\includegraphics[trim= 5 5 5 5,clip,width=0.325\linewidth]{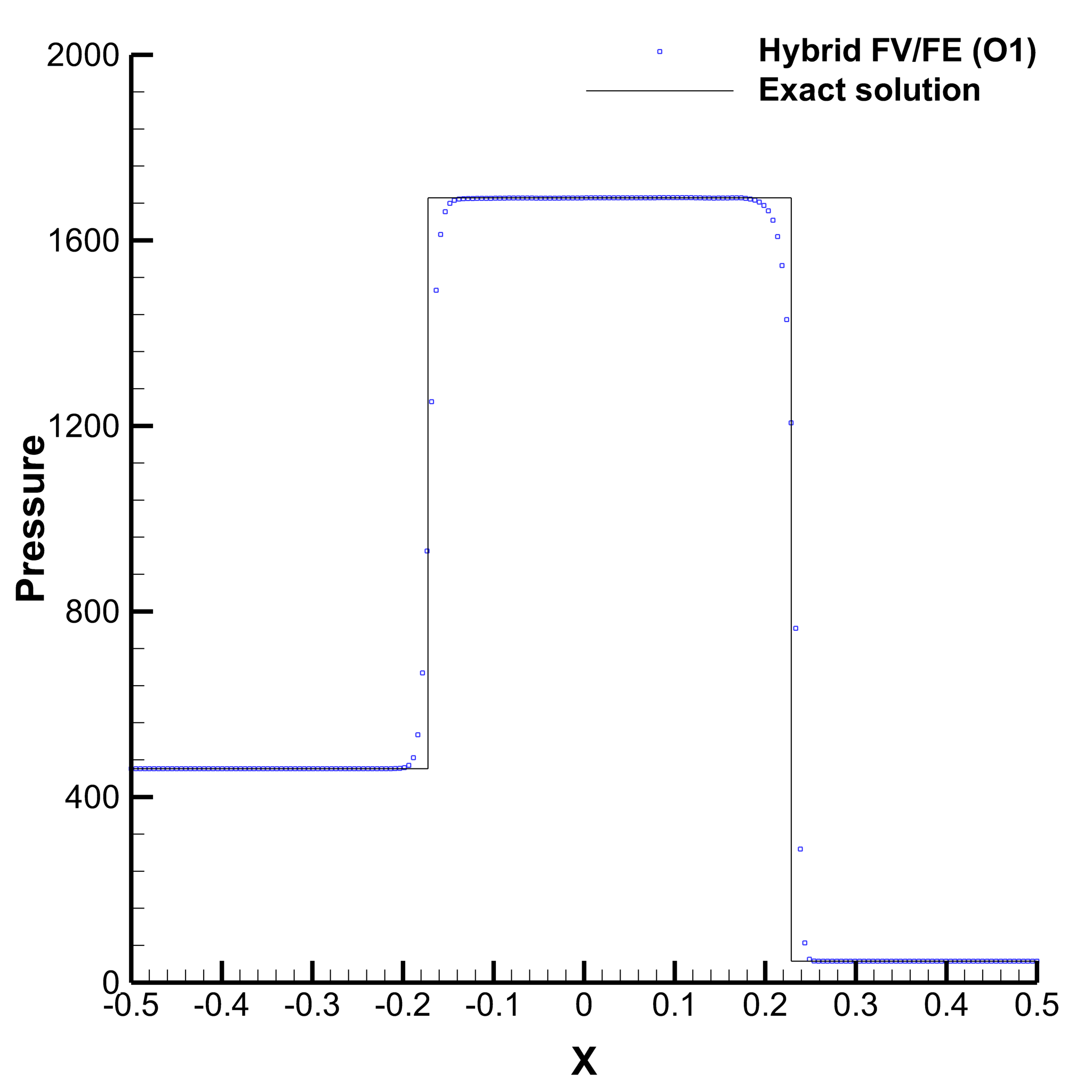}
	
	\caption{Riemann problem 4. 1D cut through the numerical results along the line $y=0$ for $\rho$, $u$ and $p$ at $t_{\mathrm{end}}=0.035$ using the first order method ($\mathrm{CFL}_{c}=0.42$,  $c_{\alpha}=5$, $M\approx 1.97$).}
	\label{fig:RP4_o1}
\end{figure}

\begin{figure}[!htbp]
	\centering
	\includegraphics[trim= 5 5 5 5,clip,width=0.325\linewidth]{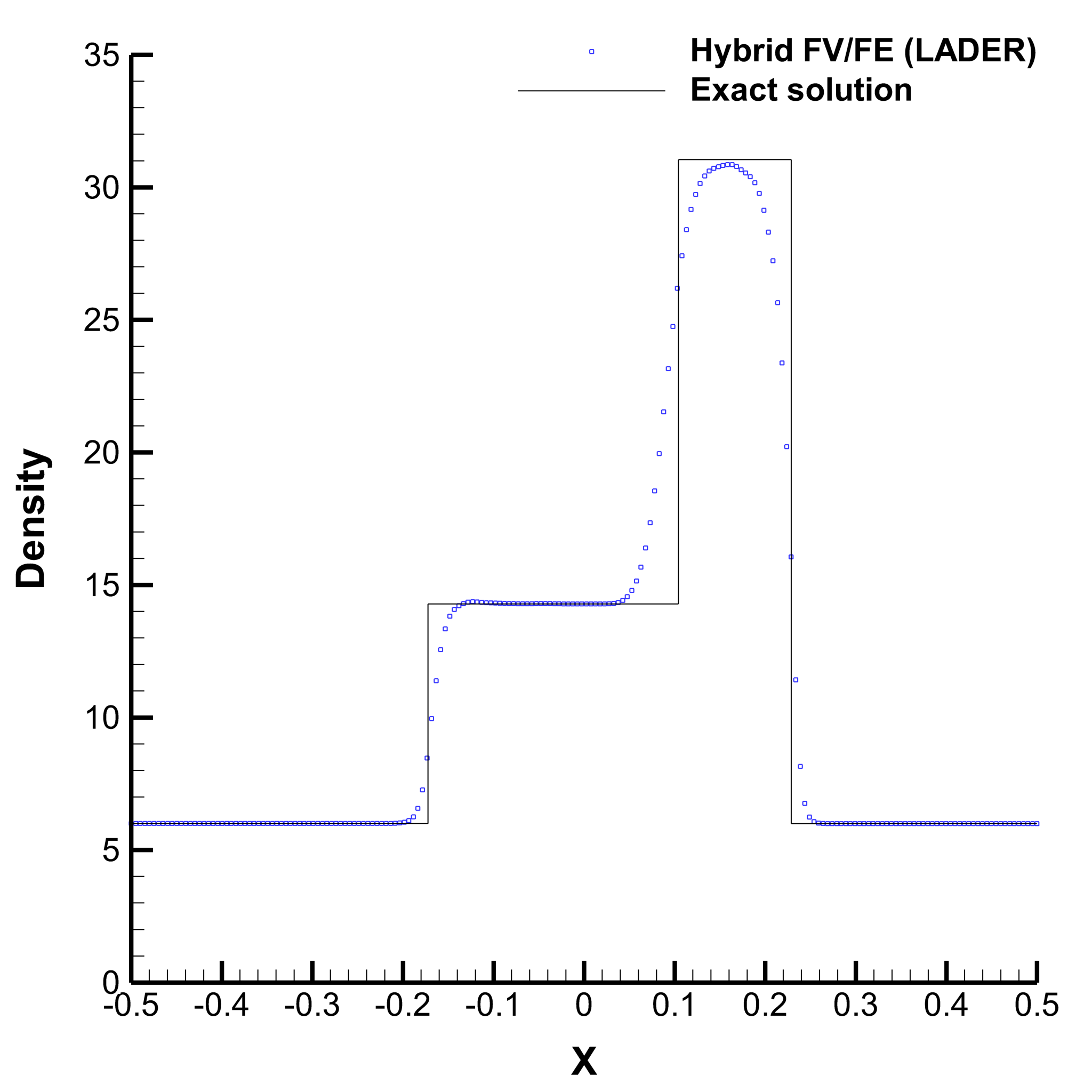}
	\includegraphics[trim= 5 5 5 5,clip,width=0.325\linewidth]{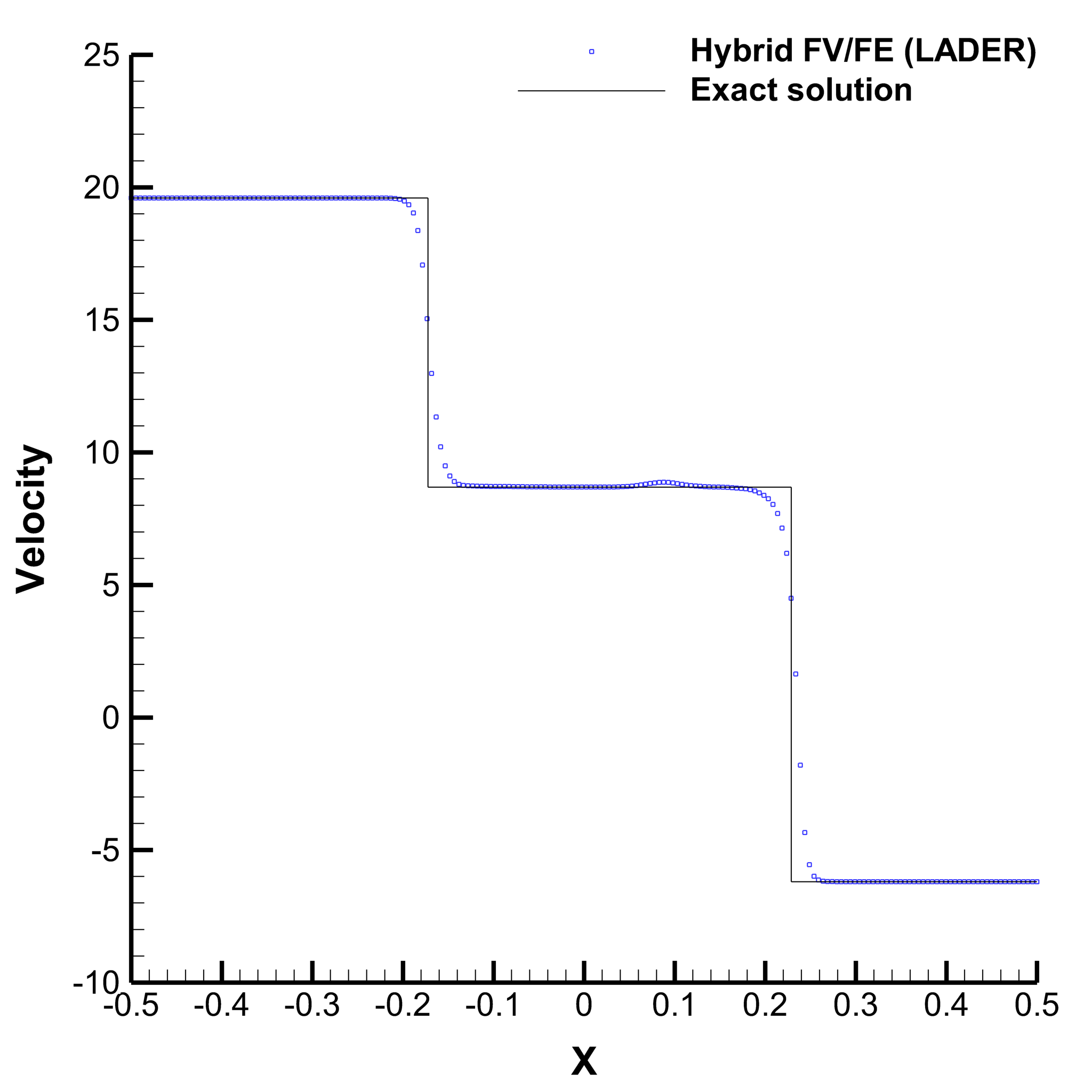}
	\includegraphics[trim= 5 5 5 5,clip,width=0.325\linewidth]{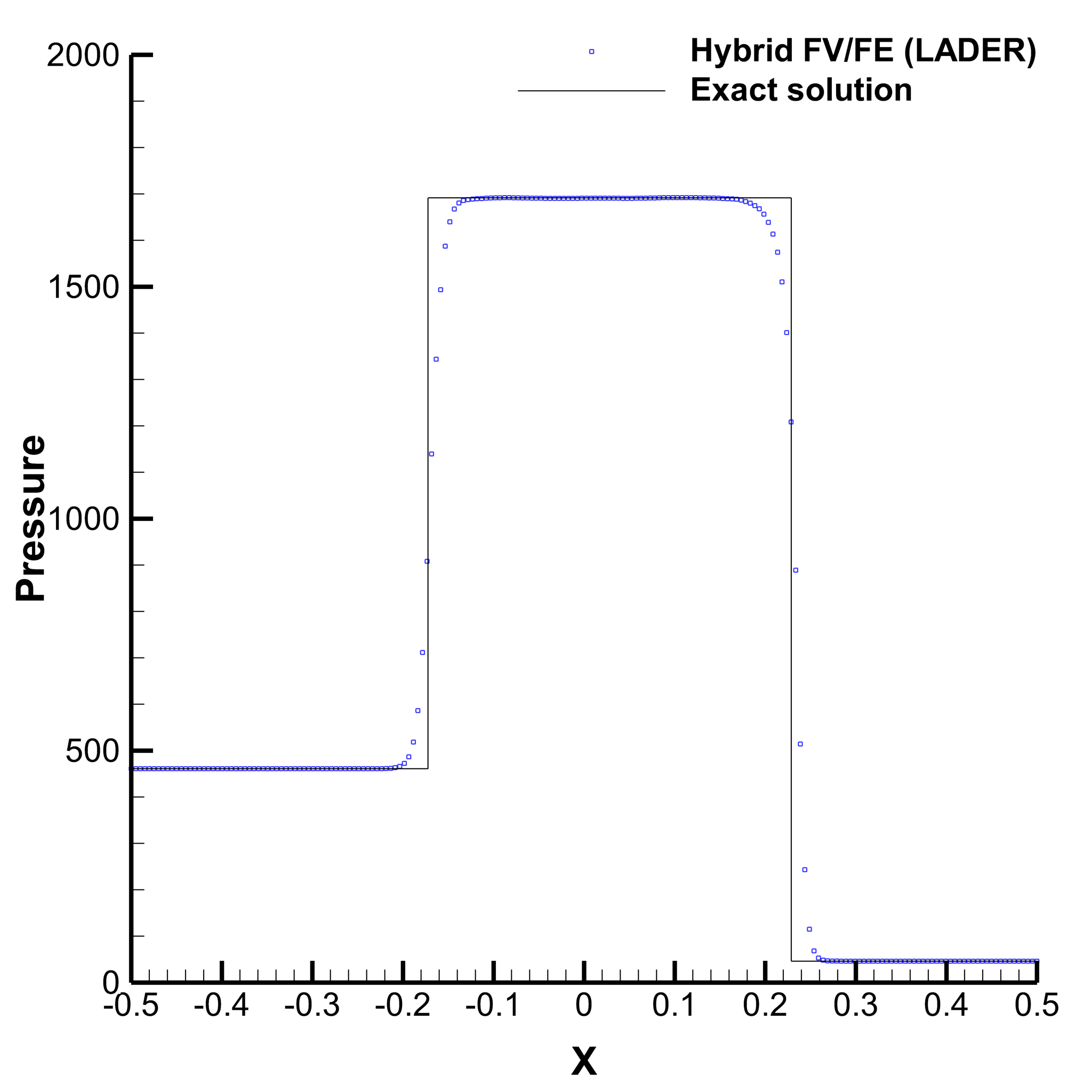}
	
	\caption{Riemann problem 4. 1D cut through the numerical results along the line $y=0$ for $\rho$, $u$ and $p$ at $t_{\mathrm{end}}=0.035$ using LADER-BJ method with reconstruction thorough primitive variables instead the conservative ones ($\mathrm{CFL}_{c}=0.42$,  $c_{\alpha}=5$, $M\approx 1.97$).}
	\label{fig:RP4_LBJ}
\end{figure}

The fourth Riemann problem, RP4, presents three strong discontinuities travelling to the right originated from two shock colliding waves. Figures \ref{fig:RP4_o1}-\ref{fig:RP4_LBJ} show the solution obtained with the first and second order schemes. Note that the highly restrictive Barth and Jespersen limiter has been employed jointly with an artificial viscosity coefficient, $c_{\alpha}=5$, to keep the stability of the scheme.

\begin{figure}[!htbp]
	\centering
	\includegraphics[trim= 5 5 5 5,clip,width=0.325\linewidth]{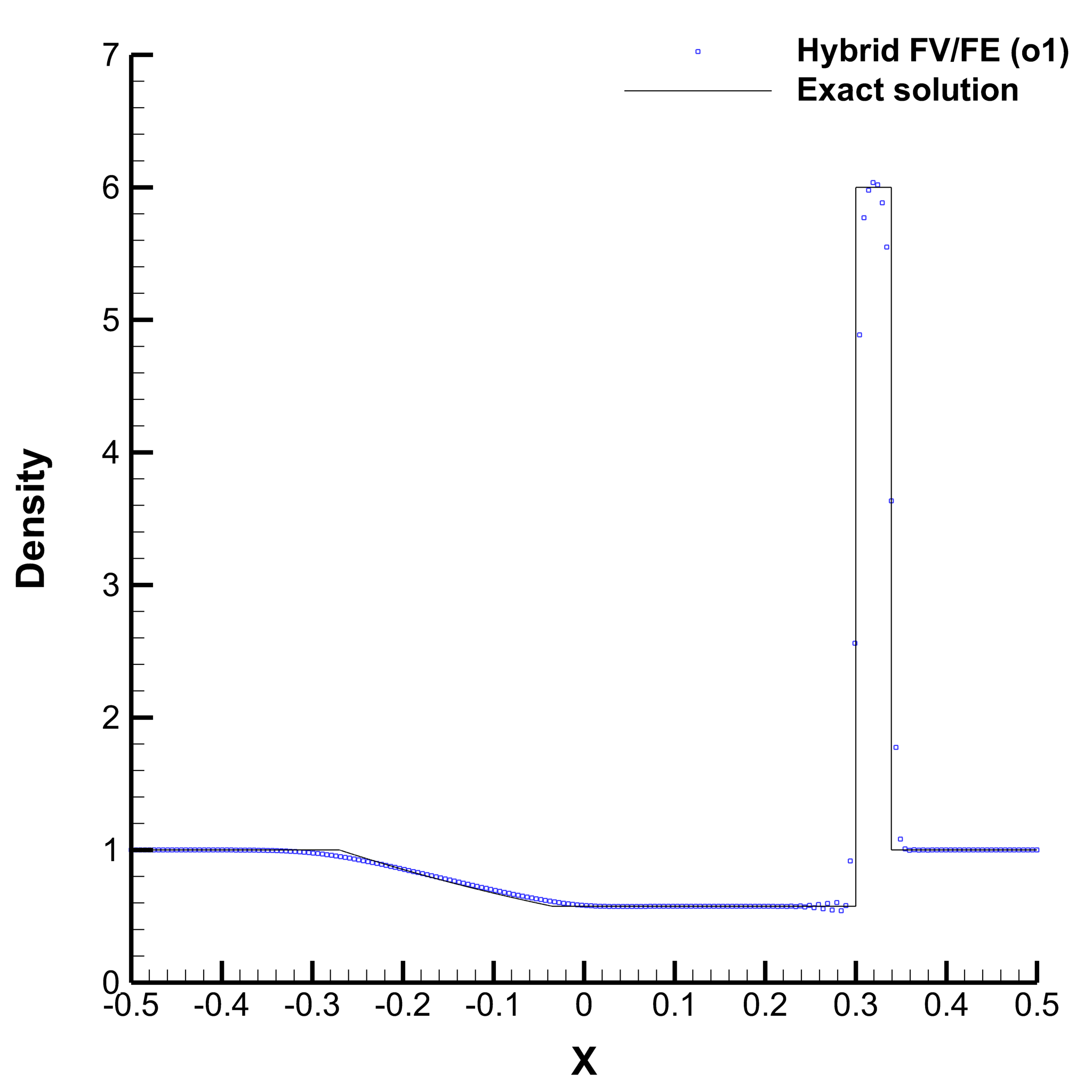}
	\includegraphics[trim= 5 5 5 5,clip,width=0.325\linewidth]{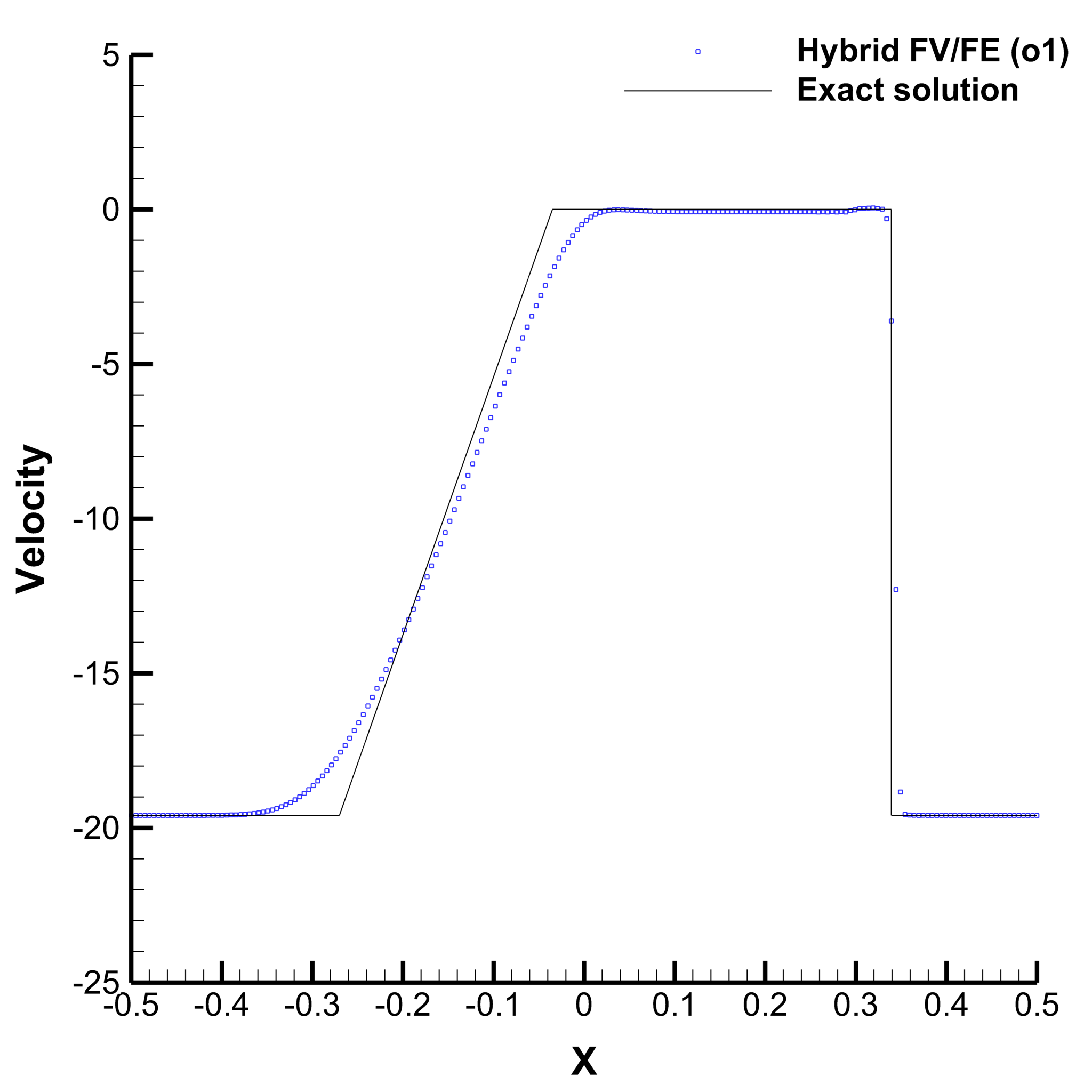}
	\includegraphics[trim= 5 5 5 5,clip,width=0.325\linewidth]{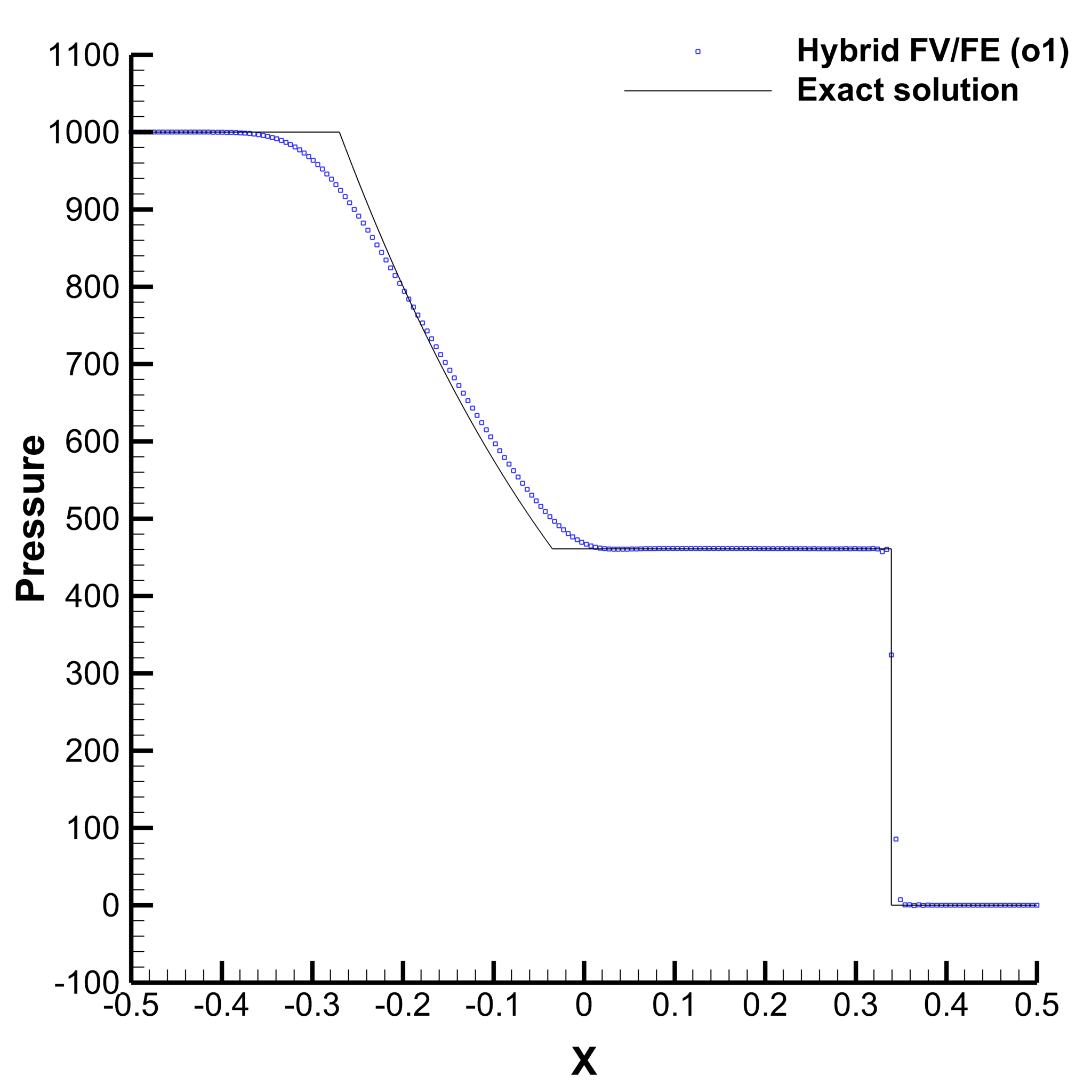}
	
	\caption{Riemann problem 5. 1D cut through the numerical results along the line $y=0$ for $\rho$, $u$ and $p$ at $t_{\mathrm{end}}=0.01$ using the first order scheme ($\mathrm{CFL}_{c}=2.7\cdot 10^{-3}$, $c_{\alpha}=2$, $M\approx 956.42$).} 
	\label{fig:RP5_o1}
\end{figure}

\begin{figure}[!htbp]
	\centering
	\includegraphics[trim= 5 5 5 5,clip,width=0.325\linewidth]{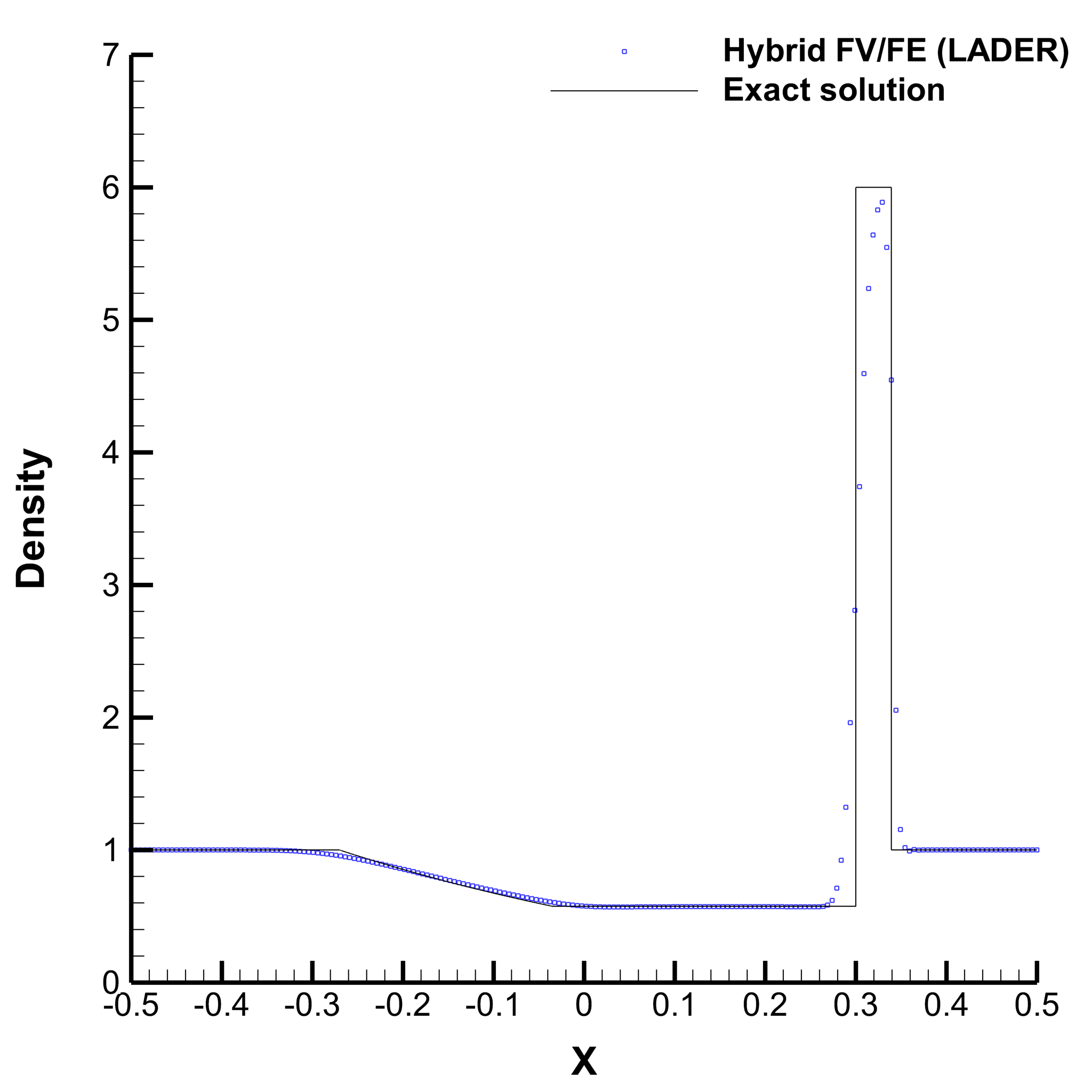}
	\includegraphics[trim= 5 5 5 5,clip,width=0.325\linewidth]{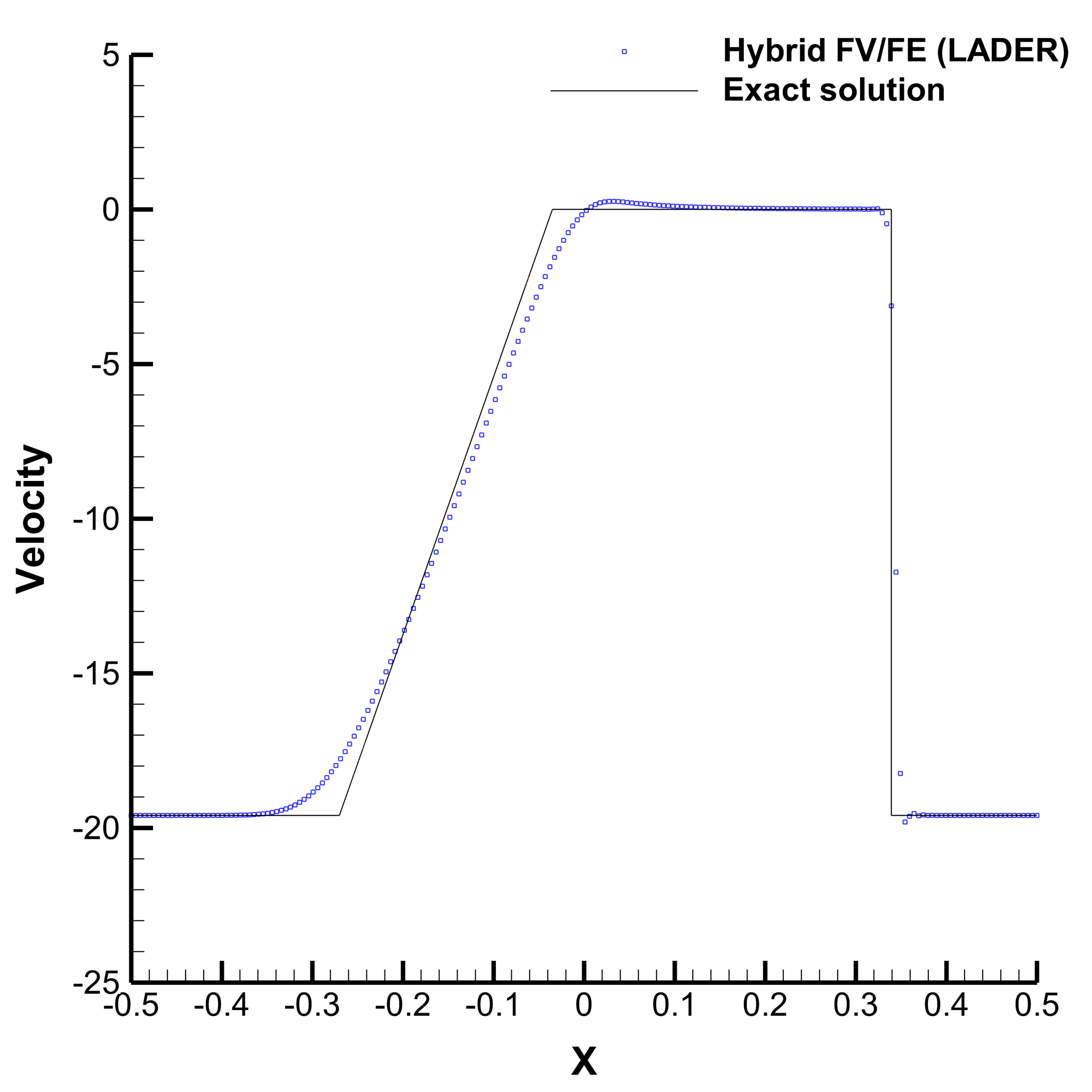}
	\includegraphics[trim= 5 5 5 5,clip,width=0.325\linewidth]{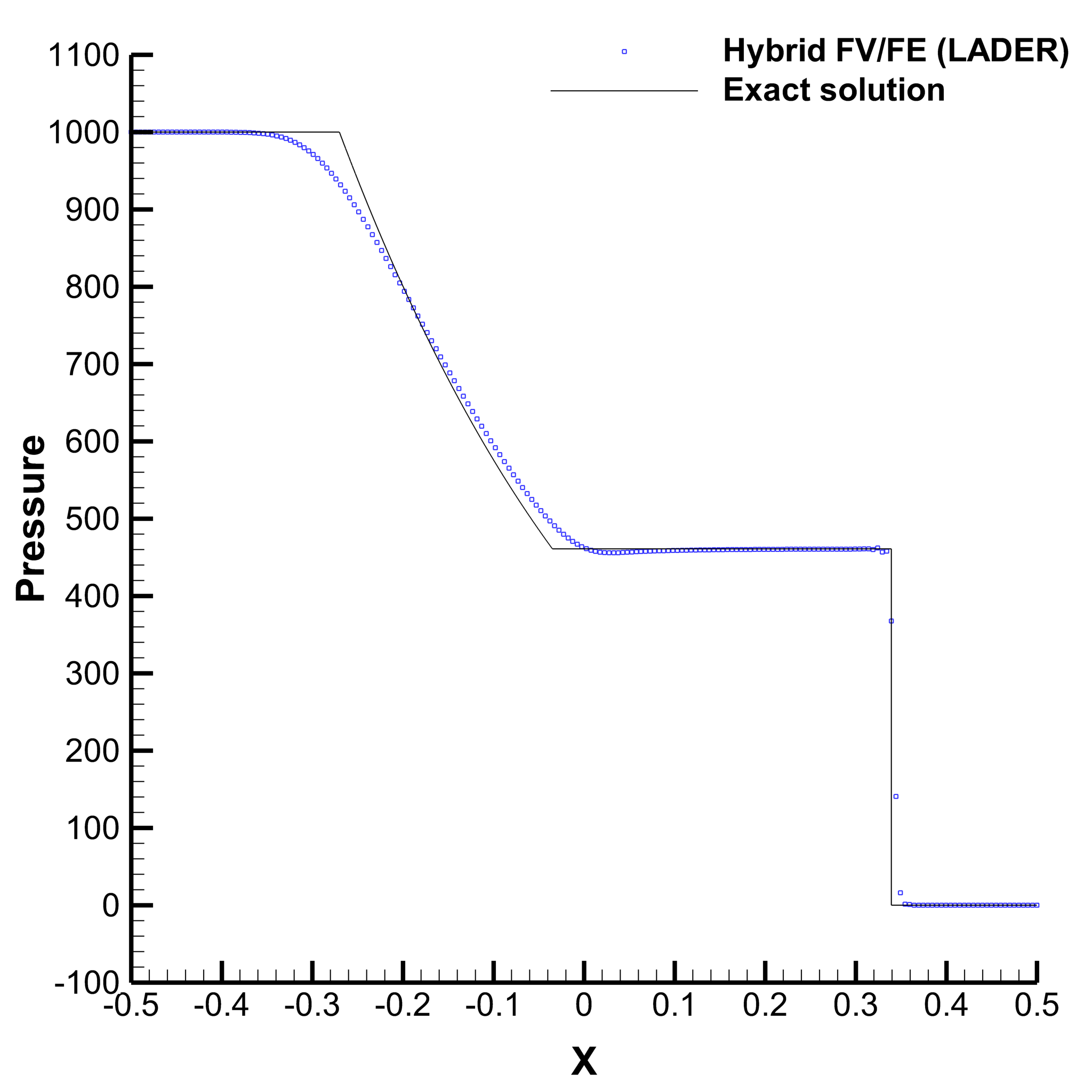}
	
	\caption{Riemann problem 5.1D cut through the numerical results along the line $y=0$ for $\rho$, $u$ and $p$ at $t_{\mathrm{end}}=0.01$ using the LADER-ENO method  ($\mathrm{CFL}_{c}=2.7\cdot 10^{-3}$, $c_{\alpha}=2$, $M\approx 513.68$).}
	\label{fig:RP5_lader1}
\end{figure}

\begin{figure}[!htbp]
	\centering
	\includegraphics[trim= 5 5 5 5,clip,width=0.325\linewidth]{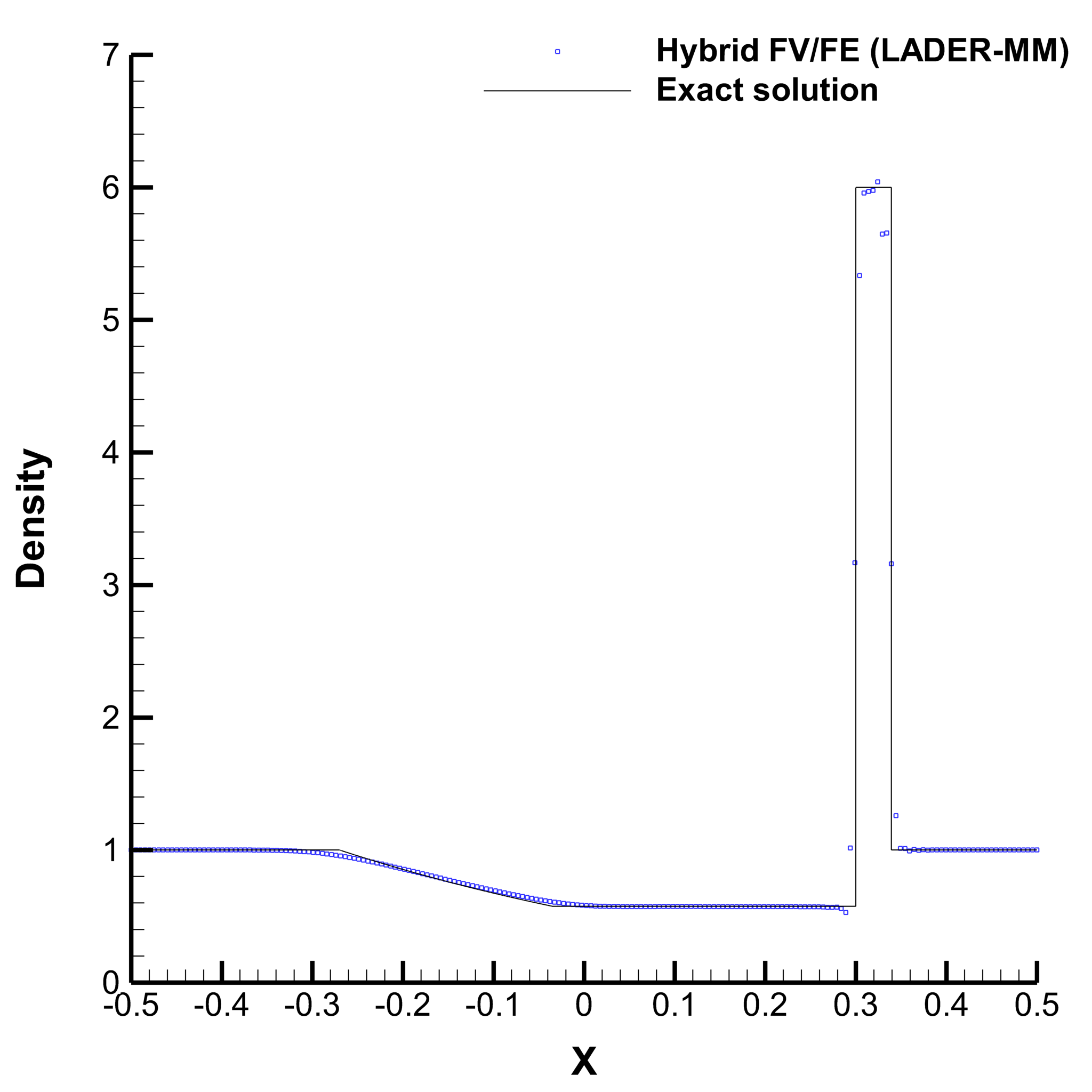}
	\includegraphics[trim= 5 5 5 5,clip,width=0.325\linewidth]{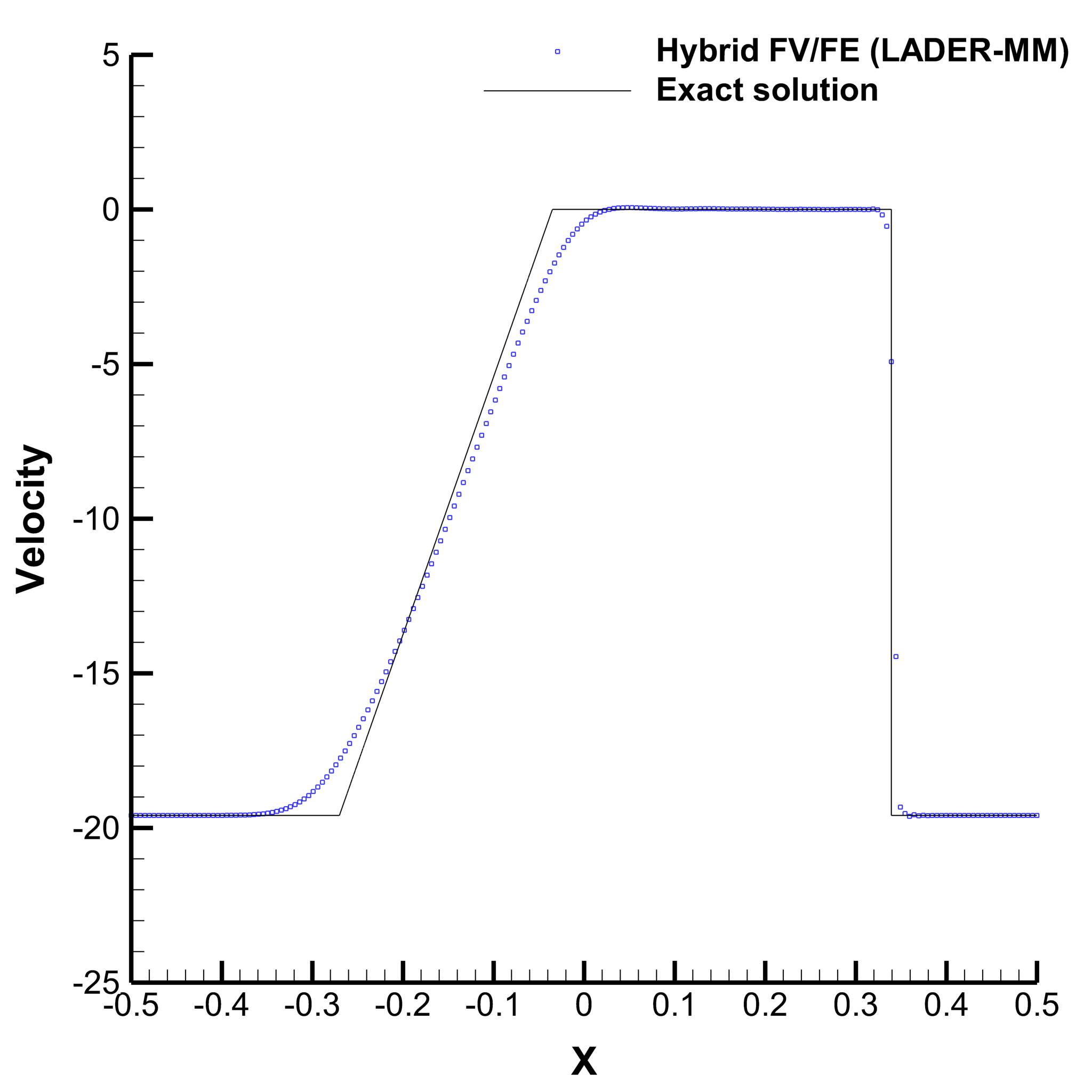}
	\includegraphics[trim= 5 5 5 5,clip,width=0.325\linewidth]{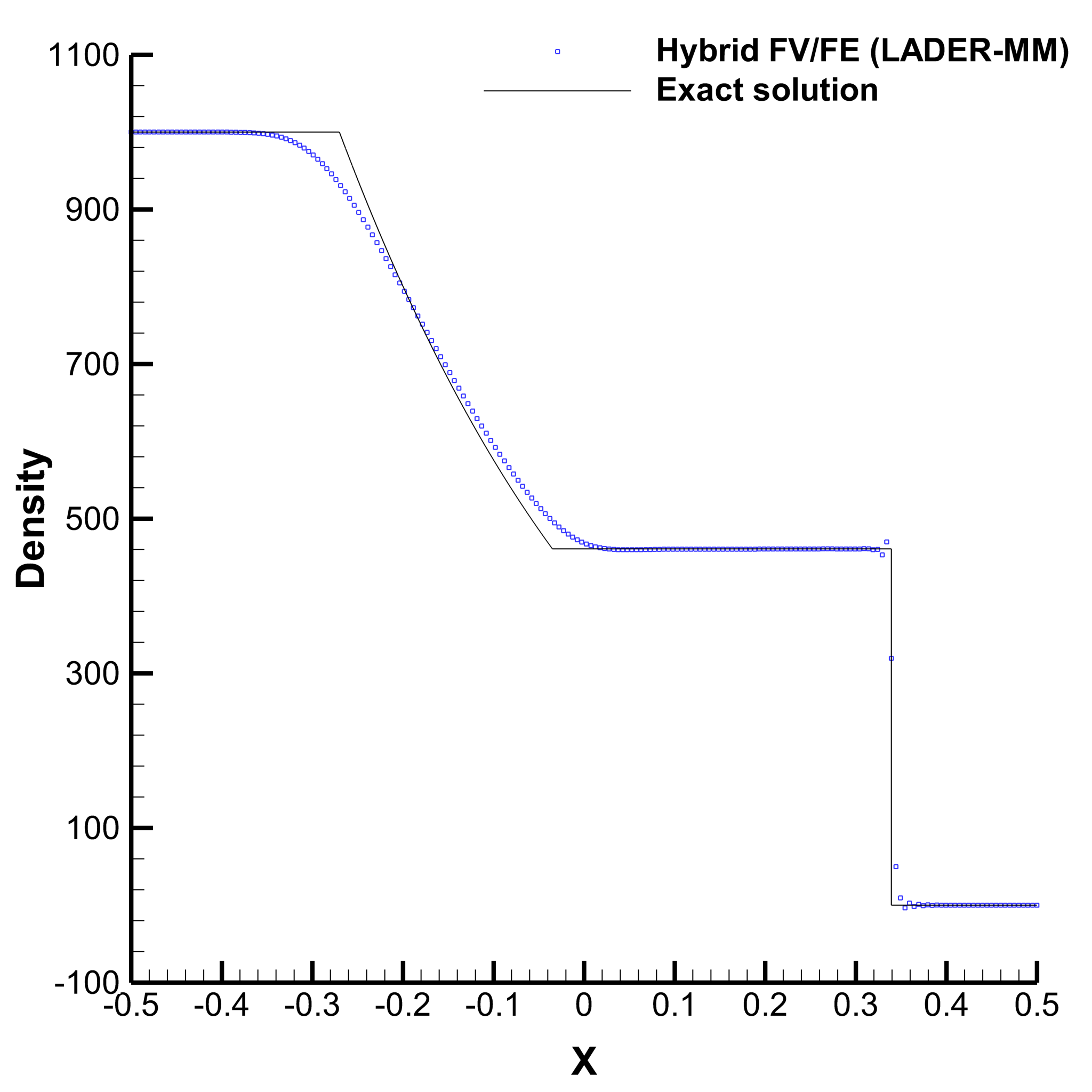}
	
	\caption{Riemann problem 5. 1D cut through the numerical results along the line $y=0$ for $\rho$, $u$ and $p$ at $t_{\mathrm{end}}=0.01$ using the LADER-minmod method ($\mathrm{CFL}_{c}=2.7\cdot 10^{-3}$, $c_{\alpha}=2$, $M\approx 691.58$).}
	\label{fig:RP5_laderMM}
\end{figure}

RP5 is a severe test defined as a modification of the left half of the blast problem introduced in \cite{WC84}. It accounts for a left rarefaction wave, a right-travelling shock wave and a stationary contact discontinuity generated by an initial large pressure jump of order $10^{5}$ and a small velocity variation. The second order scheme has been run using two different limiter strategies. We observe that the minmod limiter, Figure \ref{fig:RP5_laderMM}, is more severe on damping the oscillation appearing after the rarefaction wave on the velocity field than the ENO-based reconstruction, Figure \ref{fig:RP5_lader1}, which captures better the right shock. The results obtained with the first order scheme are reported in Figure \ref{fig:RP5_o1}. The robustness of the developed methodology and its capability to deal with slowly moving contact discontinuities, for really high Mach numbers, are clearly proven.

\begin{figure}[!htbp]
	\centering
	\includegraphics[trim= 5 5 5 5,clip,width=0.325\linewidth]{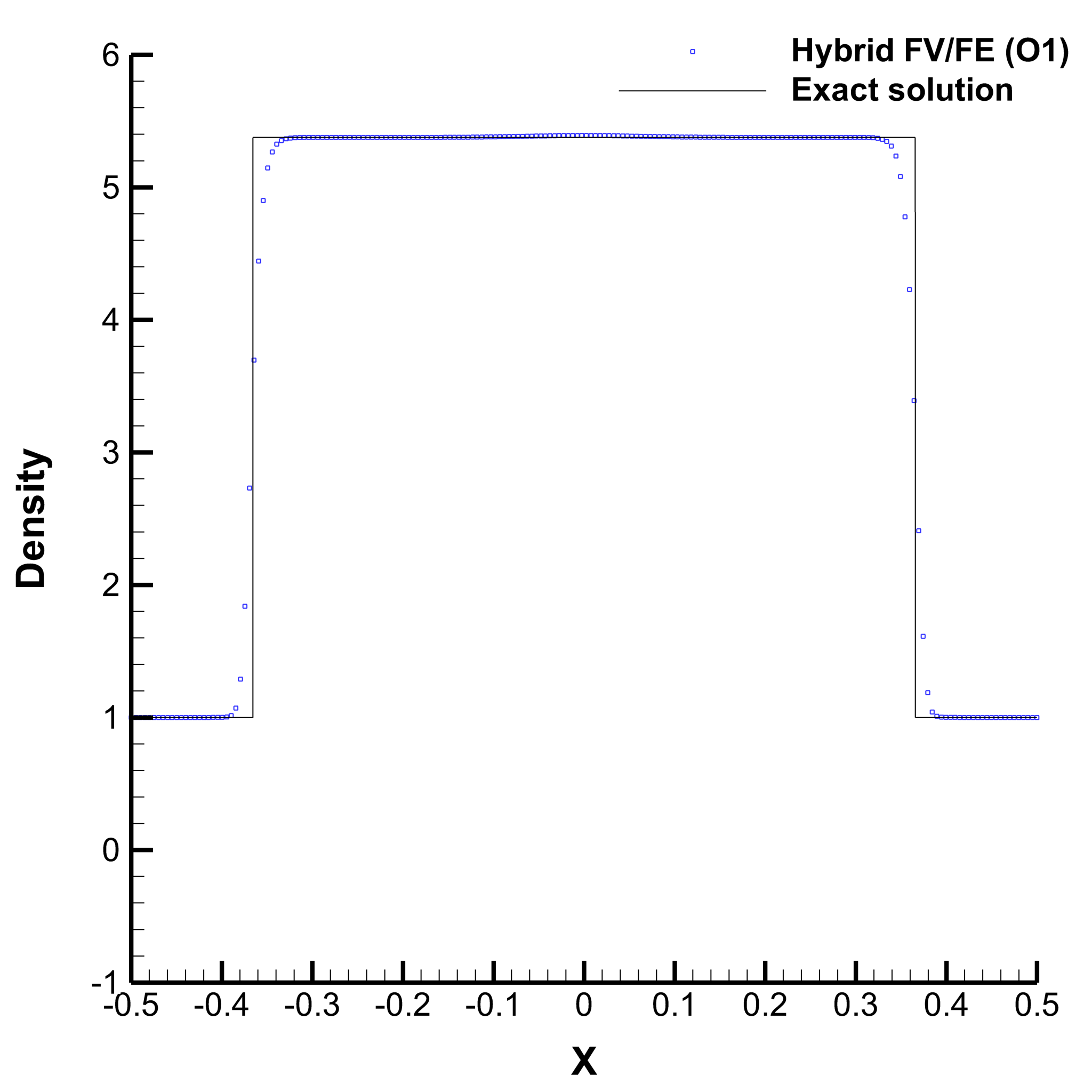}
	\includegraphics[trim= 5 5 5 5,clip,width=0.325\linewidth]{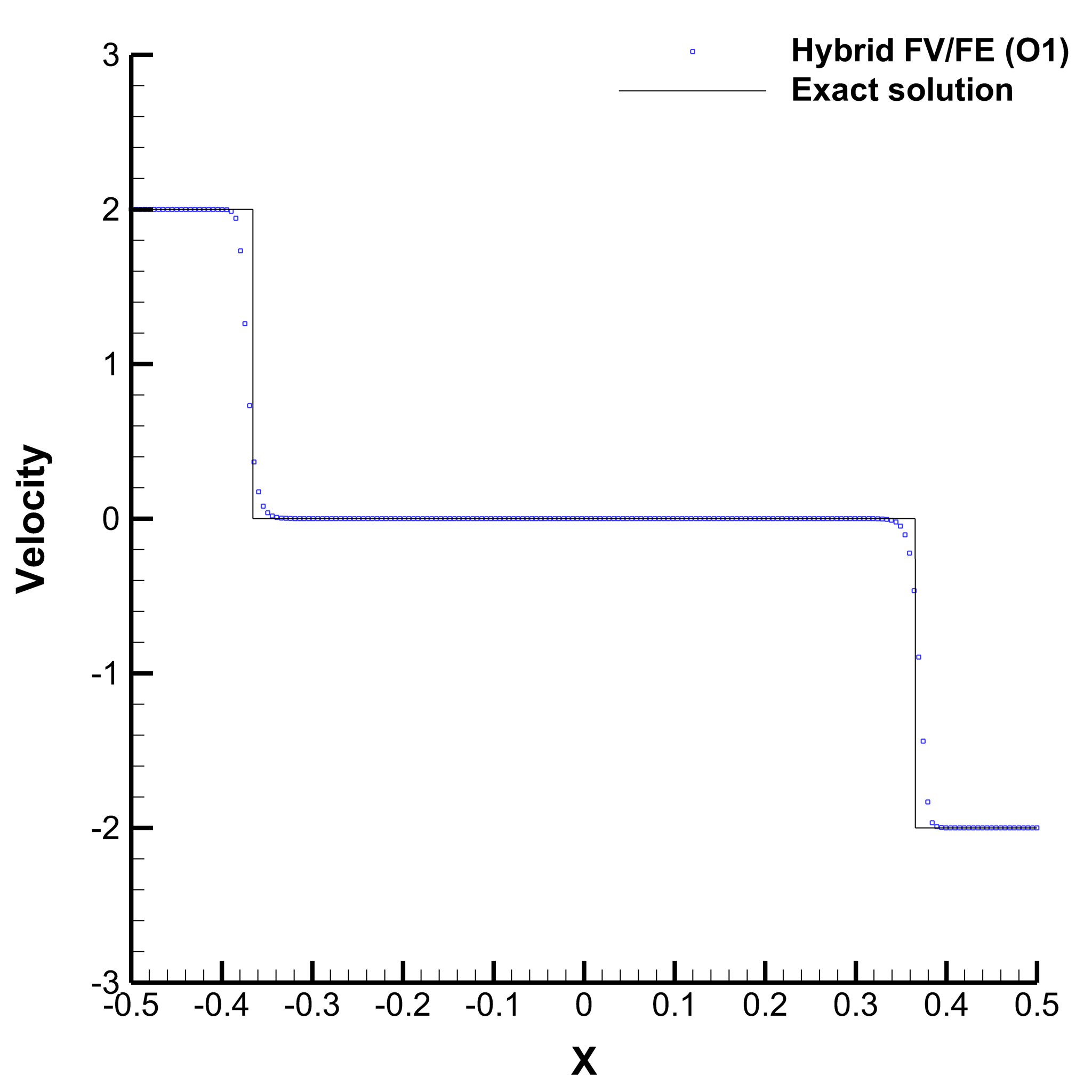}
	\includegraphics[trim= 5 5 5 5,clip,width=0.325\linewidth]{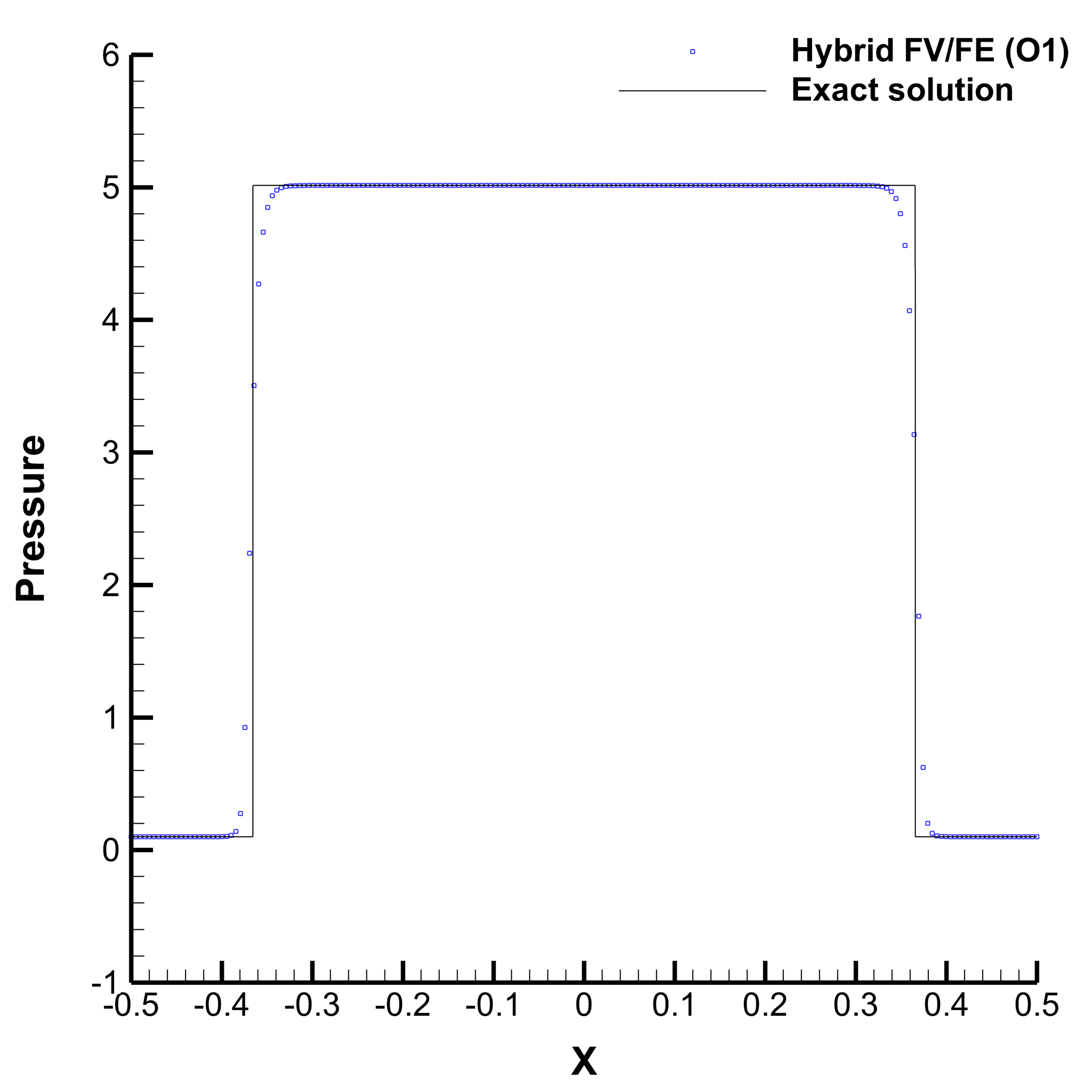}
	
	\caption{Riemann problem 6. 1D cut through the numerical results along the line $y=0$ for $\rho$, $u$ and $p$ at $t_{\mathrm{end}}=0.8$ using the first order method ($\mathrm{CFL}_{c}=9.35$, $c_{\alpha}=2$, $M\approx 7.75$).}
	\label{fig:RP6_o1}
\end{figure}

\begin{figure}[!htbp]
	\centering
	\includegraphics[trim= 5 5 5 5,clip,width=0.325\linewidth]{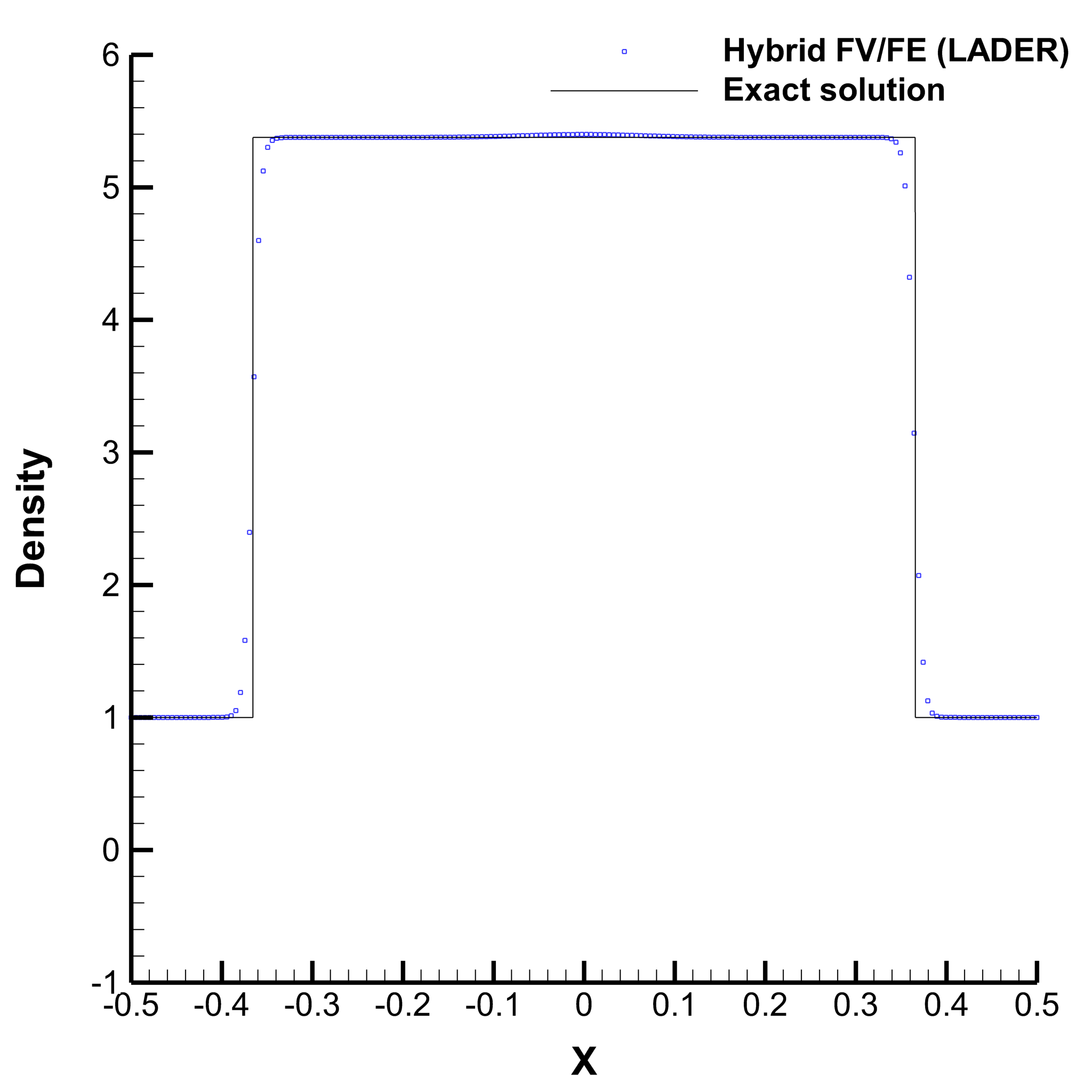}
	\includegraphics[trim= 5 5 5 5,clip,width=0.325\linewidth]{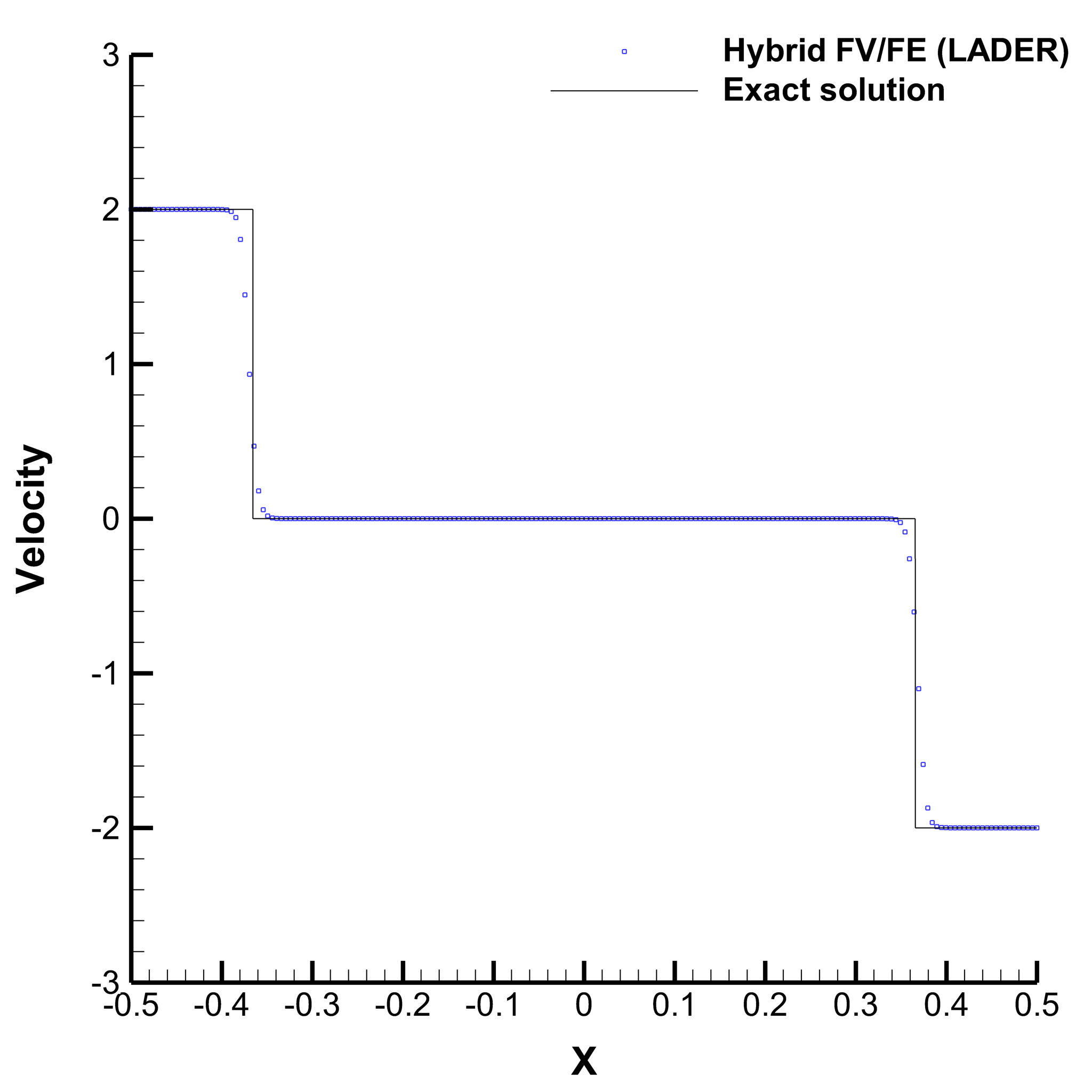}
	\includegraphics[trim= 5 5 5 5,clip,width=0.325\linewidth]{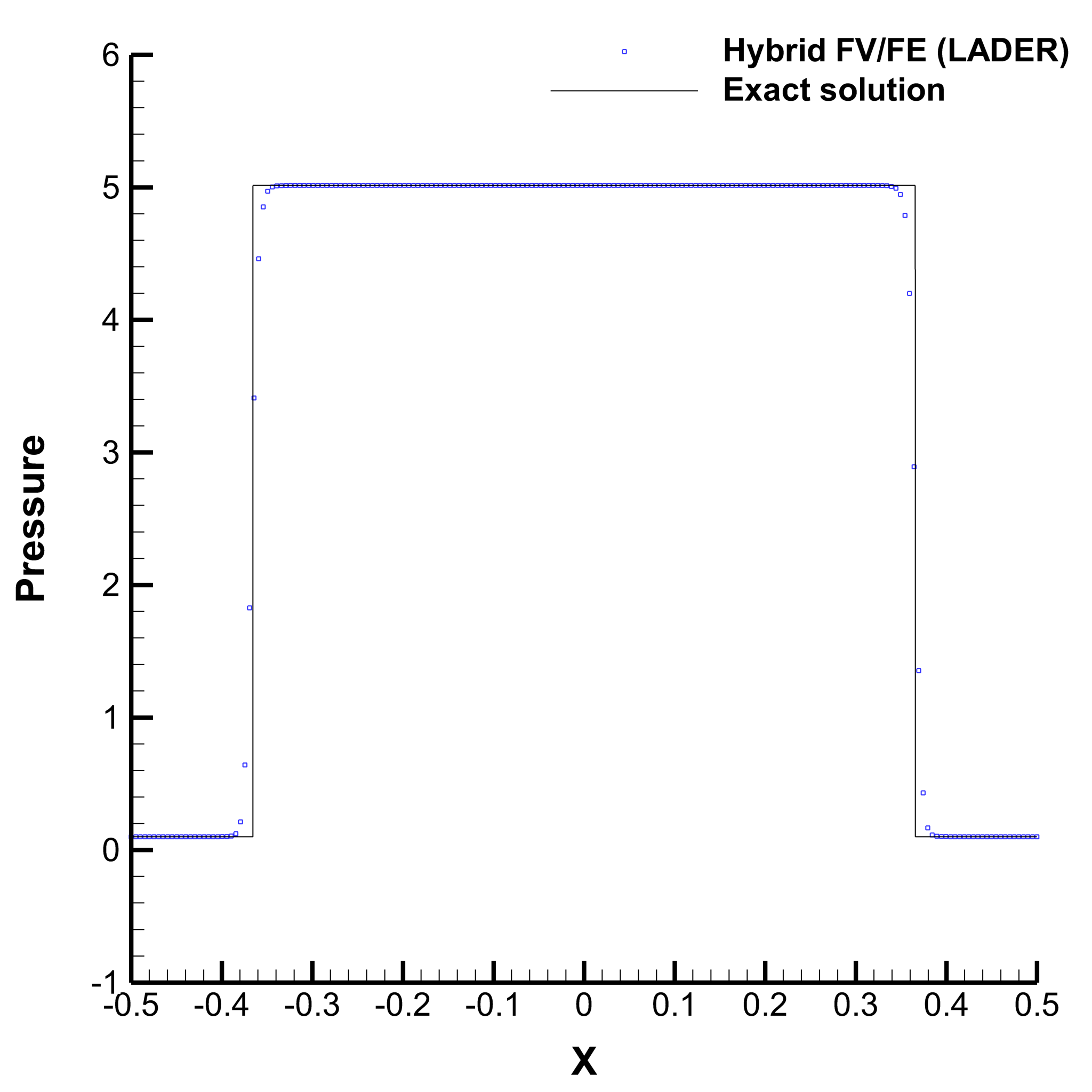}
	
	\caption{Riemann problem 6. 1D cut through the numerical results along the line $y=0$ for $\rho$, $u$ and $p$ at $t_{\mathrm{end}}=0.8$ using the first order method ($\mathrm{CFL}_{c}=9.35$, $c_{\alpha}=2$, $M\approx 7.29$).}
	\label{fig:RP6_lader1}
\end{figure}

The last Riemann problem considered, RP6, is characterised by two shock waves travelling in opposite directions. An excellent agreement with the exact solution is observed in Figures \ref{fig:RP6_o1}-\ref{fig:RP6_lader1}.

\subsection{2D circular explosion} \label{sec:2DCE}
\begin{figure}[!htbp]
	\centering
	\includegraphics[trim= 10 5 5 5,clip,width=0.45\linewidth]{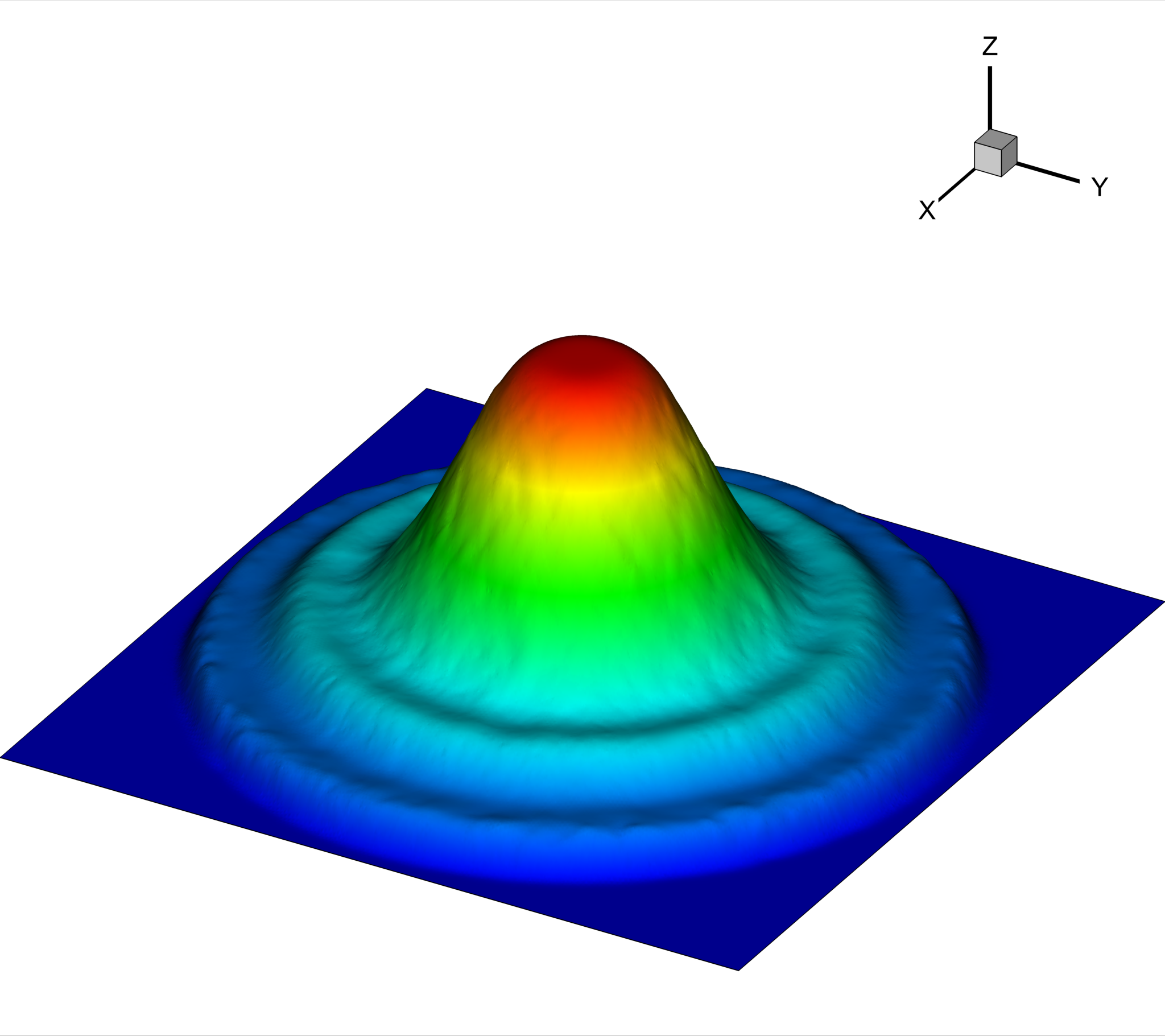}\hspace{0.05\linewidth}
	\includegraphics[trim= 10 5 5 5,clip,width=0.45\linewidth]{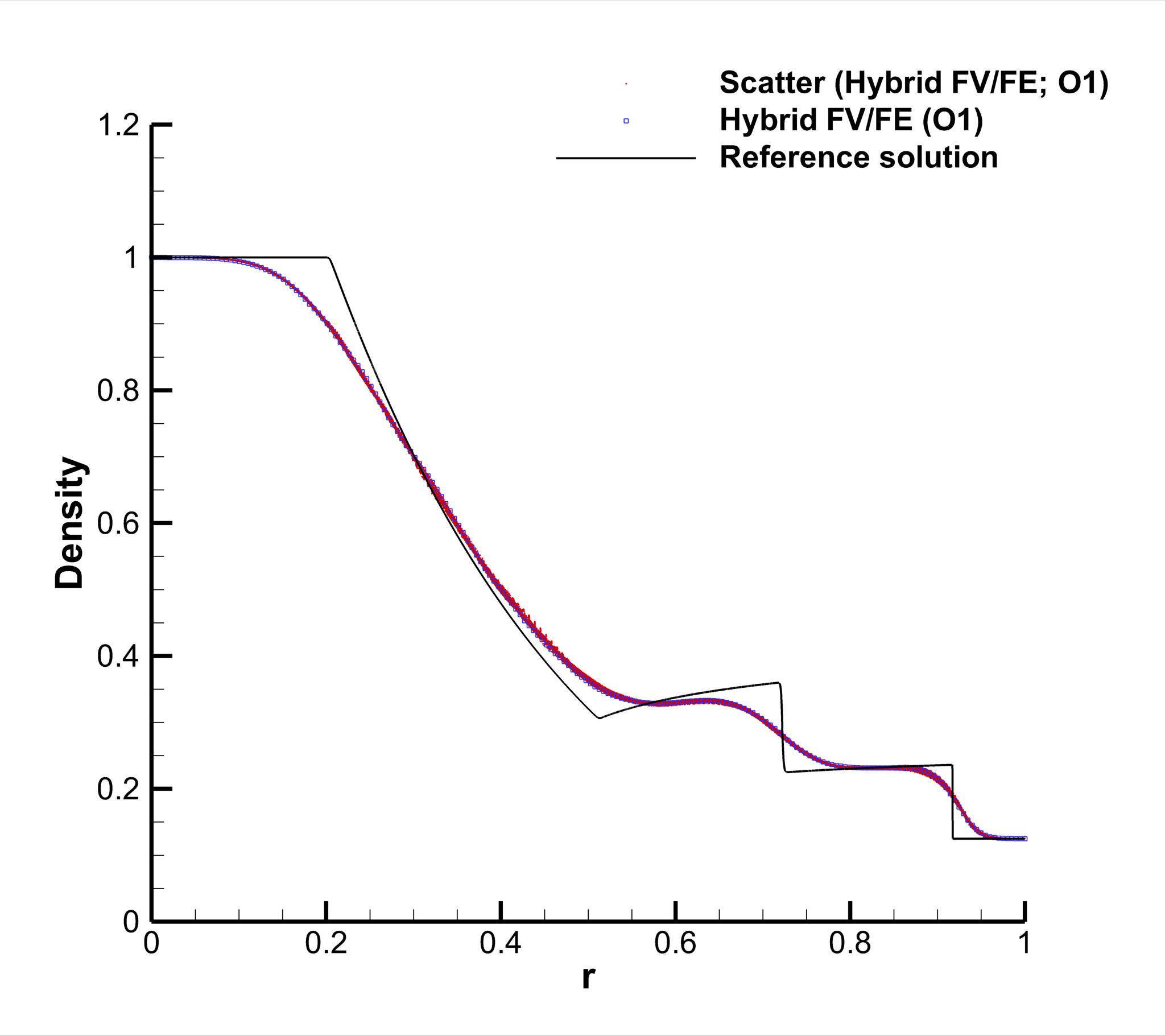}
	
	\vspace{0.05\linewidth}
	\includegraphics[trim= 10 5 5 5,clip,width=0.45\linewidth]{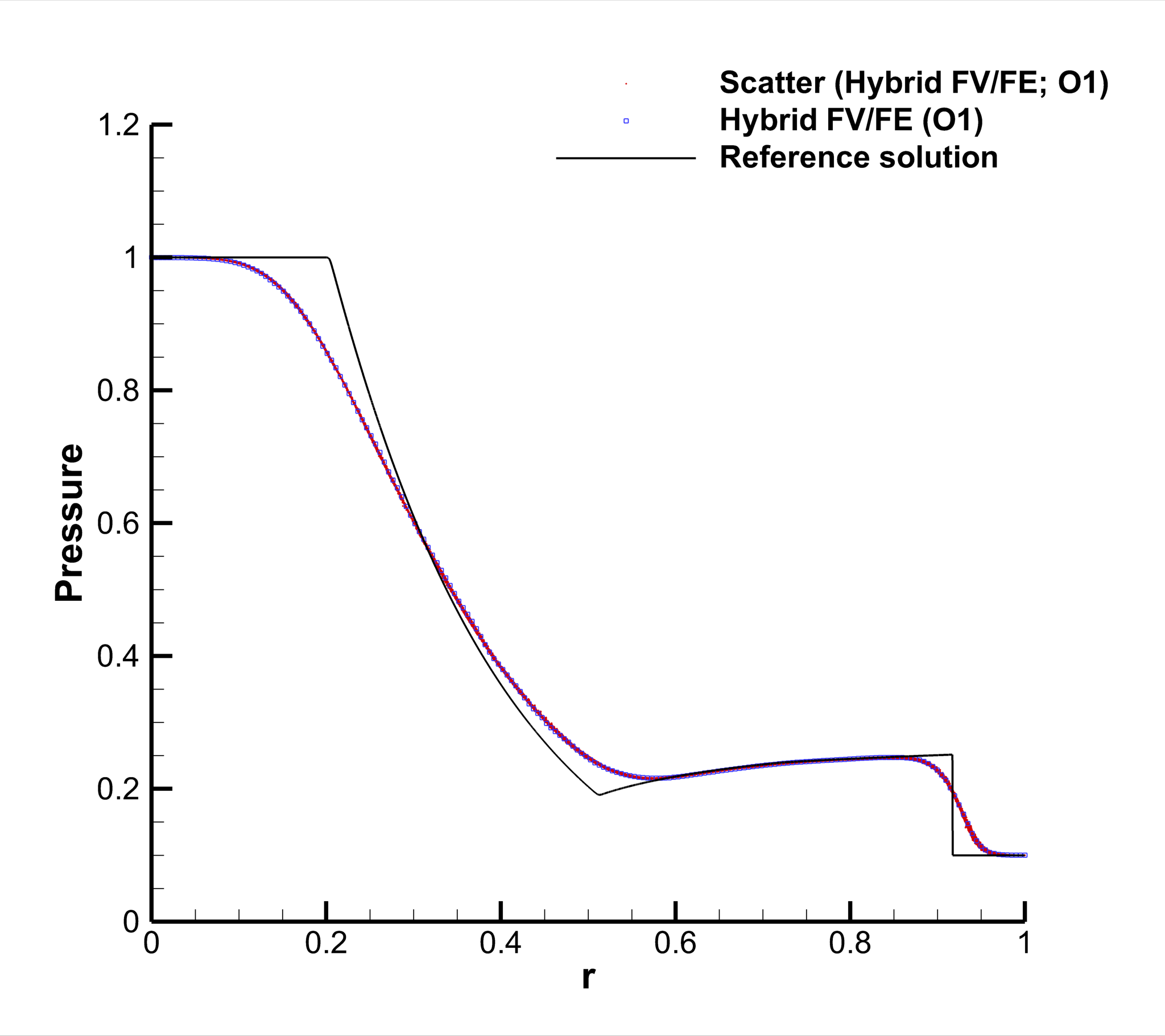}\hspace{0.05\linewidth}
	\includegraphics[trim= 10 5 5 5,clip,width=0.45\linewidth]{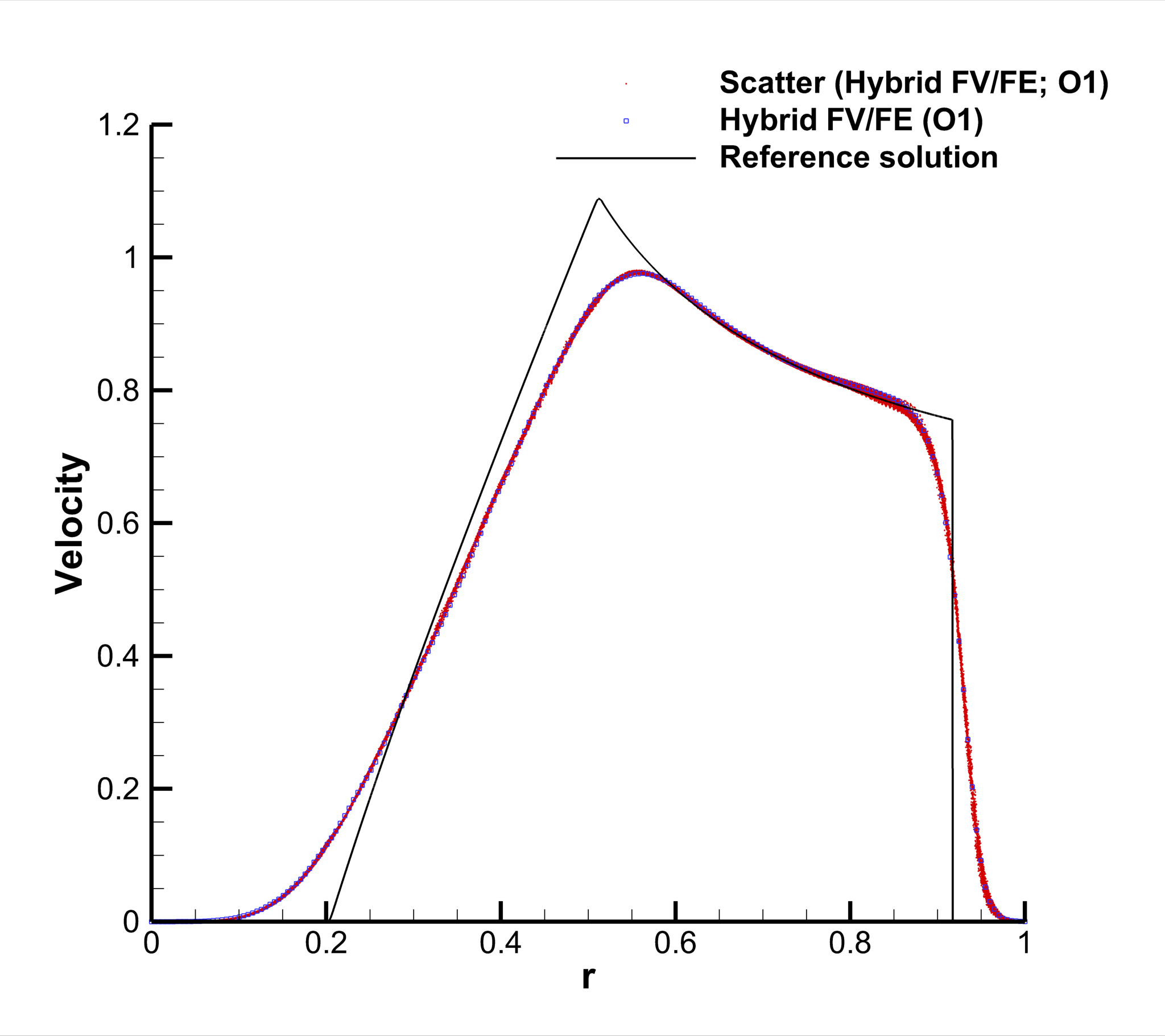}
	\caption{Circular explosion. The left top image corresponds to the 3D plot of the obtained $\rho$ at the final time whereas the 1D plots containing the reference solution (black continuous line), a 1D cut (blue squared line) and the scatter plot (red dots), correspond to the $\rho$, $\press$, and $\left|\vel\right|$ fields obtained using the first order scheme (\textcolor{black}{$c_{\alpha }=1$}).} 
	\label{CE85_o1_t025}
\end{figure}

\begin{figure}[!htbp]
	\centering
	\includegraphics[trim= 10 5 5 5,clip,width=0.45\linewidth]{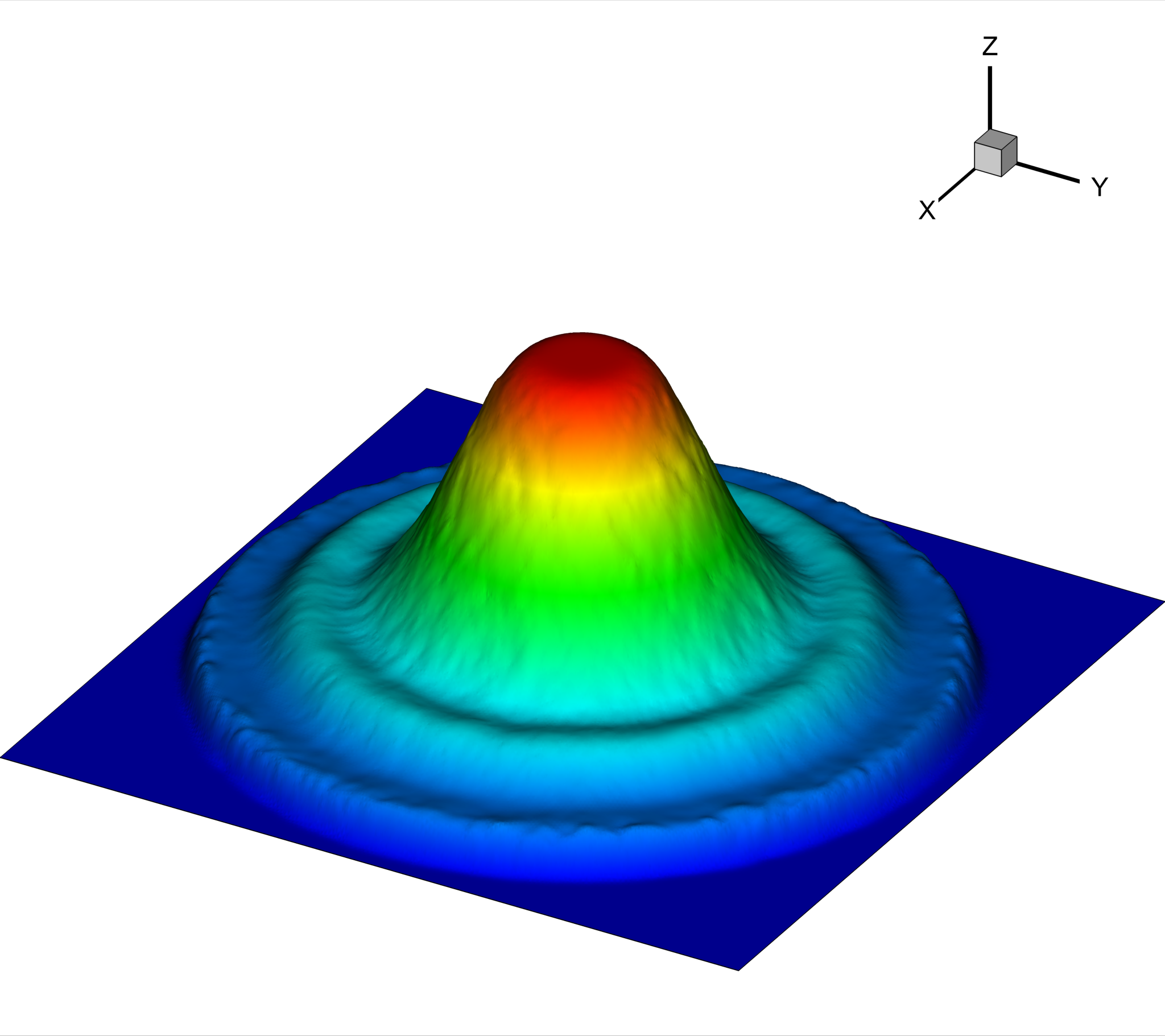}\hspace{0.05\linewidth}
	\includegraphics[trim= 10 5 5 5,clip,width=0.45\linewidth]{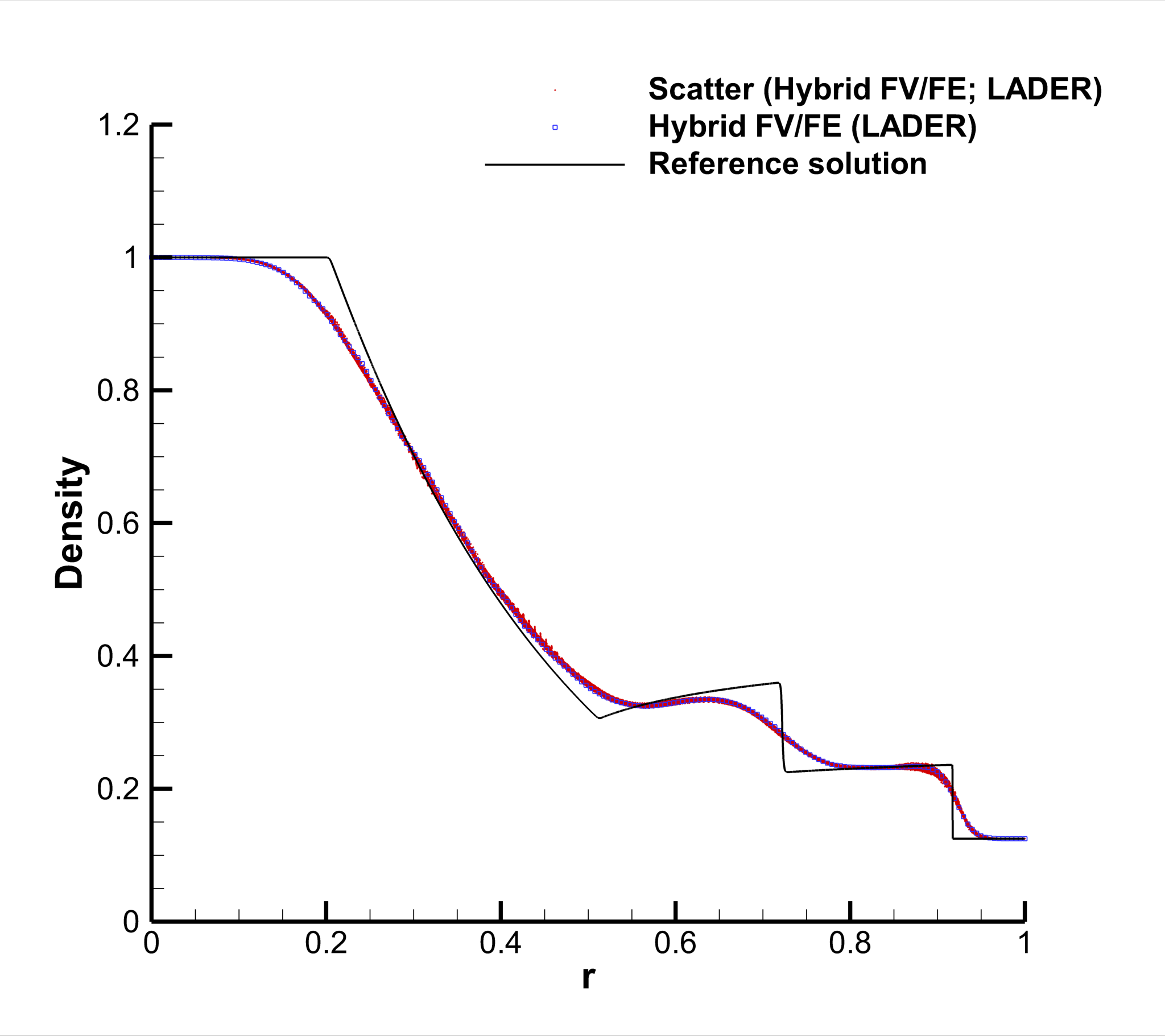}
	
	\vspace{0.05\linewidth}
	\includegraphics[trim= 10 5 5 5,clip,width=0.45\linewidth]{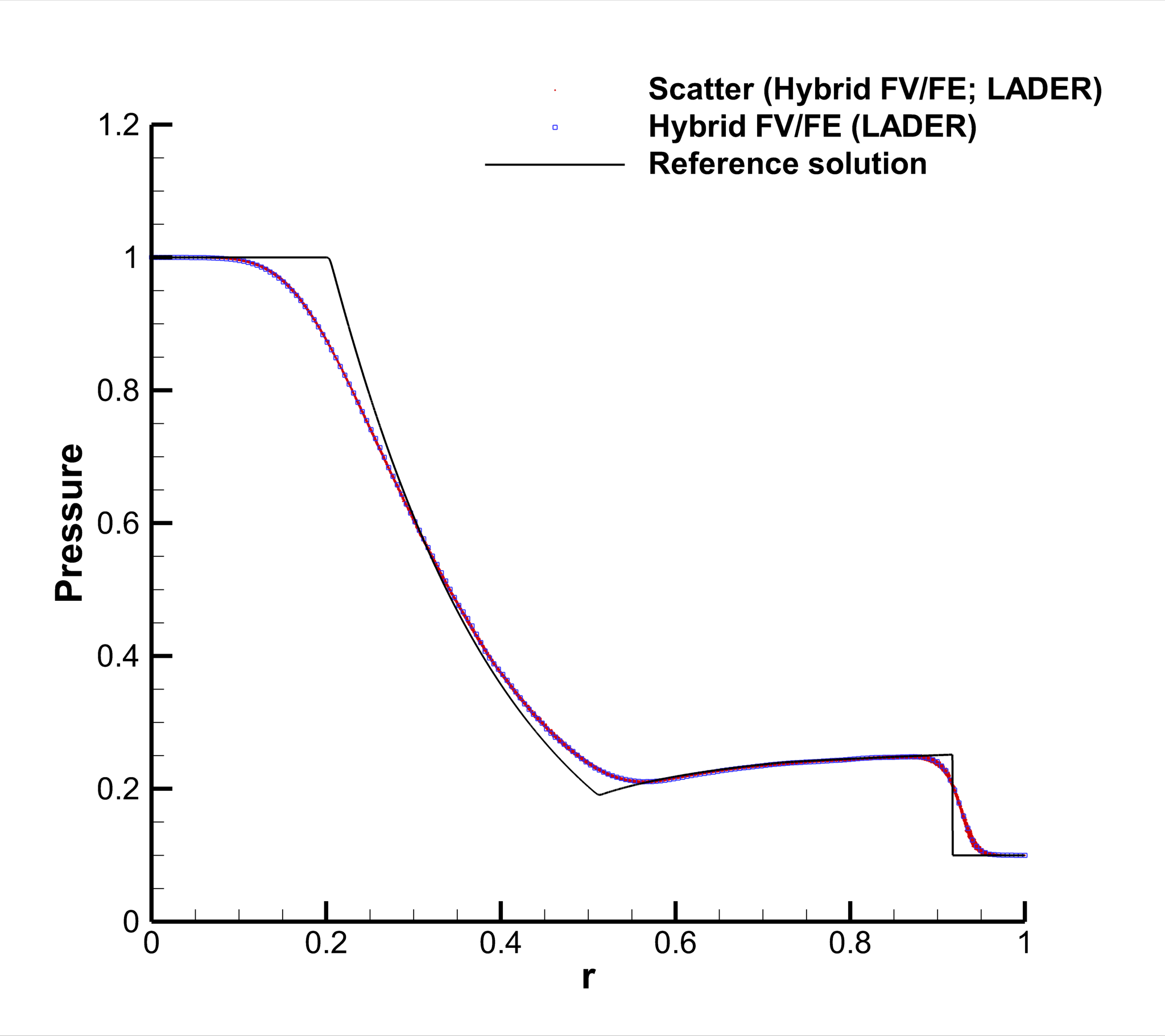}\hspace{0.05\linewidth}
	\includegraphics[trim= 10 5 5 5,clip,width=0.45\linewidth]{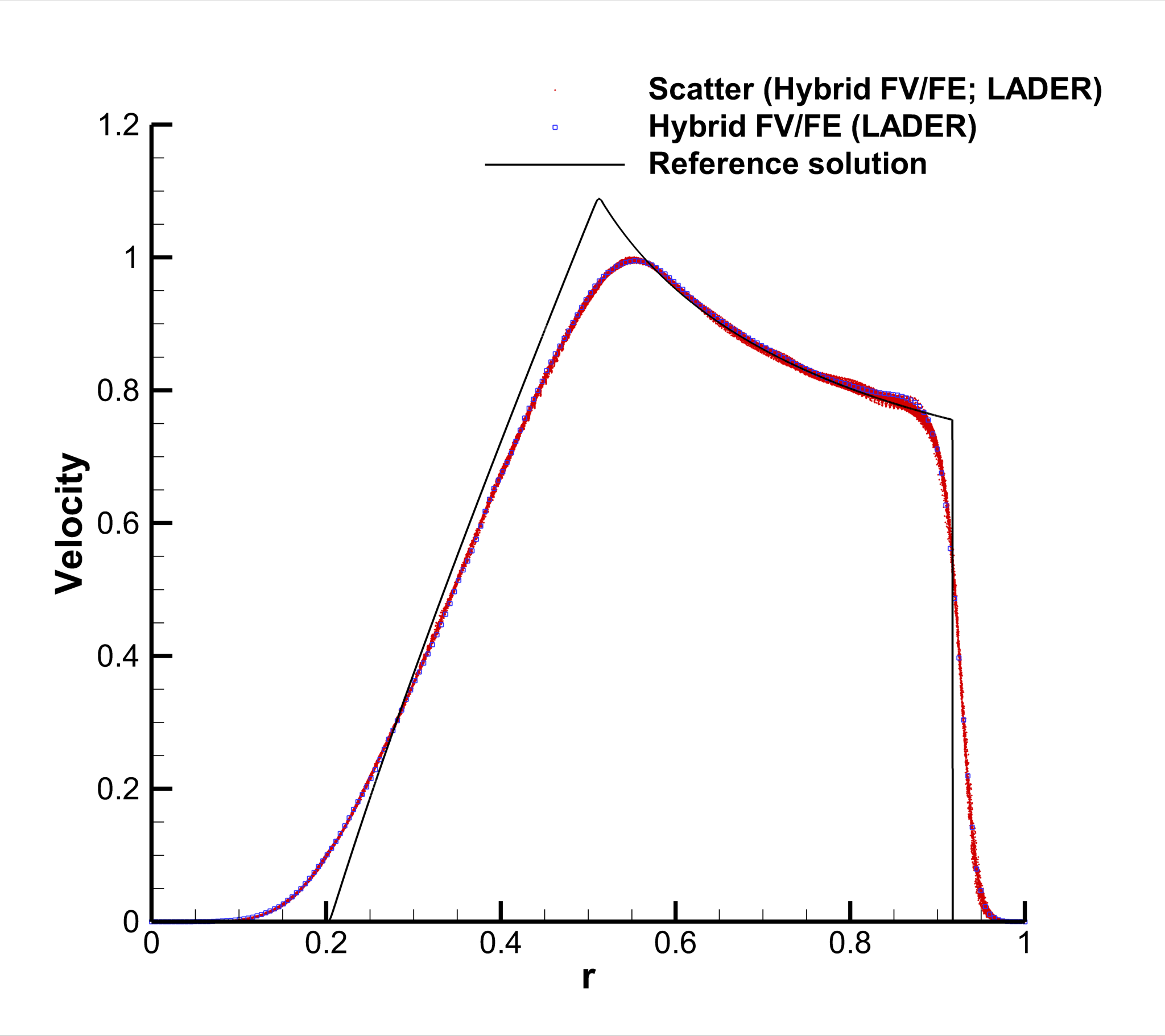}
	\caption{Circular explosion. The left top image corresponds to the 3D plot of the obtained $\rho$ at the final time whereas the 1D plots containing the reference solution (black continuous line), a 1D cut (blue squared line) and the scatter plot (red dots), correspond to the $\rho$, $\press$, and $\left|\vel\right|$ fields obtained using the LADER-ENO scheme \textcolor{black}{($c_{\alpha }=1$)}.} 
	\label{CE85_ader_t025}
\end{figure}

\begin{figure}[!htbp]
	\centering
	\includegraphics[trim= 10 5 5 5,clip,width=0.33\linewidth]{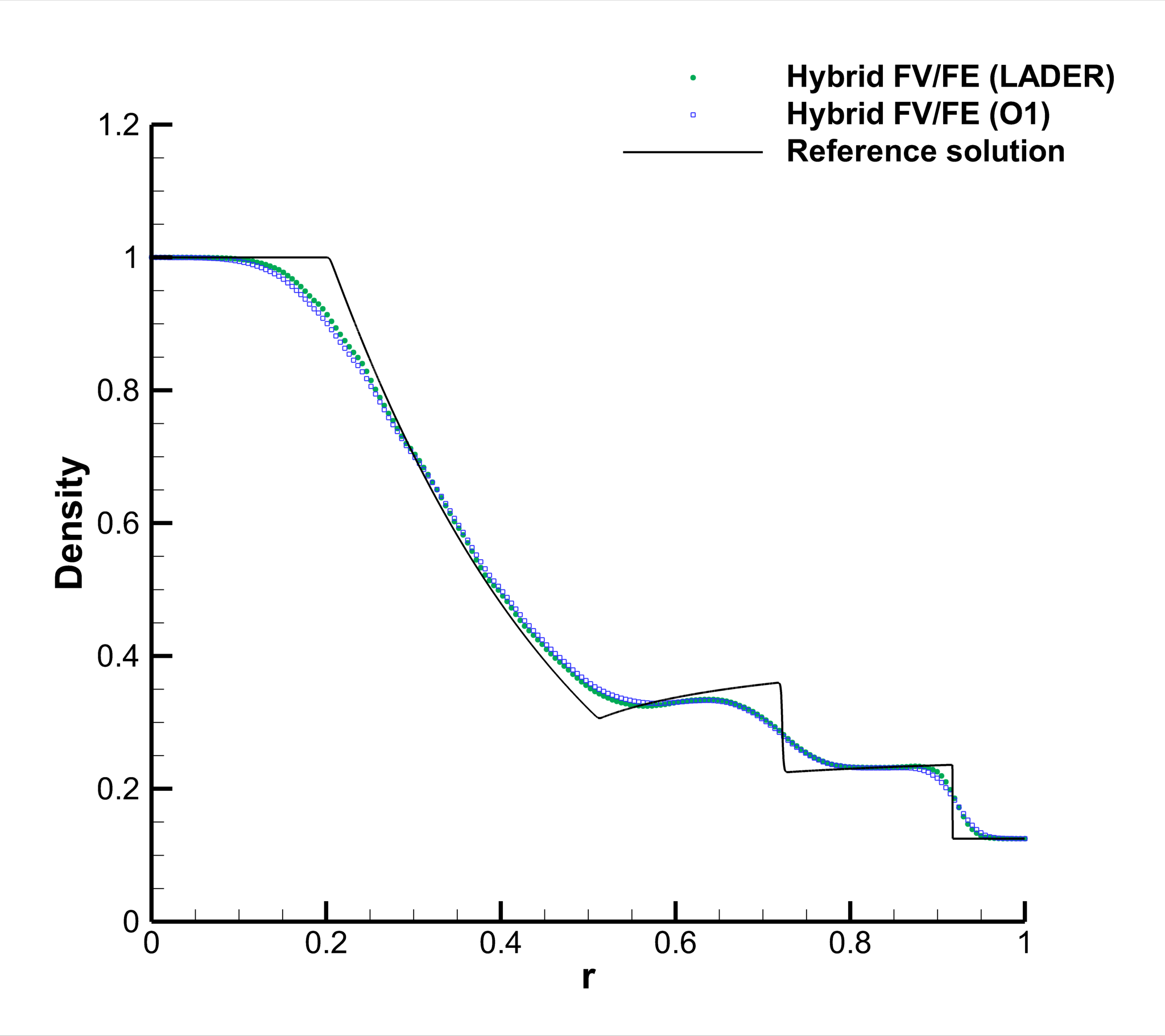}\hfill
	\includegraphics[trim= 10 5 5 5,clip,width=0.33\linewidth]{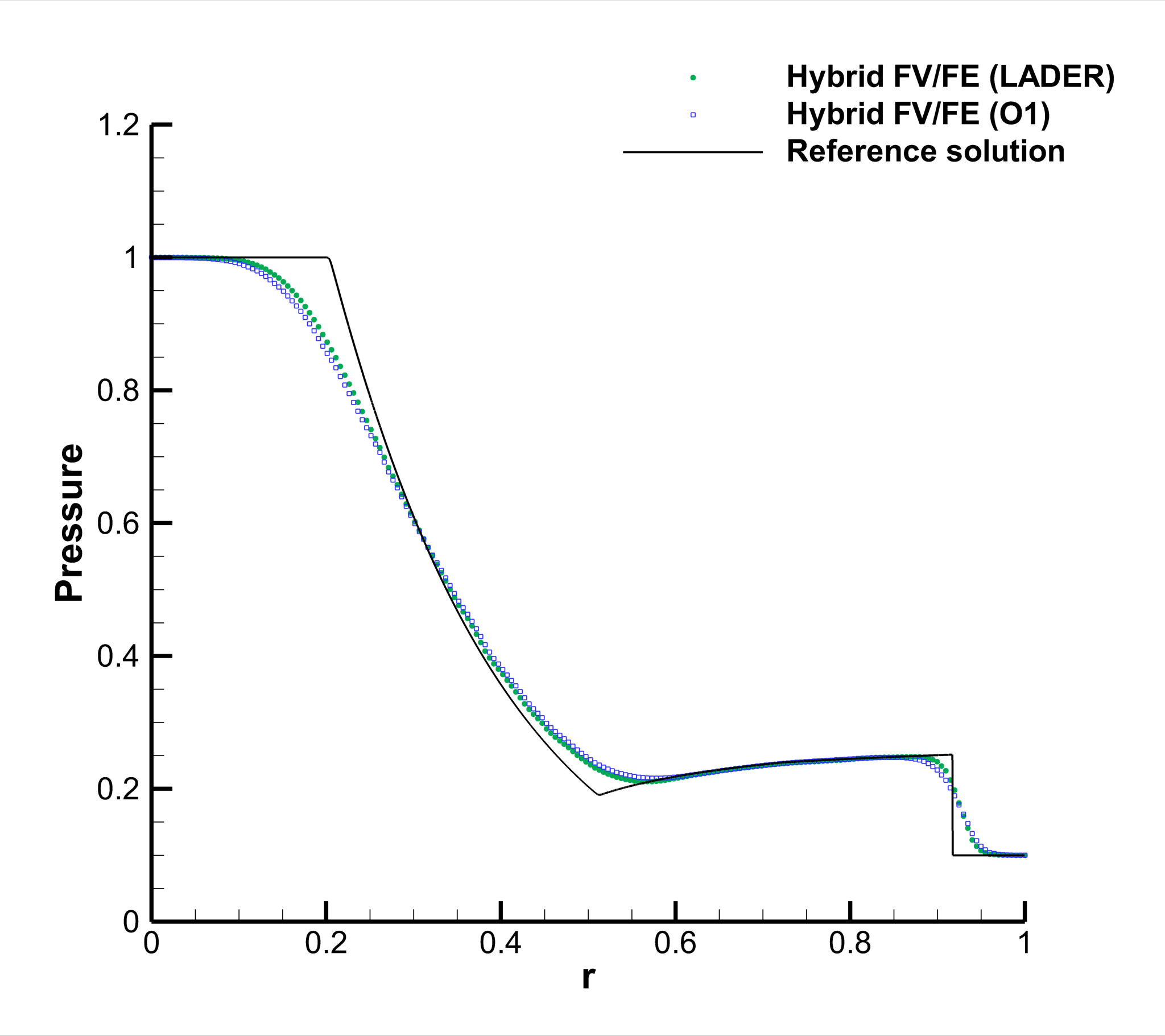}\hfill
	\includegraphics[trim= 10 5 5 5,clip,width=0.33\linewidth]{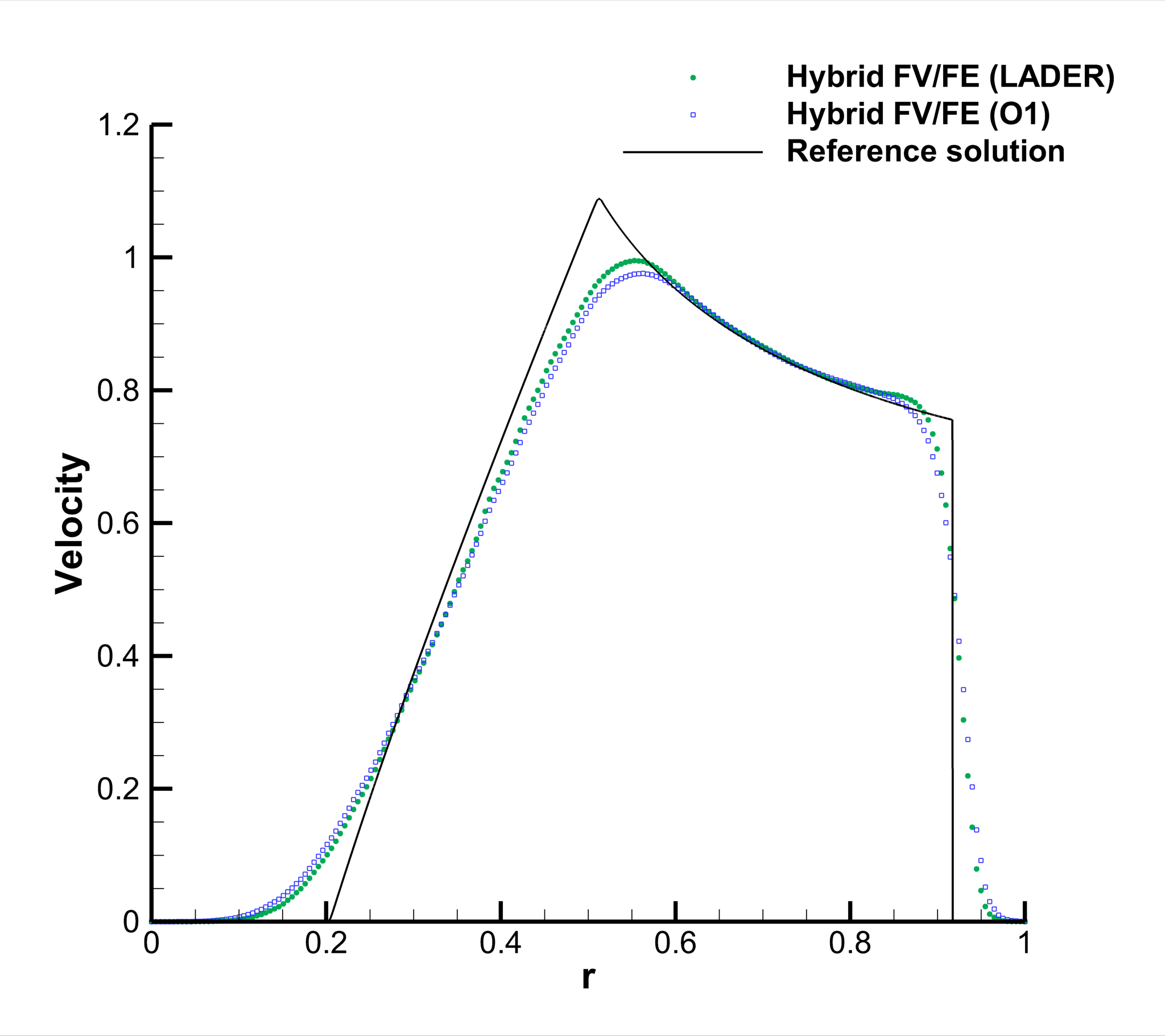}
	\caption{\textcolor{black}{Circular explosion. Comparison between the numerical solution obtained with the first order scheme (blue squares) and the LADER-ENO approximation (green circles).}} 
	\label{CE85_comparative_t025}
\end{figure}

The circular explosion problem presented here is based on an initial radial solution given by the Sod shock tube
\begin{equation}
\rho^{0}\left(\mathbf{x}\right) =  \left\lbrace \begin{array}{lr}
1 & \mathrm{ if } \; r \le 0.5,\\
0.125 & \mathrm{ if } \; r > 0.5,
\end{array}\right. \qquad
\mathbf{u}^{0} \left(\mathbf{x}\right) = 0, \qquad
\press^{0} \left(\mathbf{x}\right) = \left\lbrace \begin{array}{lr}
1 & \mathrm{ if } \; r \le 0.5,\\
0.1 & \mathrm{ if } \; r > 0.5,
\end{array}\right.
\end{equation}
see \cite{Toro,TT05,DPRZ16}. We consider the computational domain $\Omega=[-1,1]\times[-1,1]$ and periodic boundary conditions everywhere. The simulation is run until time $t_{\mathrm{end}}=0.25$ on a primal triangular mesh of $85344$ elements. 
To get a reference solution, a one-dimensional PDE in the radial direction obtained from the compressible Euler equations when using convenient geometrical source terms, \cite{Toro}, is solved using a second order TVD scheme on a very fine mesh made of $10000$ elements. The results obtained with the first order scheme and the LADER-ENO methodology, Figures \ref{CE85_o1_t025}-\ref{CE85_ader_t025},  present a good agreement with the reference solution. \textcolor{black}{Figure \ref{CE85_comparative_t025} allows for a direct comparison of the solution obtained with both schemes along a 1D cut. The second order LADER method with ENO reconstruction provides a better approximation of the solution compared to the first order scheme, as expected.}

\subsection{3D spherical explosion}
In this section, we study the behaviour of the method for a 3D spherical explosion benchmark based on the Sod problem. 
The computational domain is defined to be the sphere of unit radius centered at the origin. Initial conditions read
\begin{equation}
\rho^{0}\left(\mathbf{x}\right) =  \left\lbrace \begin{array}{lr}
1 & \mathrm{ if } \; r \le \frac{1}{2},\\[6pt]
0.125 & \mathrm{ if } \; r > \frac{1}{2},
\end{array}\right. \qquad
\press^{0} \left(\mathbf{x}\right) = \left\lbrace \begin{array}{lr}
1 & \mathrm{ if } \; r \le \frac{1}{2},\\[6pt]
0.1 & \mathrm{ if } \; r > \frac{1}{2},
\end{array}\right. \qquad
\mathbf{u}^{0} \left(\mathbf{x}\right) = 0, \label{eq:IC_CE3D}
\end{equation}
with $r=\sqrt{x^{2}+y^{2}+z^{2}}$. Dirichlet boundary conditions are imposed and the domain is covered by $2280182$ tetrahedra. 

\begin{figure}[!htbp]
	\centering 	
	\includegraphics[trim= 5 0 5 0,clip,width=0.45\linewidth]{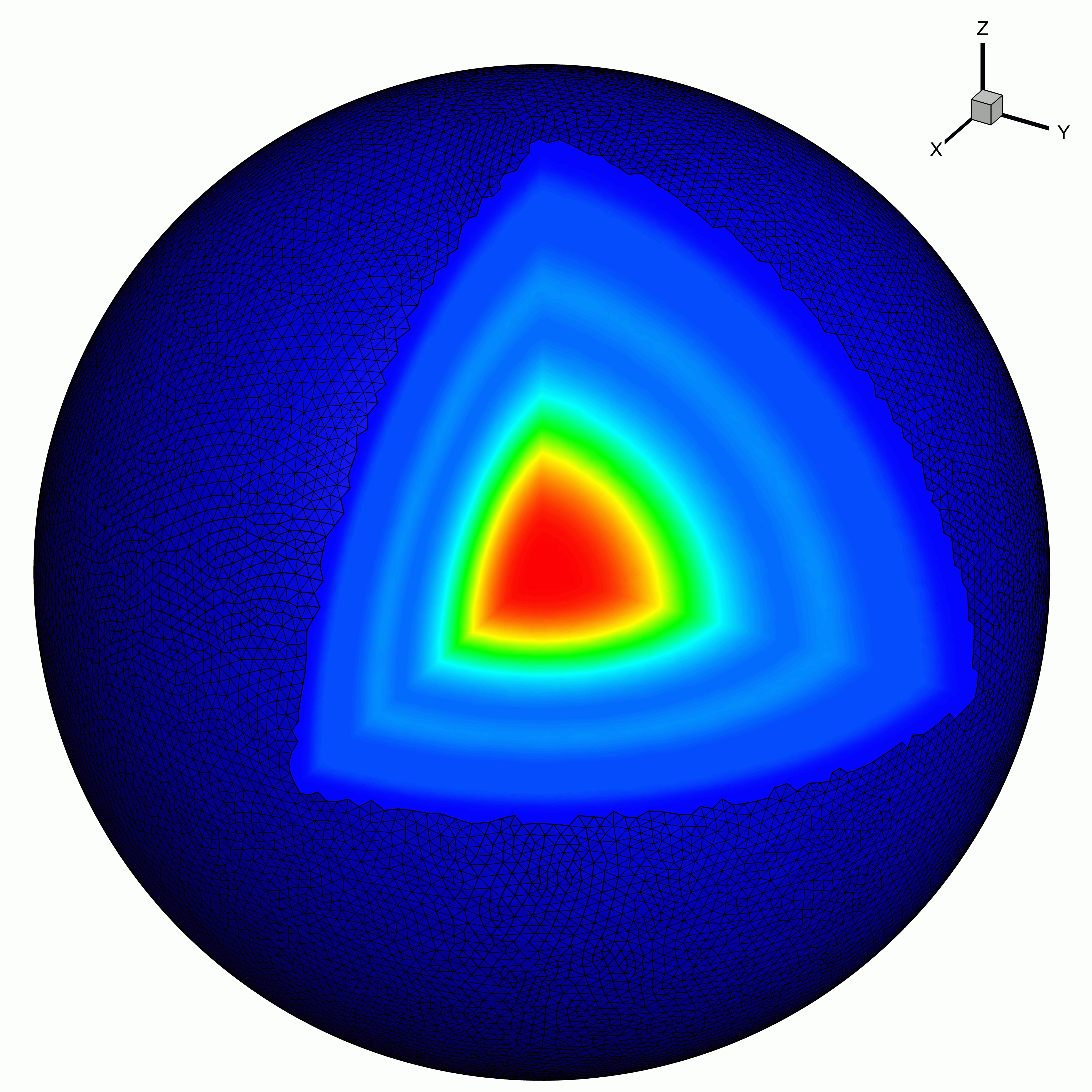}\hspace{0.05\linewidth}
	\includegraphics[trim= 5 0 5 0,clip,width=0.45\linewidth]{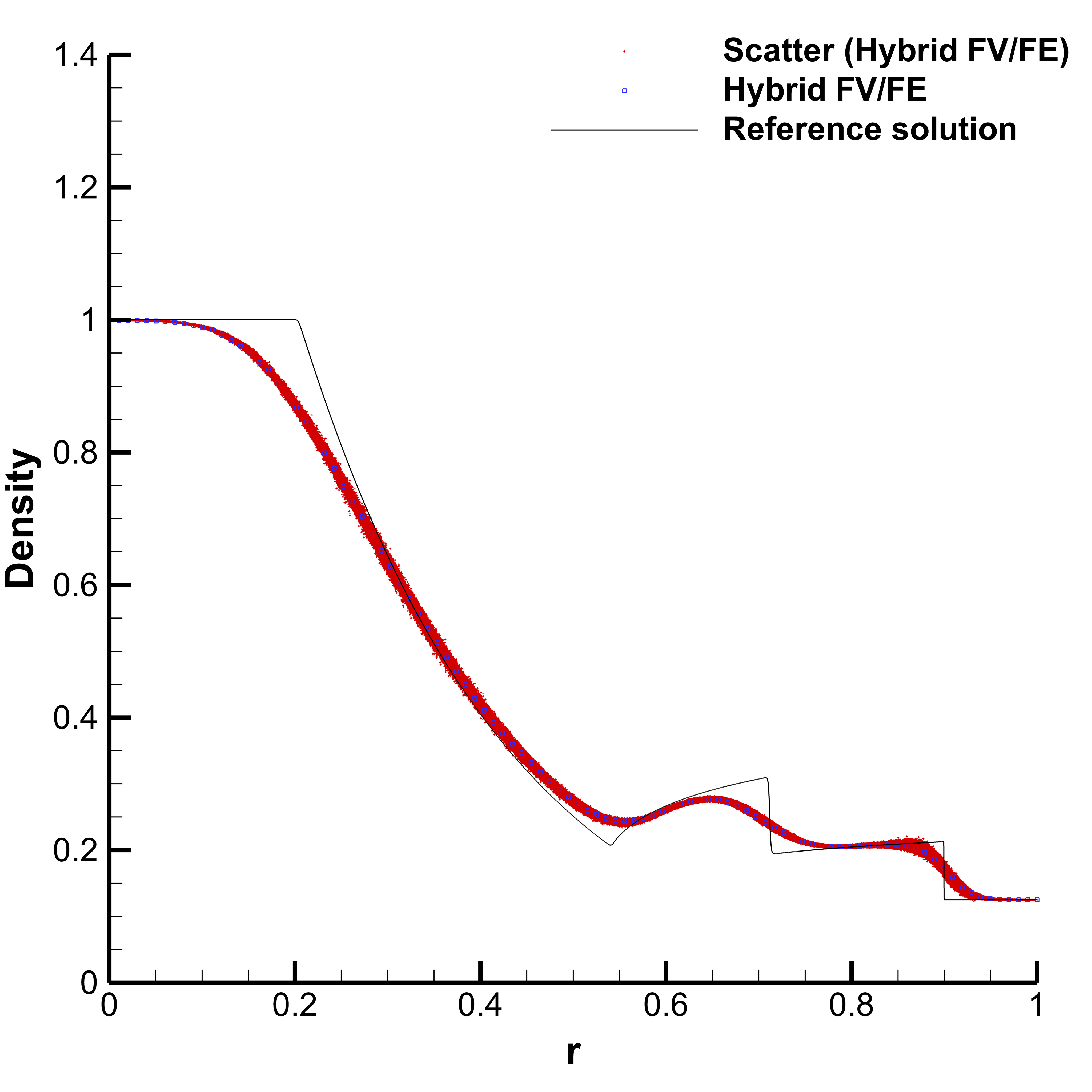}
	
	\vspace{0.05\linewidth}
	\includegraphics[trim= 5 0 5 0,clip,width=0.45\linewidth]{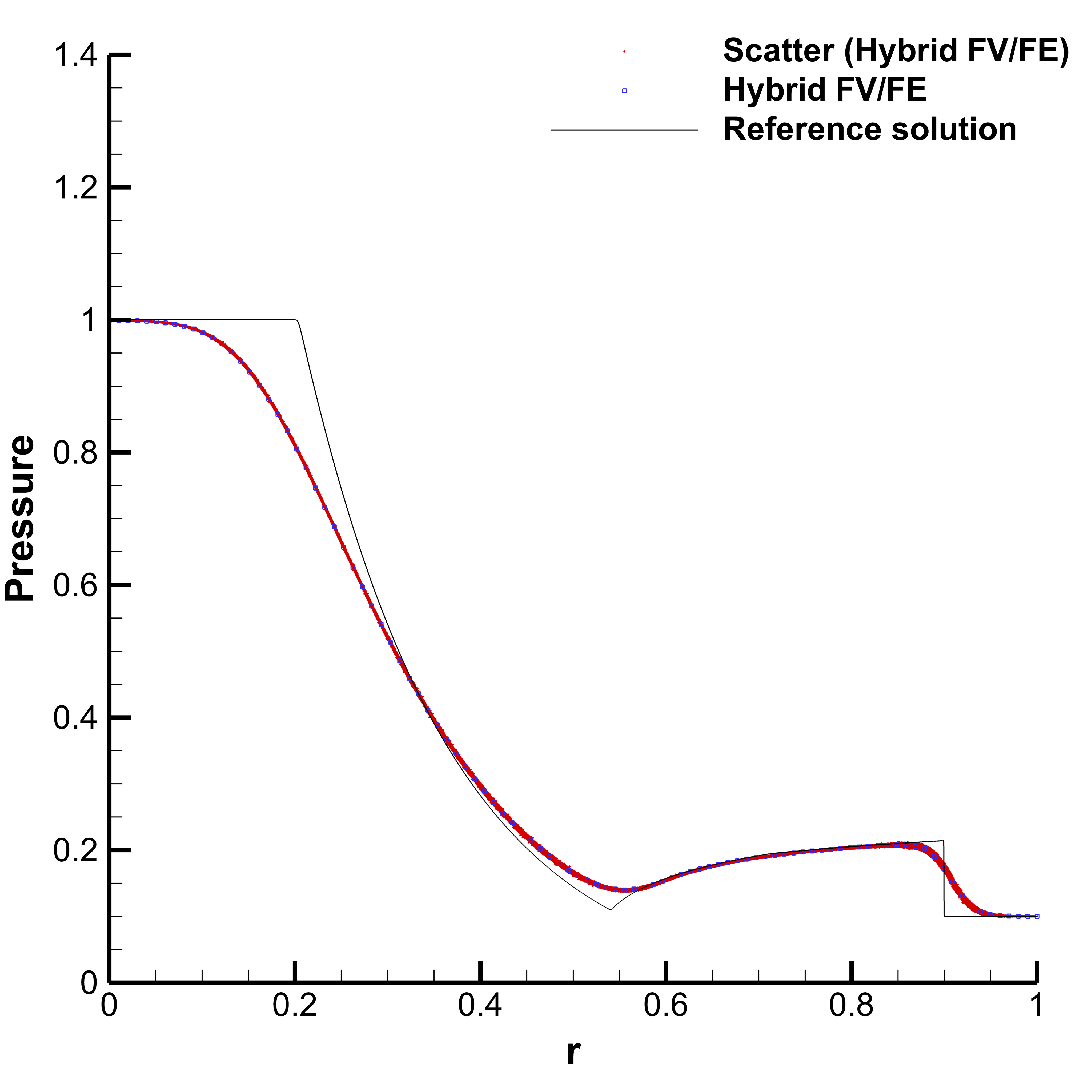}\hspace{0.05\linewidth}
	\includegraphics[trim= 5 0 5 0,clip,width=0.45\linewidth]{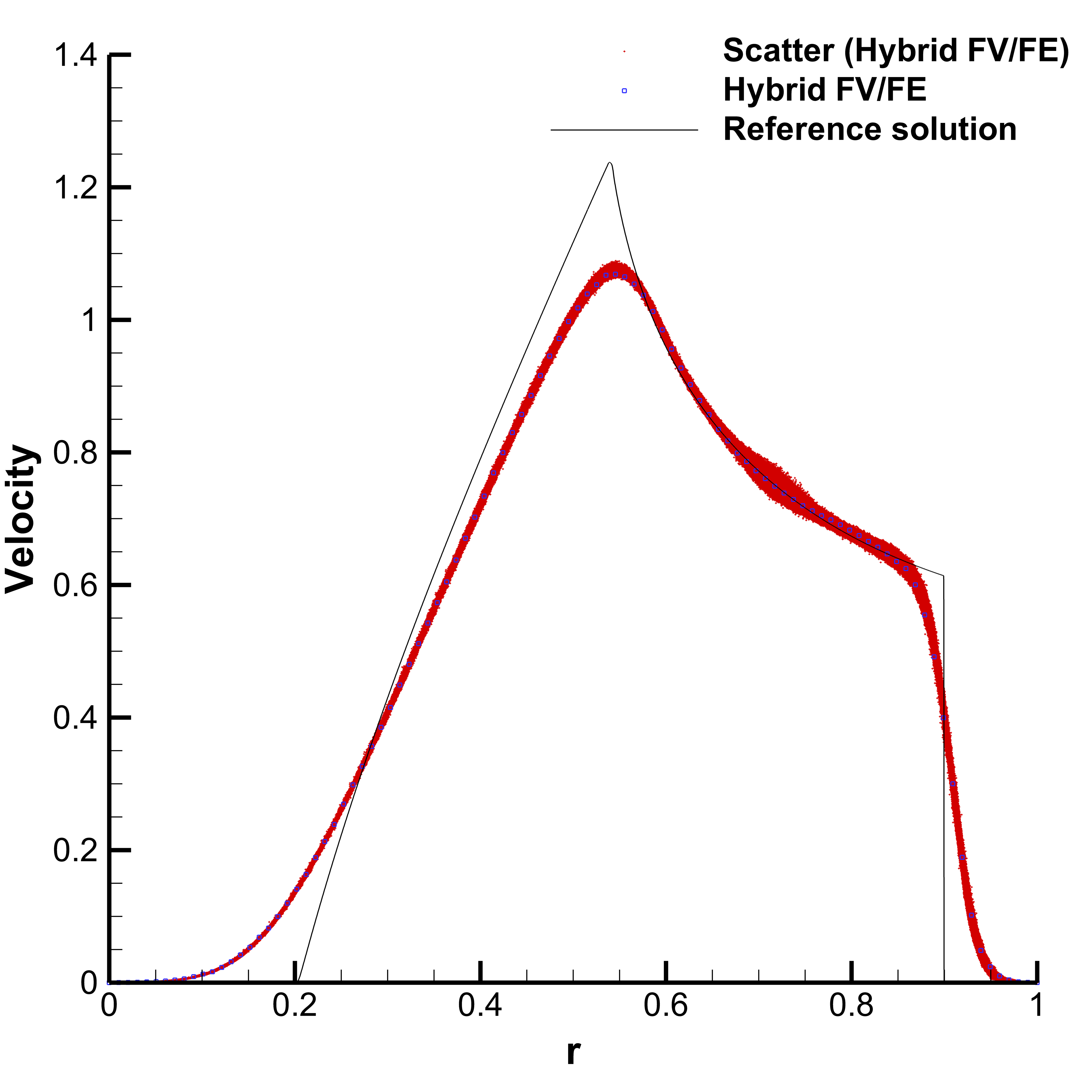}
	\caption{LADER-ENO solution for the 3D spherical explosion test at final time. Surface mesh on the boundary and density contours on the interior surfaces obtained after taking away the first quadrant and
	1D plots for density, velocity magnitude and pressure fields: 1D cut on  $x\in\left[0,1\right], \, y=z=0$ (blue squares), scatter plot (red dots), reference solution (black line).}
	\label{fig:CE3D}
\end{figure}

The solution obtained using the LADER-ENO scheme with $\mathrm{CFL}=1$, \textcolor{black}{$c_{\alpha}=3$}, up to $t_{\mathrm{end}}=0.25$ is depicted in Figure \ref{fig:CE3D}. As reference solution we solve again the 1D code for Euler equations introduced in Section \ref{sec:2DCE} updated with appropriate source terms to account for three-dimensional effects.
The agreement observed for the 1D cuts of density, velocity magnitude and pressure prove the capability of the method to handle three-dimensional problems.

\subsection{First problem of Stokes}
\begin{figure}[!htbp]
	\centering
	\includegraphics[width=0.325\linewidth]{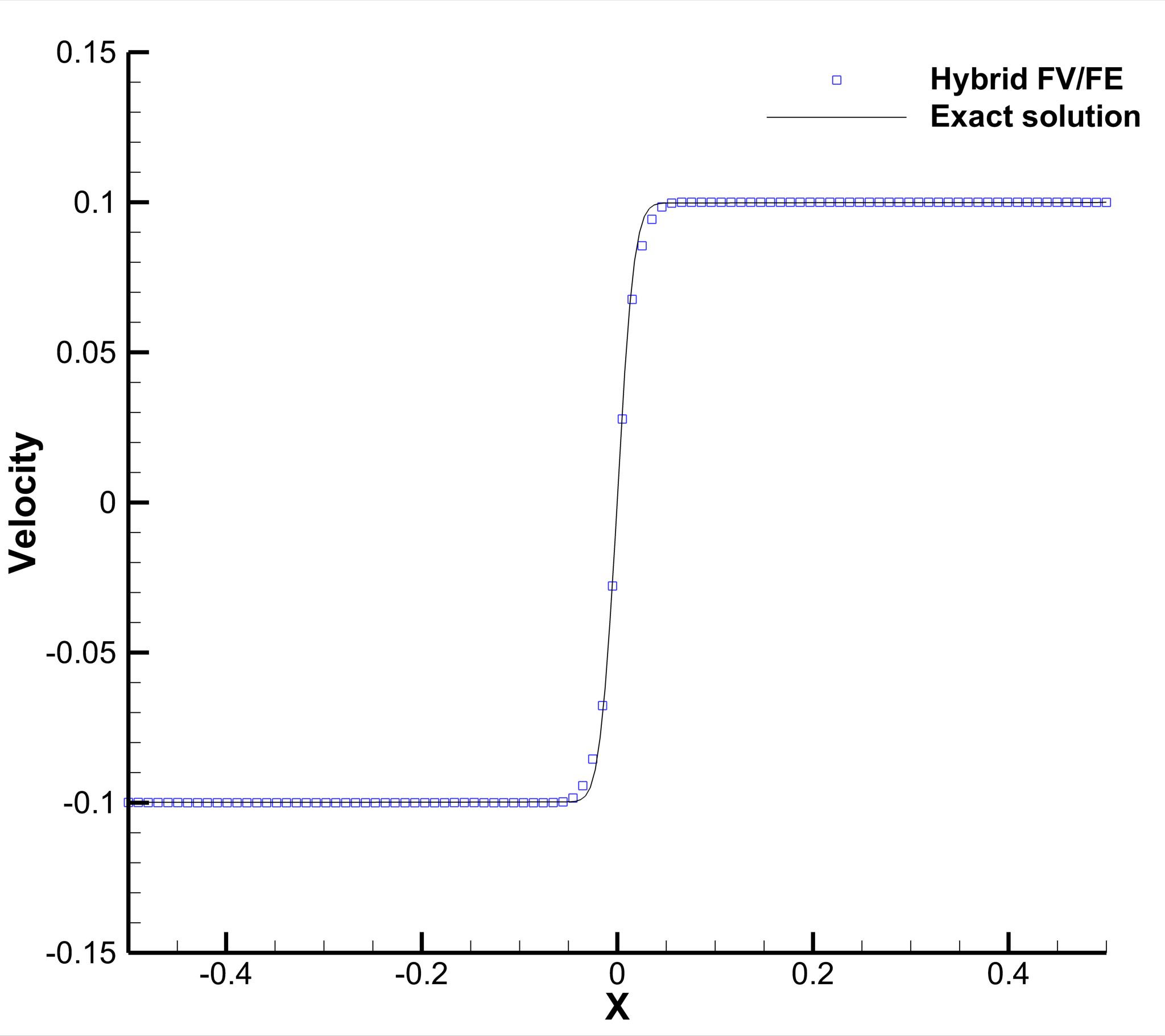}
	\includegraphics[width=0.325\linewidth]{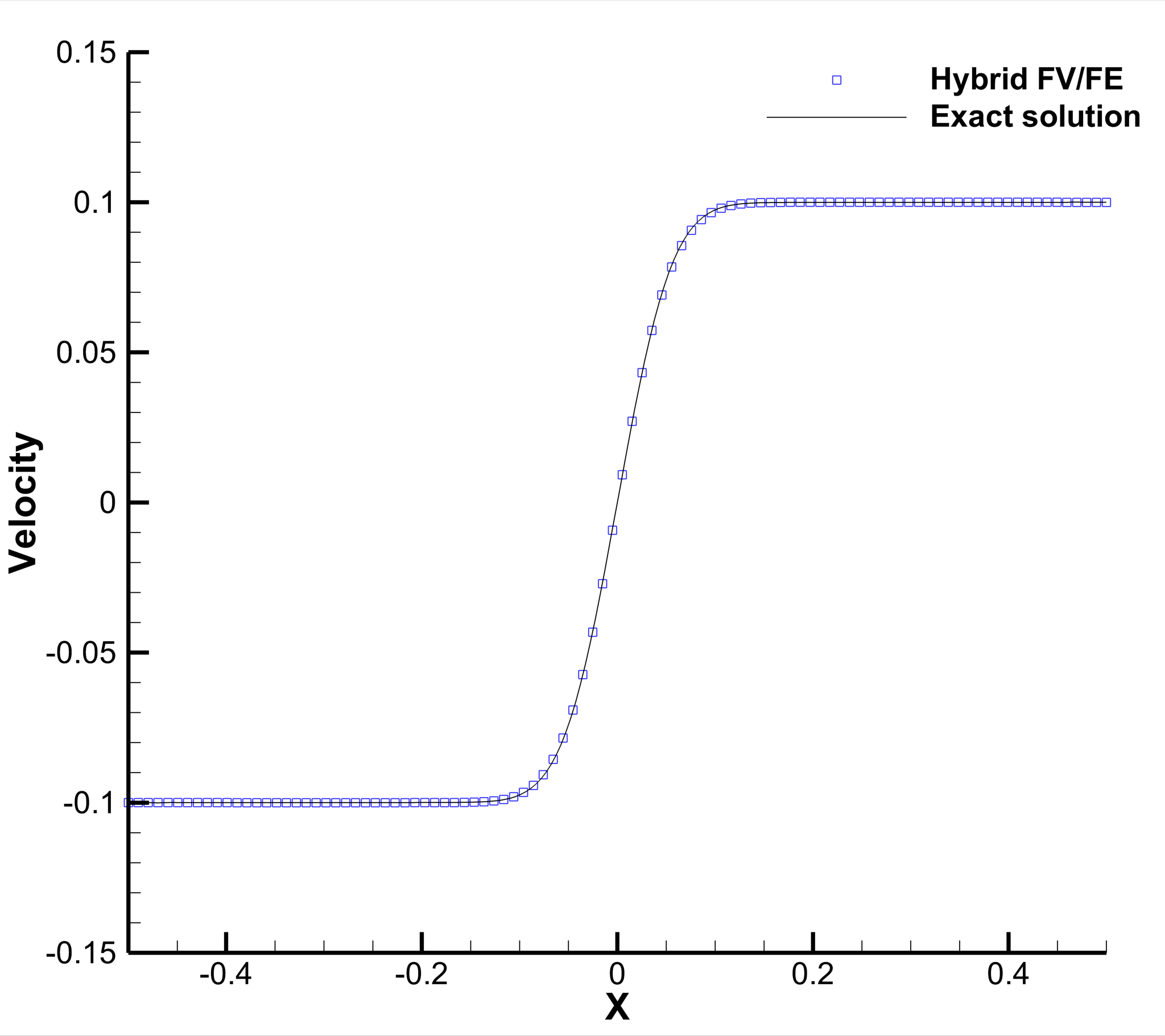}
	\includegraphics[width=0.325\linewidth]{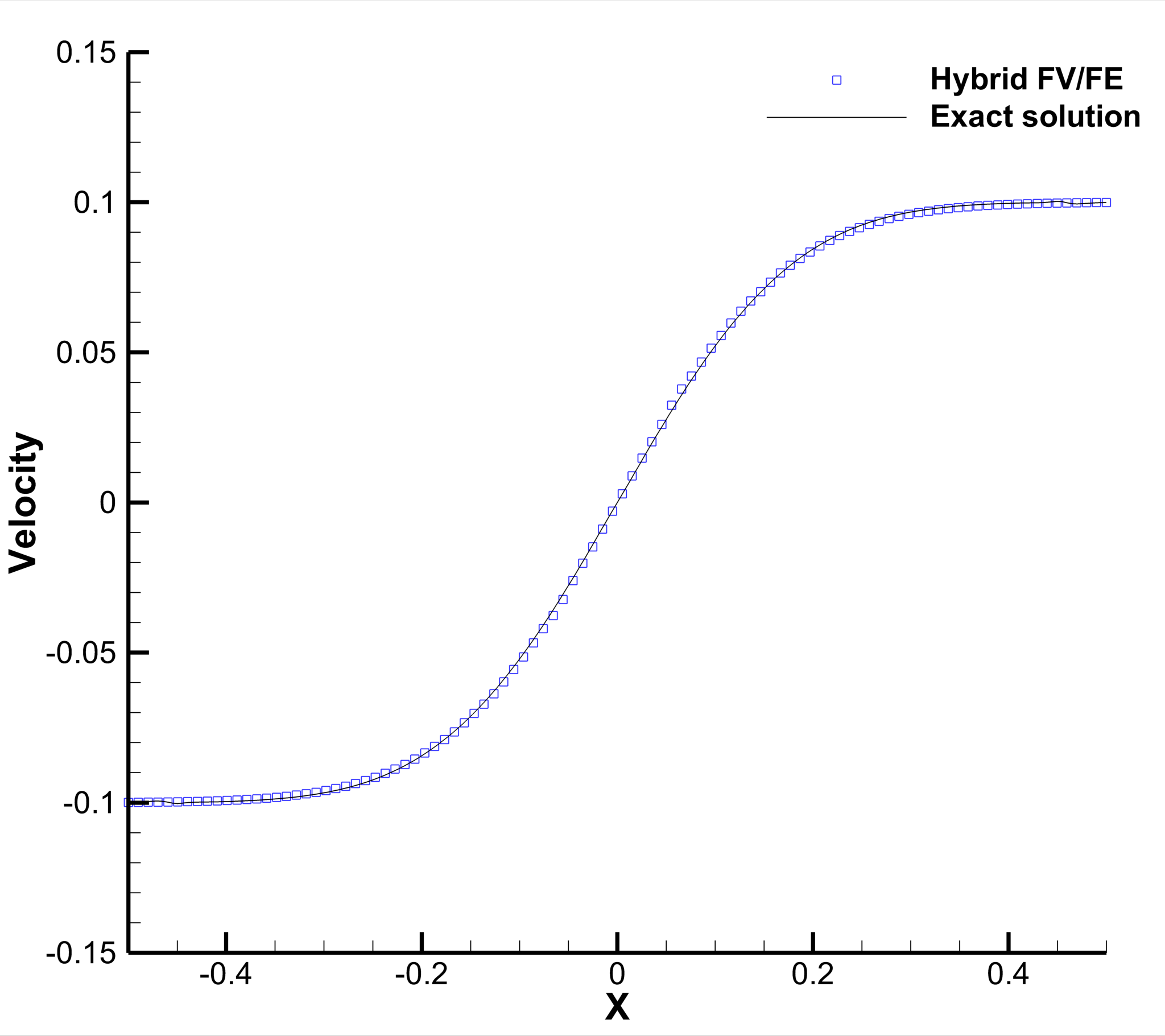}
	\caption{Comparison of the vertical velocity  $u_{2}$ along the 1D cut $y=0$ obtained for the first Stokes problem using LADER scheme against the exact solution at $t_{\mathrm{end}}=1$. $\mu=10^{-4}$ LADER-EN (left), $\mu=10^{-3}$ LADER without limiters (center), $\mu=10^{-2}$ LADER without limiters (right).}
	\label{FSP_vvelocity}
\end{figure}

To further analyse the behaviour of the developed method in the incompressible limit, we now consider the first problem of Stokes, \cite{SG16}.
The initial condition, defined in $\Omega=[-0.5,0.5]\times[-0.5,0.5]$, reads
\begin{equation}
\rho^{0}\left(\mathbf{x}\right) = 1,\qquad
\press^{0} \left(\mathbf{x}\right) = \frac{1}{\gamma}, \qquad
{u}_{1}^{0} \left(\mathbf{x}\right) = 0, \qquad
{u}_{2}^{0} \left(\mathbf{x}\right) = \left\lbrace \begin{array}{lc}
-0.1 & \mathrm{ if } \; y \le 0,\\
0.1 & \mathrm{ if } \; y > 0
\end{array}\right.
\end{equation}
In the incompressible limit, this test case has an exact analytical solution for $u_{2}$ given by
\begin{equation}
{u}_{2} \left(\mathbf{x},t\right) = \frac{1}{10} \mathrm{erf}\left( \frac{x}{2\sqrt{\mu t}}\right).
\end{equation}
To complete the physical set up, we define $\gamma= c_{\press} = 1.4$, $\lambda=0$, leading to \mbox{$M\approx 10^{-1}$}.
\textcolor{black}{Regarding boundary conditions, we impose the exact values for density and velocity in the $x$-direction, while on the top and bottom boundaries, we set periodic boundary conditions in $y$-direction.} Meanwhile, the exact values for density and velocity are employed in the remaining boundaries.
Finally, three different simulations are run attending to the value for the viscosity coefficient: $\mu= 10^{-2}$, $\mu=10^{-3}$, and $\mu=10^{-4}$. The simulations are run on a triangular primal mesh made of $1000$ elements up to time $t_{\mathrm{end}}=1$. The vertical velocity along  $y=0$ is plotted in Figure~\ref{FSP_vvelocity} against the exact solution. We observe a good agreement between both curves for all three viscosities. Let us note that $\mu=10^{-4}$ is the only simulation run using the ENO reconstruction so that we completely avoid the small bump arising after the discontinuity if any limiting strategy is employed. In the other cases, such reconstruction can be neglected due to the high physical viscosity considered.

\subsection{Viscous shock} 

Here we analyse a steady viscous shock with $M_s=2$ the shock Mach number. 
Considering the particular case Pr$=0.75$, with Pr the Prandtl number, it is possible to find an exact solution of the compressible Navier-Stokes equations, derived by Becker in 1923, see \cite{Becker1923,BonnetLuneau,GPRmodel} for all the necessary details to setup this test case.  

The computational domain $\Omega=[-0.5,0.5] \times [0,0.1]$ is discretized with 12500 triangular elements 
of characteristic mesh spacing $h=1/250$. The shock wave is centered at $x=0$. 

\begin{figure}[!htbp]
	\begin{center}
		\includegraphics[width=0.32\textwidth]{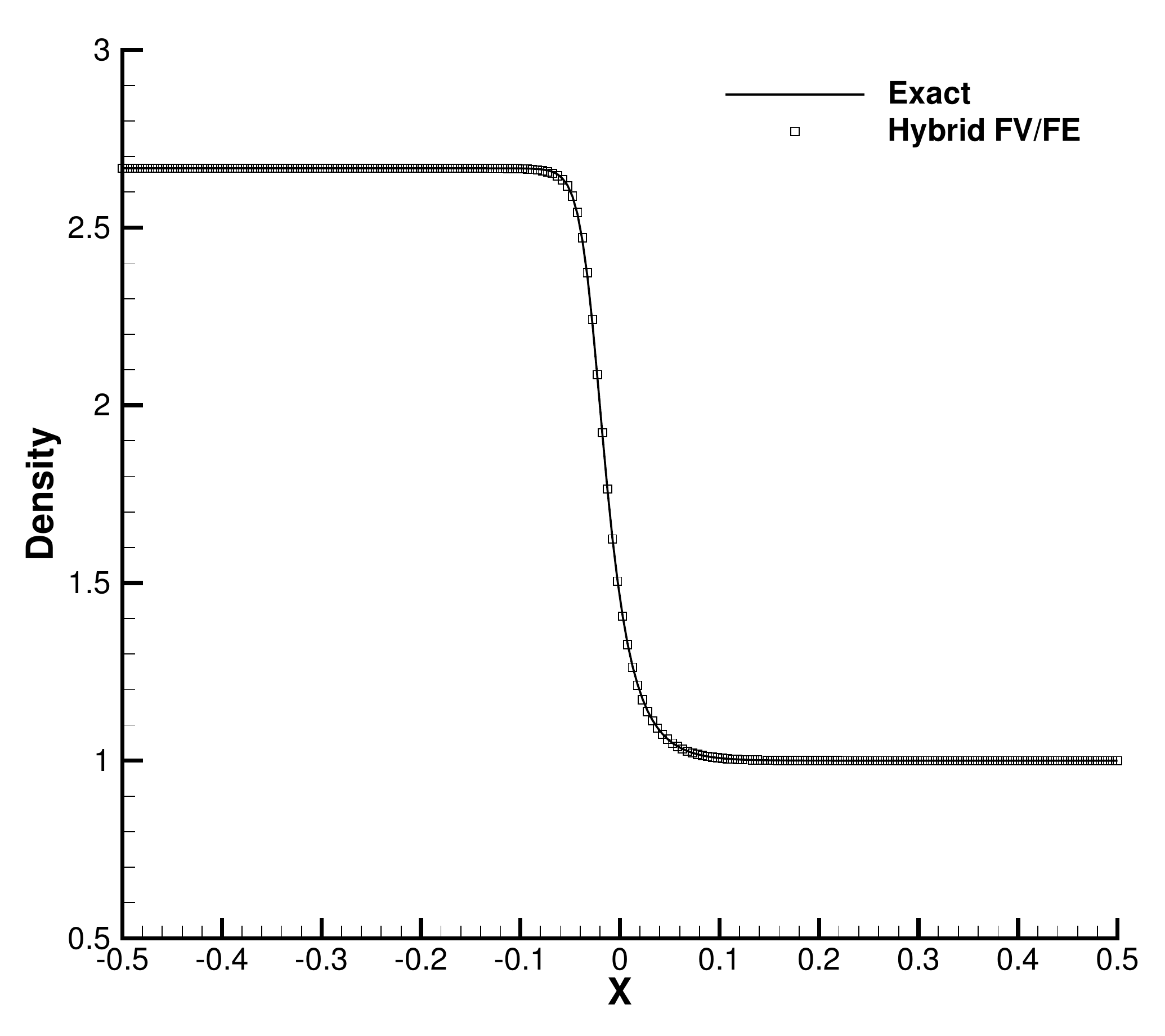}   
		\includegraphics[width=0.32\textwidth]{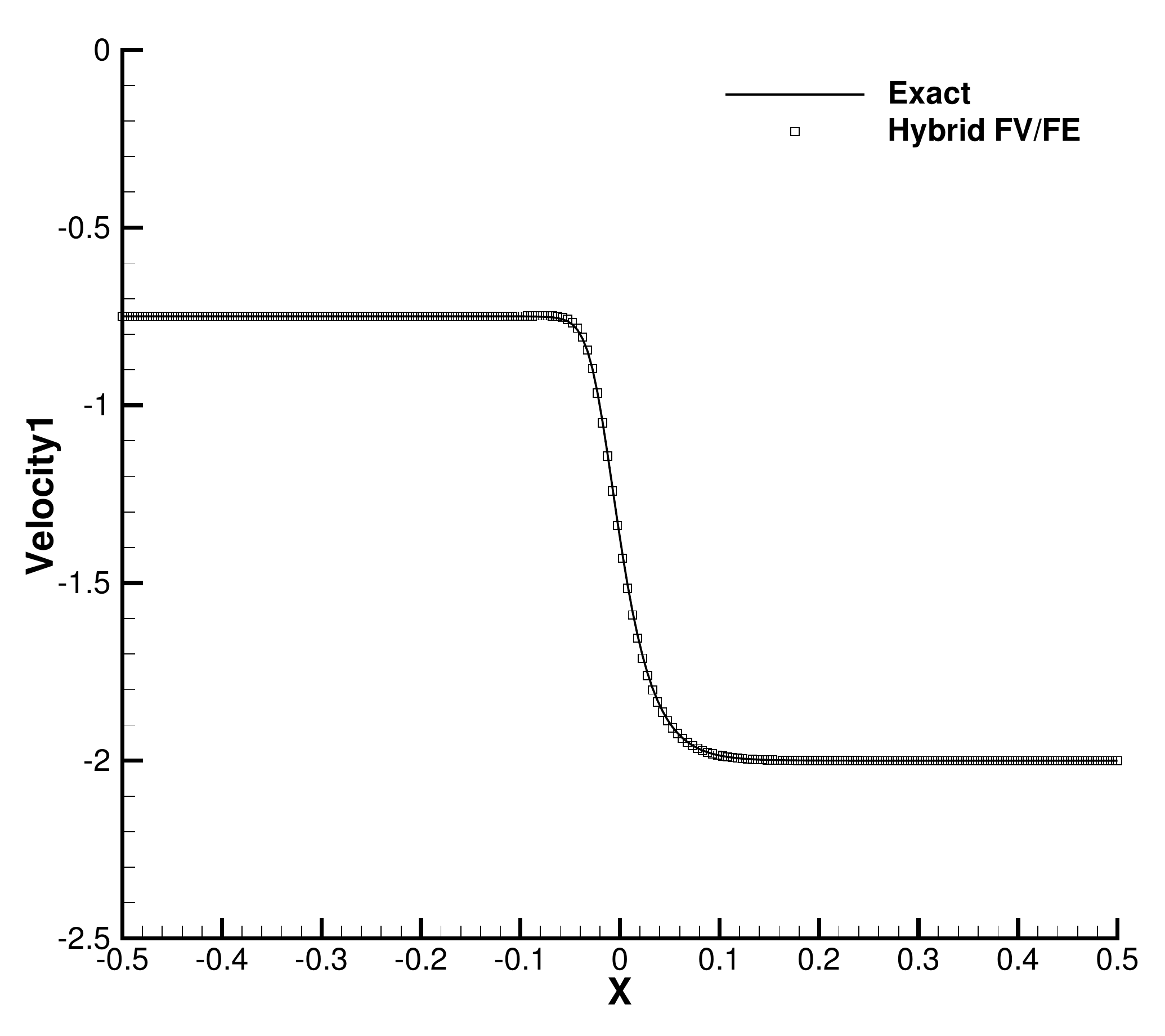}   
		\includegraphics[width=0.32\textwidth]{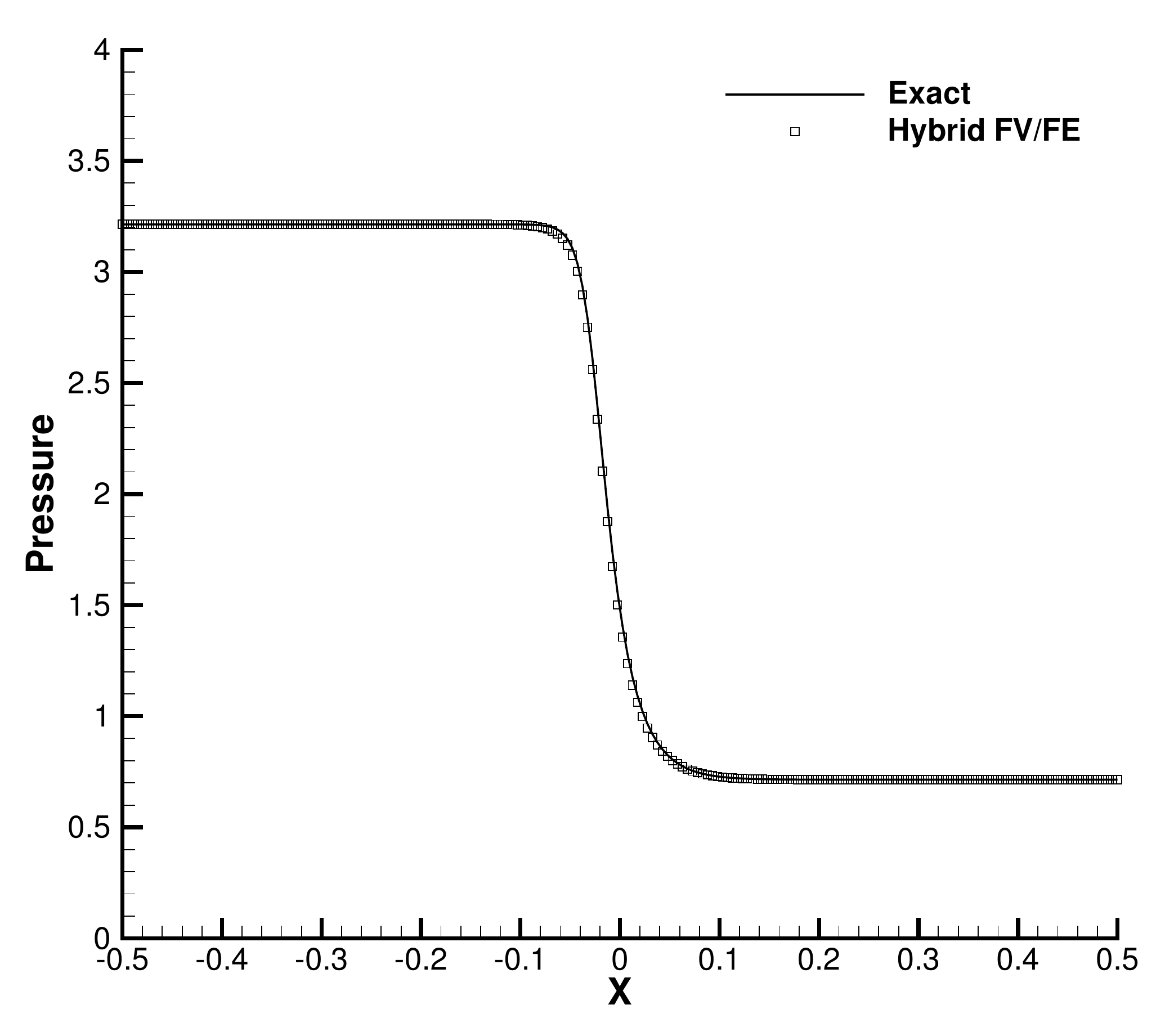}   
		\caption{Numerical solution for the steady viscous shock obtained with the hybrid FV/FE scheme and exact solution of the compressible Navier-Stokes equations, Pr$=0.75$, $M=2$ and $Re=100$ at time $t=0.025$.} 
		\label{fig.vshock}
	\end{center}
\end{figure}

The values of the fluid in front of the shock wave are given by $\rho_0 =1$, $u_0=-2$, $v_0=w_0=0$, and $p_0=1/\gamma$ so that the corresponding sound speed is $c_0 = 1$ and the fluid is moving into the shock from the right to the left at shock Mach number $M_s=2$. 
The Reynolds number based on a unitary reference 
length ($L=1$) and on the flow speed $u_0$ is given by $Re_s=\frac{\rho_0 \, c_0 \, M_s \, L }{\mu}$.  
The fluid parameters are chosen as $\gamma = 1.4$, $c_v = 2.5$, $\mu=2 \cdot 10^{-2}$  and $\lambda = 9 \frac{1}{3} \cdot 10^{-2}$, hence the corresponding shock Reynolds number is $Re_s=100$. 
The simulation with the new hybrid FV/FE scheme proposed in this paper is run until time $t_{\mathrm{end}}=0.025$, \textcolor{black}{setting $c_\alpha=3$.}   
The comparison between the numerical solution and the exact solution is shown in Figure \ref{fig.vshock} for the density $\rho$, the velocity $u$, and the pressure $p$. For all quantities, one can note a very good agreement. 

\subsection{Lid-driven cavity flow} 
\begin{figure}[!htbp]
	\centering
	\includegraphics[trim= 10 5 5 5,clip,width=0.44\linewidth]{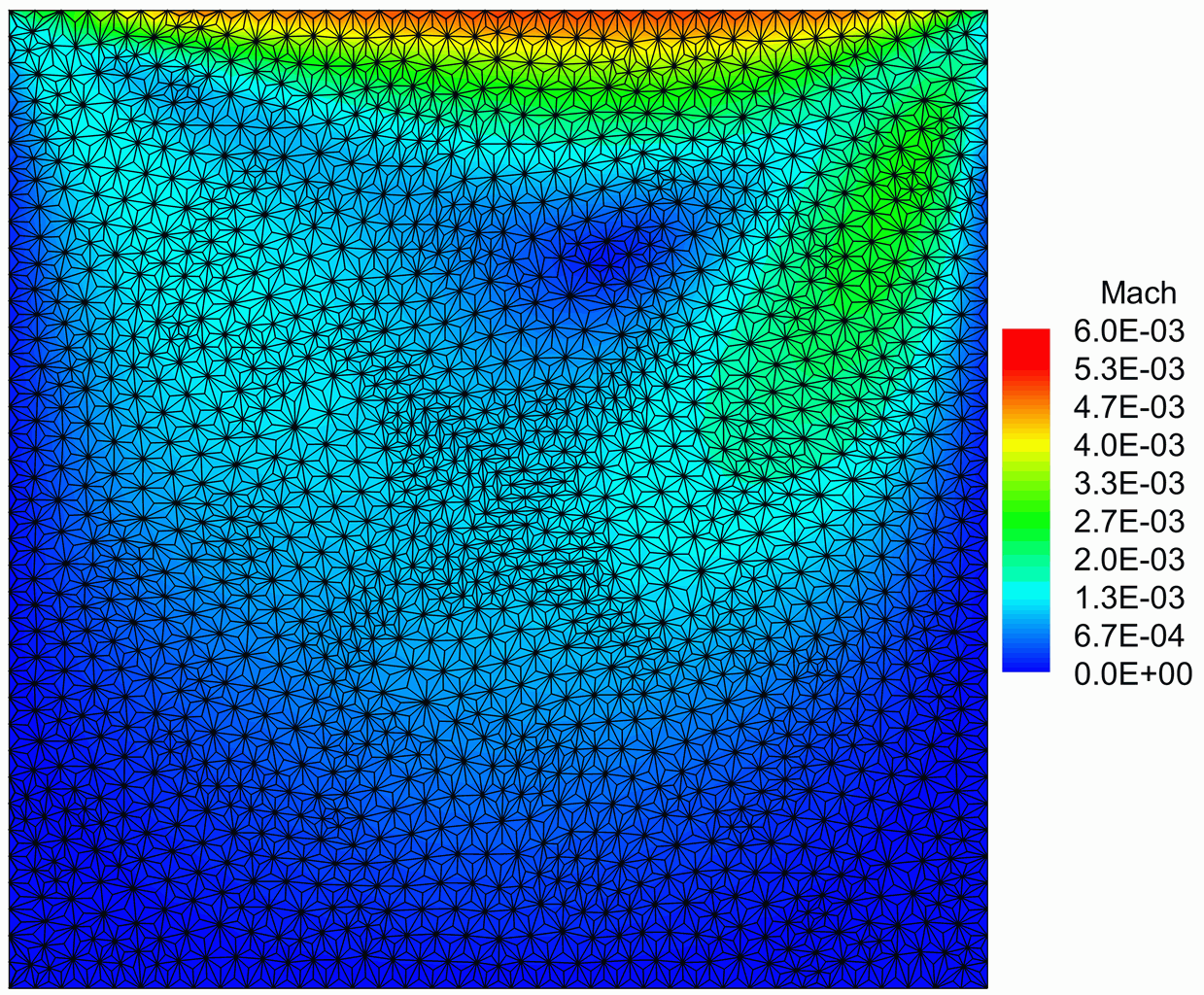}\hspace{0.1\linewidth}
	\includegraphics[trim= 10 5 5 5,clip,width=0.42\linewidth]{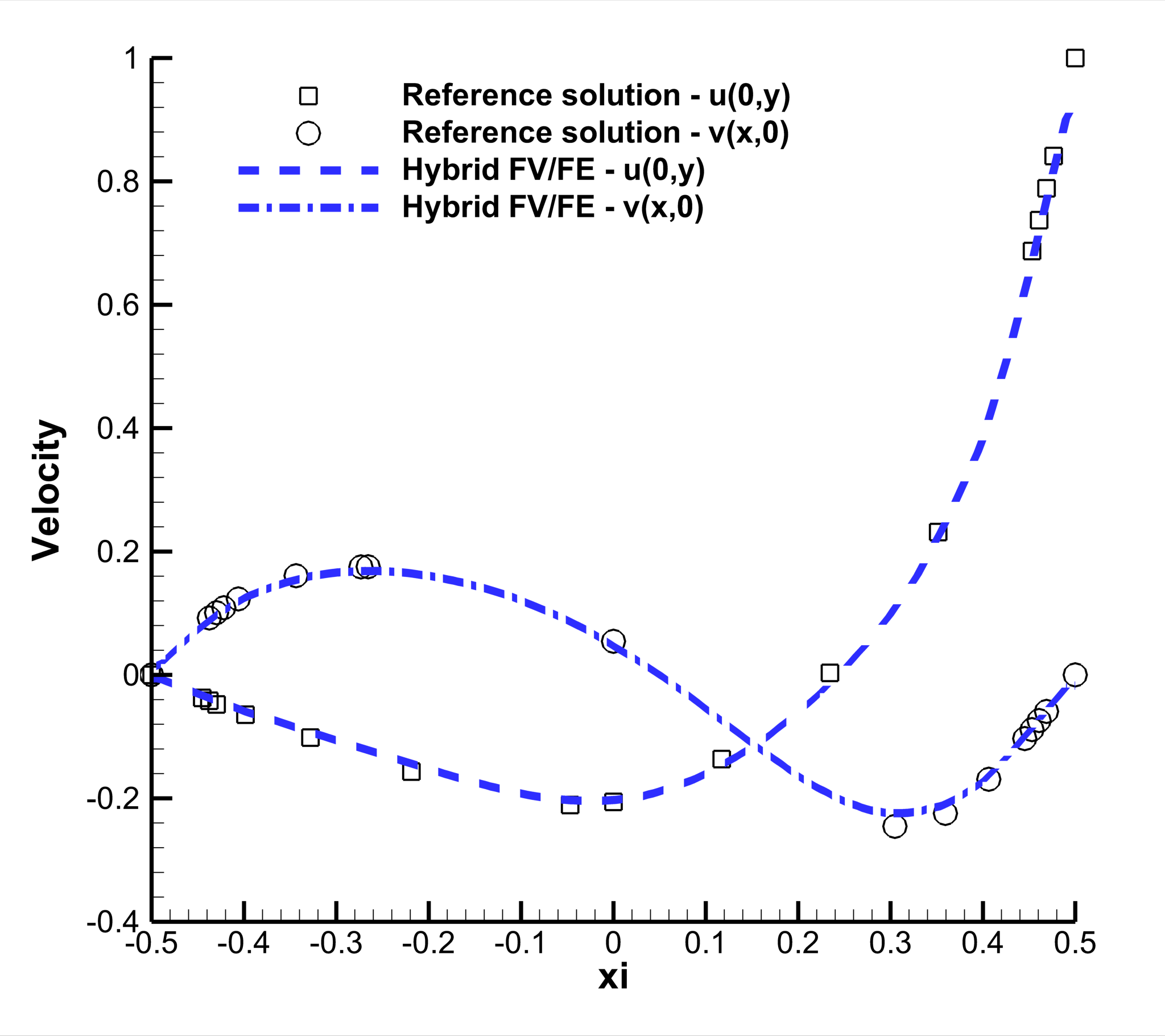}
	\caption{Solution of the lid-driven cavity test with $Re=10^{2}$. In the left subfigure the dual half grid considered is depicted over the Mach number contour plot. The right plot shows the obtained 1D velocity cuts along $x=0.5$ and $y=0.5$ (blue dashed line) and the reference solution (black circles) given by \cite{GGS82}.}
	\label{LDC_figure}
\end{figure}

A classical test for incompressible flows is the lid-driven cavity benchmark, which accounts for a well-known reference solution, \cite{GGS82}. Therefore this test may be the optimal candidate to assess the behaviour of the method in the incompressible limit. We define a square computational domain of unit length and set wall boundary conditions everywhere. In particular, we fix a purely horizontal velocity at the top boundary $u_{1}=1,\, u_{2}=0$ and consider homogeneous no-slip boundary conditions on the bottom and lateral boundaries. As initial conditions, we consider a unit density, $\rho=1$, the pressure $\press=10^4$ and an initial fluid at rest. The viscosity is set to $\mu=10^{-2}$ so that $Re=100$ and $M \approx 8 \cdot 10^{-3}$ are the characteristic Reynolds and Mach numbers of this test attending to the lid velocity. \textcolor{black}{The artificial viscosity coefficient has been set to $c_{\alpha}=2$.} 
In the left plot of Figure \ref{LDC_figure}, we show the Mach contour plot of the solution obtained with the LADER-ENO scheme overlapped by a sketch of the half dual elements employed. 
The right plot reports the comparison between the approximated and the reference solution 
for the horizontal and vertical velocities along the vertical and horizontal 1D cuts in the middle of the domain. An almost perfect match is observed.

\subsection{Double shear layer} 
\begin{figure}[!htbp]
	\centering
	\includegraphics[trim= 100 100 100 100,clip,width=0.4\linewidth]{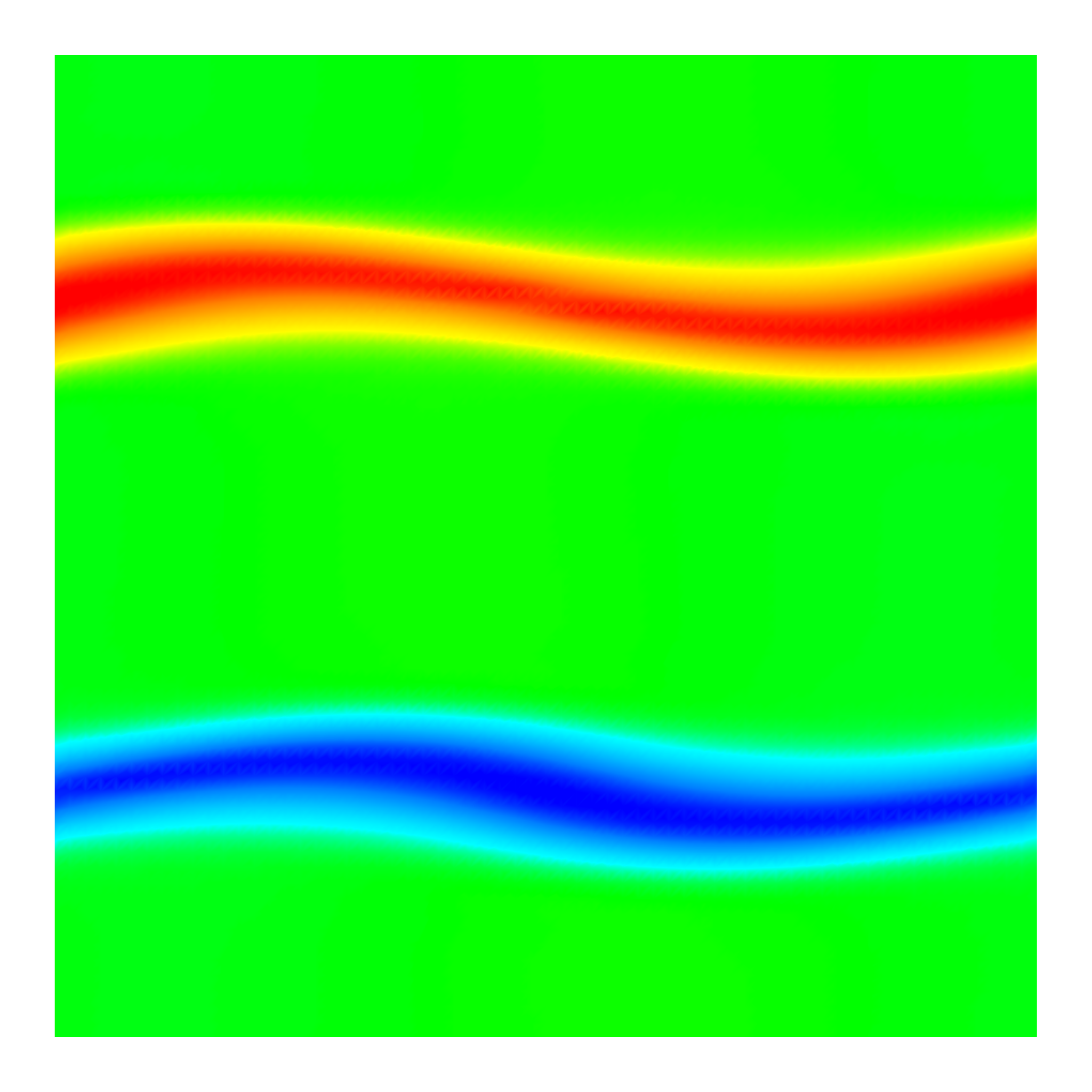}\hspace{0.05\linewidth}
	\includegraphics[trim= 100 100 100 100,clip,width=0.4\linewidth]{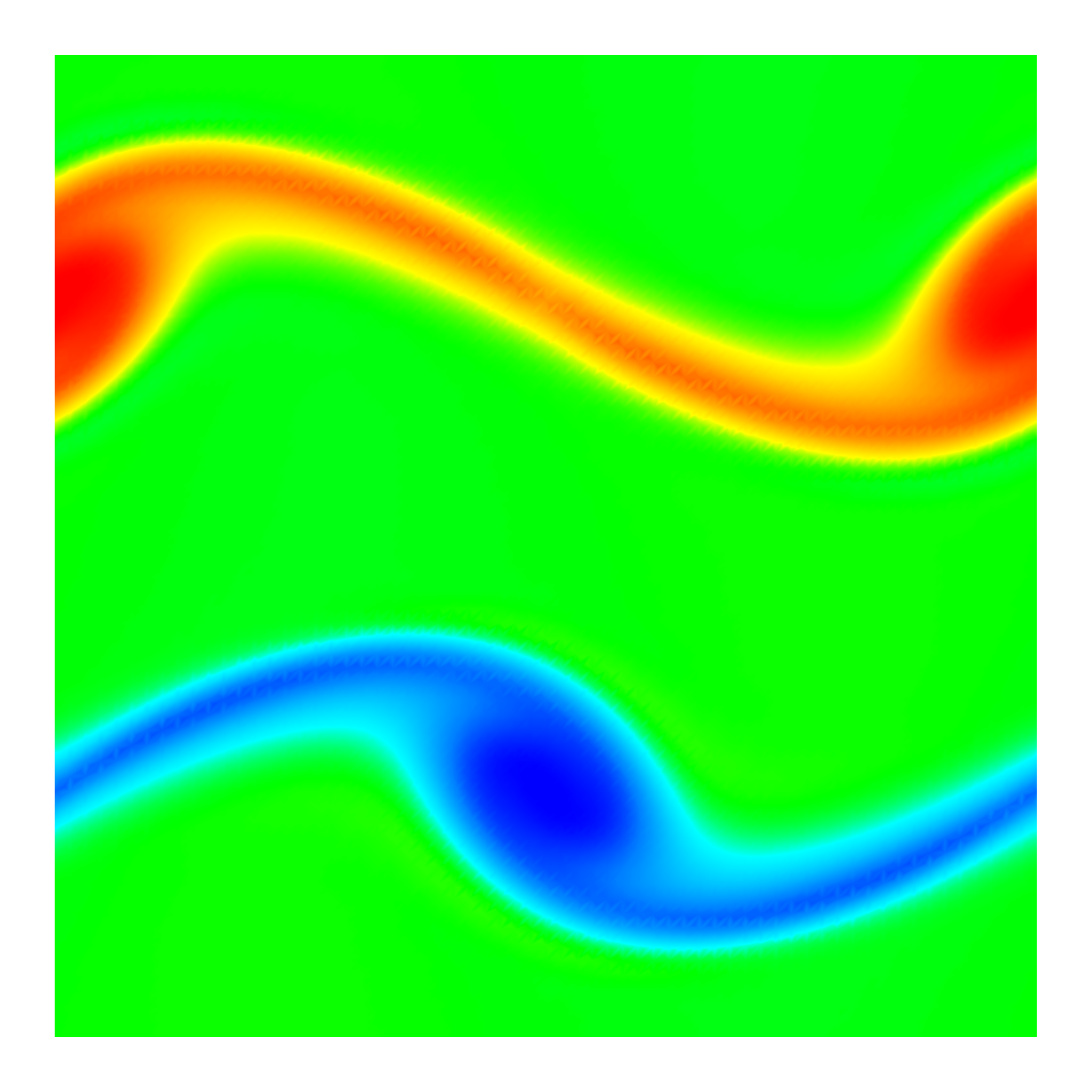}
	
	\vspace{0.05\linewidth}
	\includegraphics[trim= 100 100 100 100,clip,width=0.4\linewidth]{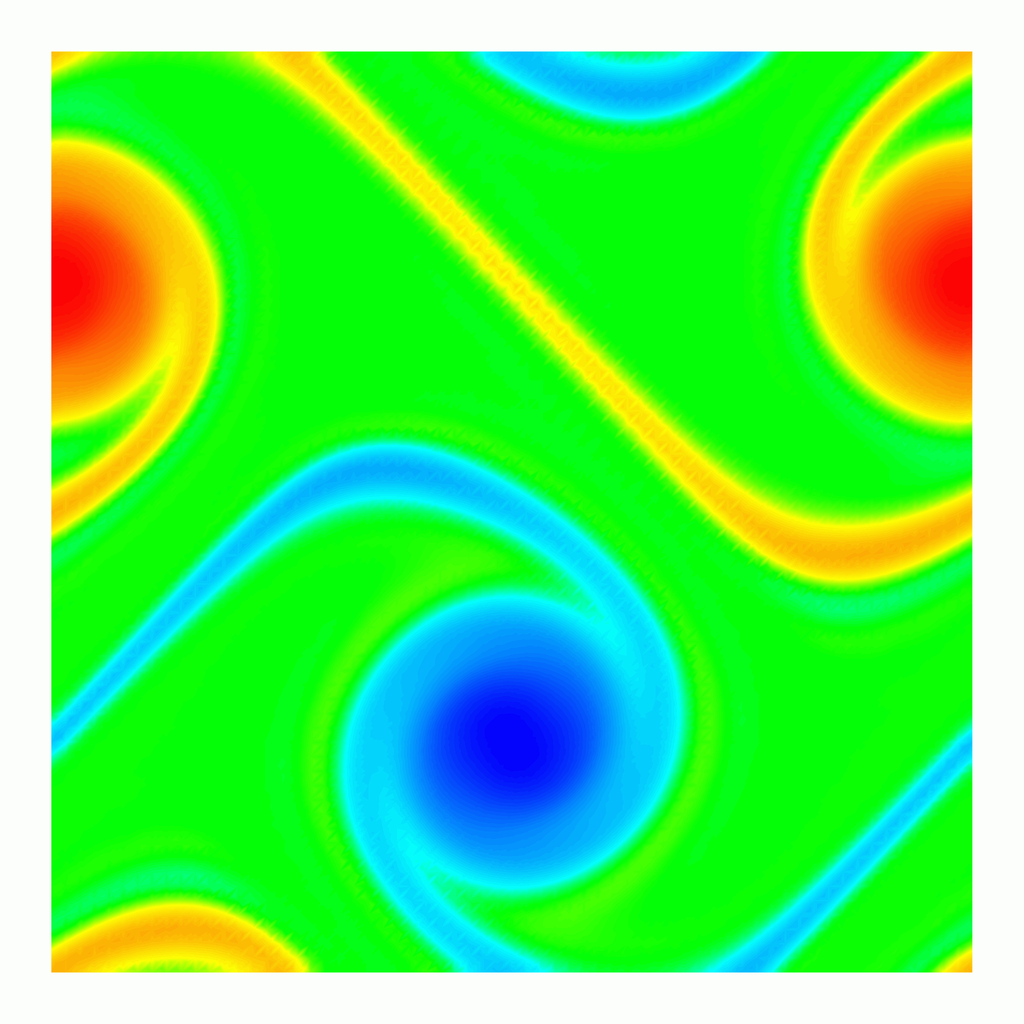}\hspace{0.05\linewidth}
	\includegraphics[trim= 100 100 100 100,clip,width=0.4\linewidth]{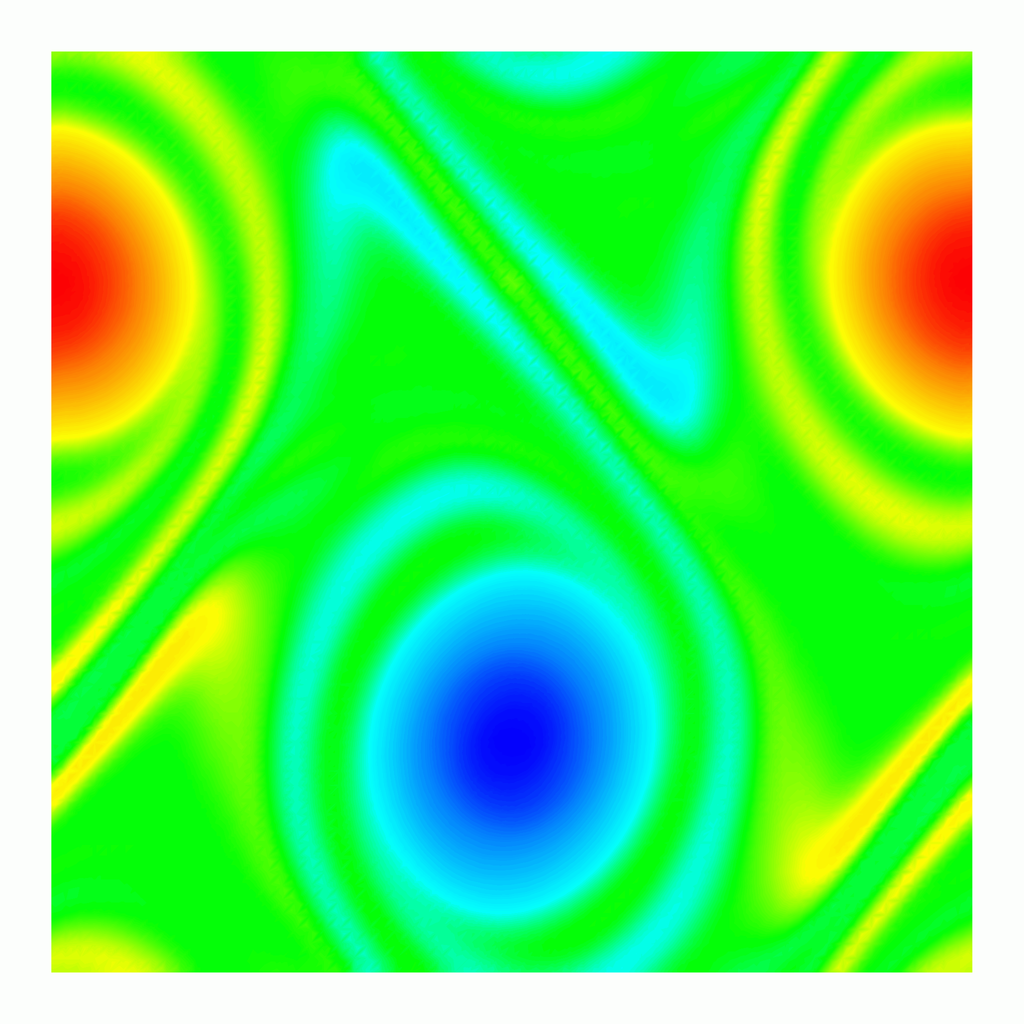}
	\caption{Contour plots of vorticity at times  $t\in \left\lbrace 0.8,1.6,2.4,3.6 \right\rbrace $ (from top left to bottom right) for the double shear layer test problem. The second order LADER scheme has been employed for the discretization of nonlinear convective terms. } 
	\label{DSL64_Lnl}
\end{figure}

In this section we apply the new hybrid finite volume / finite element method for all Mach number flows developed in this paper to the well-known double shear layer problem, see e.g. \cite{BCG89,TD15,Hybrid1}. The computational domain is given by $\Omega=[-1,1]^2$ and the initial condition reads  
\begin{equation}
\rho^{0}\left(\mathbf{x}\right) = 1,\qquad
{u}_{1}^{0} \left(\mathbf{x}\right) = \left\lbrace \begin{array}{lc}
\tanh \left[\hat{\rho}(\hat{y}-\frac{1}{4})\right] & \mathrm{ if } \; \hat{y} \le \frac{1}{2},\\[6pt]
\tanh \left[\hat{\rho}(\frac{3}{4}-\hat{y})\right] & \mathrm{ if } \; \hat{y} > \frac{1}{2},
\end{array}\right.\quad
{u}_{2}^{0} \left(\mathbf{x}\right) = \delta \sin \left(2\pi \hat{x}\right), \quad 
\press^{0} \left(\mathbf{x}\right) = \frac{10^{5}}{\gamma}, 
\end{equation}
with the abbreviations $\hat{x}= \frac{x+1}{2}$ and $\hat{y}= \frac{y+1}{2}$. The remaining parameters of the setup of this test case are chosen as $\hat{\rho} = 30$, $\mu= 2\cdot 10^{-4}$, $\lambda =0$ and $\delta =0.05$, see also \cite{TD15,Hybrid1}. The characteristic Mach number of this test case is $M \approx 2 \cdot 10^{-3}$, hence we are again in the low Mach number regime. The domain is covered with $8192$ primal elements and the boundary conditions are periodic everywhere. In Figure~\ref{DSL64_Lnl}, we show contour plots of the vorticity at times $t\in\left\lbrace 0.8, 1.6, 2.4, 3.6\right\rbrace$. 
Comparing our numerical solution with the one obtained in \cite{DPRZ16,TD15} we note a very good agreement, although the scheme presented in this paper is only second order accurate, while in \cite{DPRZ16,TD15} high order schemes have been employed.  

\subsection{Single Mach reflection problem} 

Let us now consider the single Mach reflection problem that can be found in \cite{Toro} and for which experimental reference data are available in \cite{Toro} and \cite{albumfluidmotion} under the form of Schlieren images. 
The test problem consists in a shock wave that is initially located in $x=-0.04$ and that travels to the 
right at a shock Mach number of $M=1.7$, hitting a wedge that forms an angle of $\varphi=25^\circ$ with the $x$-axis. 
The pressure, $\press$, and density, $\rho$, ahead of the shock are set to $\rho_0=1$ and $p_0 = 1/\gamma_g$, respectively, 
while for $x>-0.04$ we consider a fluid at rest. The post-shock values can be easily obtained 
from the Rankine-Hugoniot relations of the inviscid compressible Euler equations.  

\begin{figure}[!htbp]
	\begin{center} 
			\includegraphics[width=0.85\textwidth]{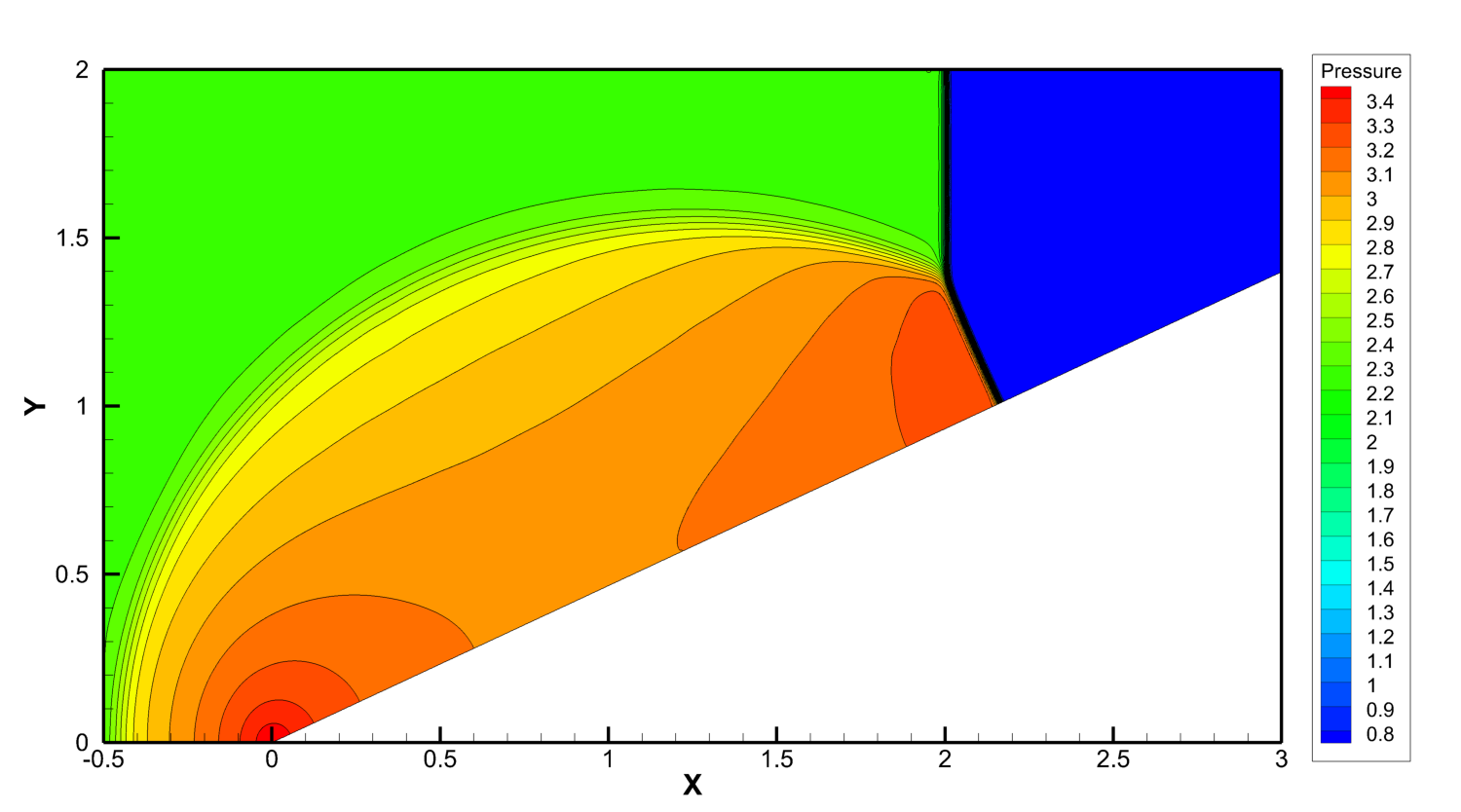}  
	\end{center} 
	\caption{Pressure contours obtained at $t_{\mathrm{end}}=1.2$ for the single Mach reflection problem using the new all Mach number hybrid FV/FE scheme presented in this paper. The shock is in the correct location at $x=2$. } 
	\label{fig.smr} 
\end{figure}

The computational domain is $\Omega = [0,3] \times [0,2]$, from which the $25^\circ$ wedge is subtracted. 
It is discretized using 1237328 triangular elements of characteristic mesh spacing $h=0.003$. \textcolor{black}{For this test we set $c_\alpha=1$.} 
The pressure field obtained for a simulation run until $t_{\mathrm{end}}=1.2$ is 
depicted in Figure~\ref{fig.smr}. 
The flow field obtained with the novel hybrid FV/FE scheme agrees well with the numerical and experimental reference solutions shown in \cite{Toro}. 
The shock wave is properly resolved and located in the correct position at $x=2$. 

\subsection{Shock-wedge interaction problem}
\label{sec.shockwedge}

\begin{figure}
	\begin{center} 
		\includegraphics[width=0.6\textwidth]{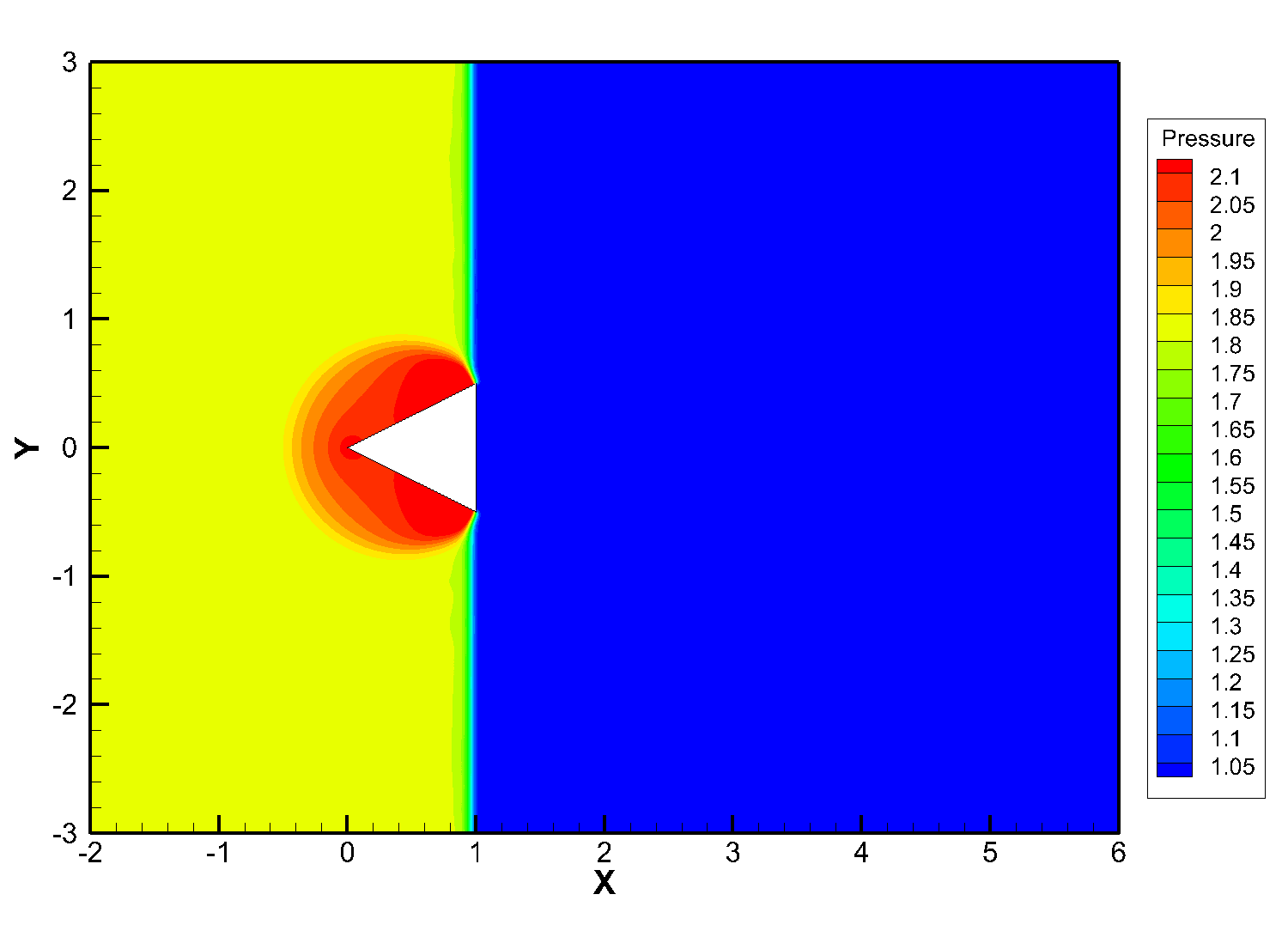}  
		\includegraphics[width=0.6\textwidth]{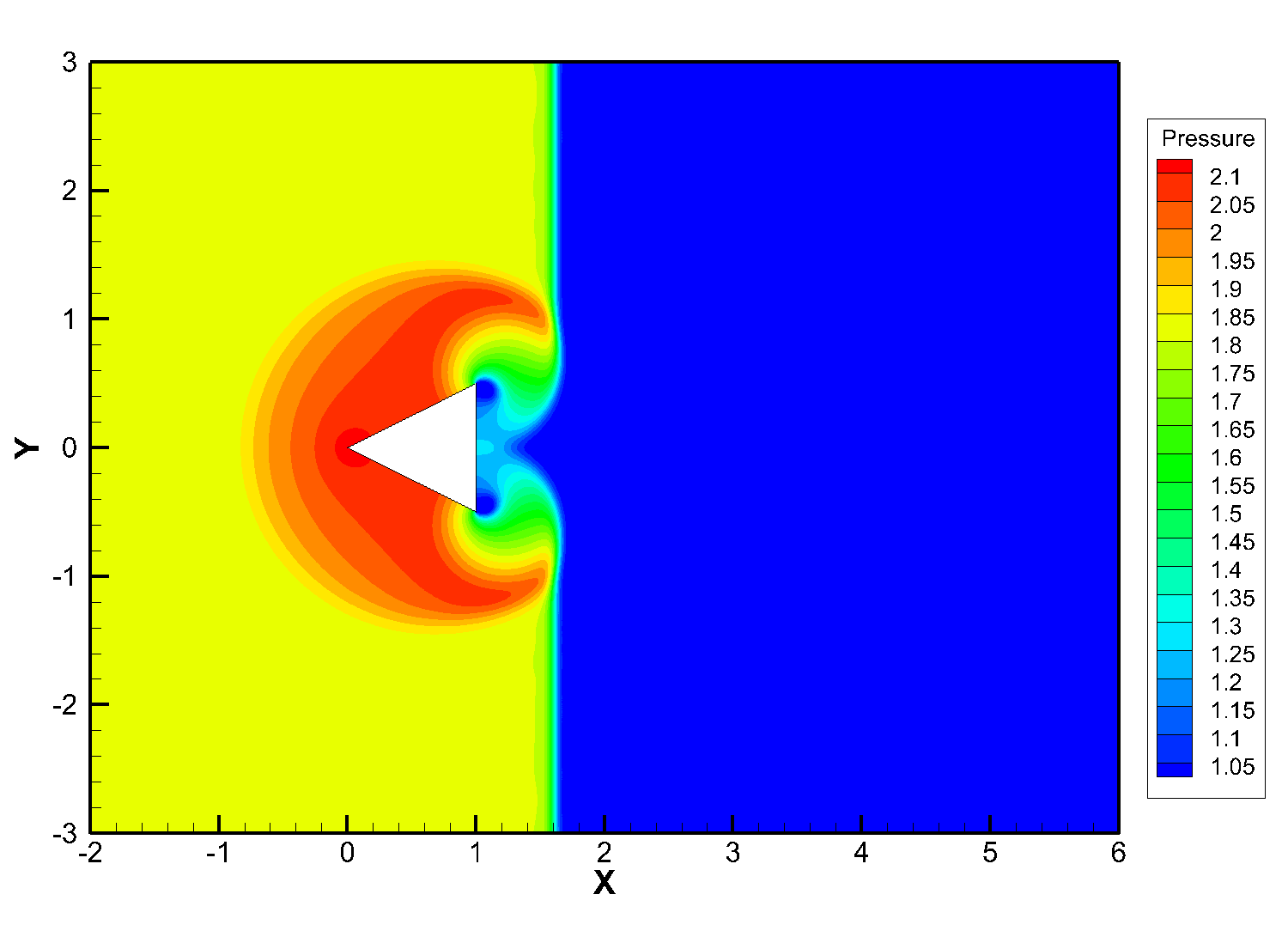}  
		\includegraphics[width=0.6\textwidth]{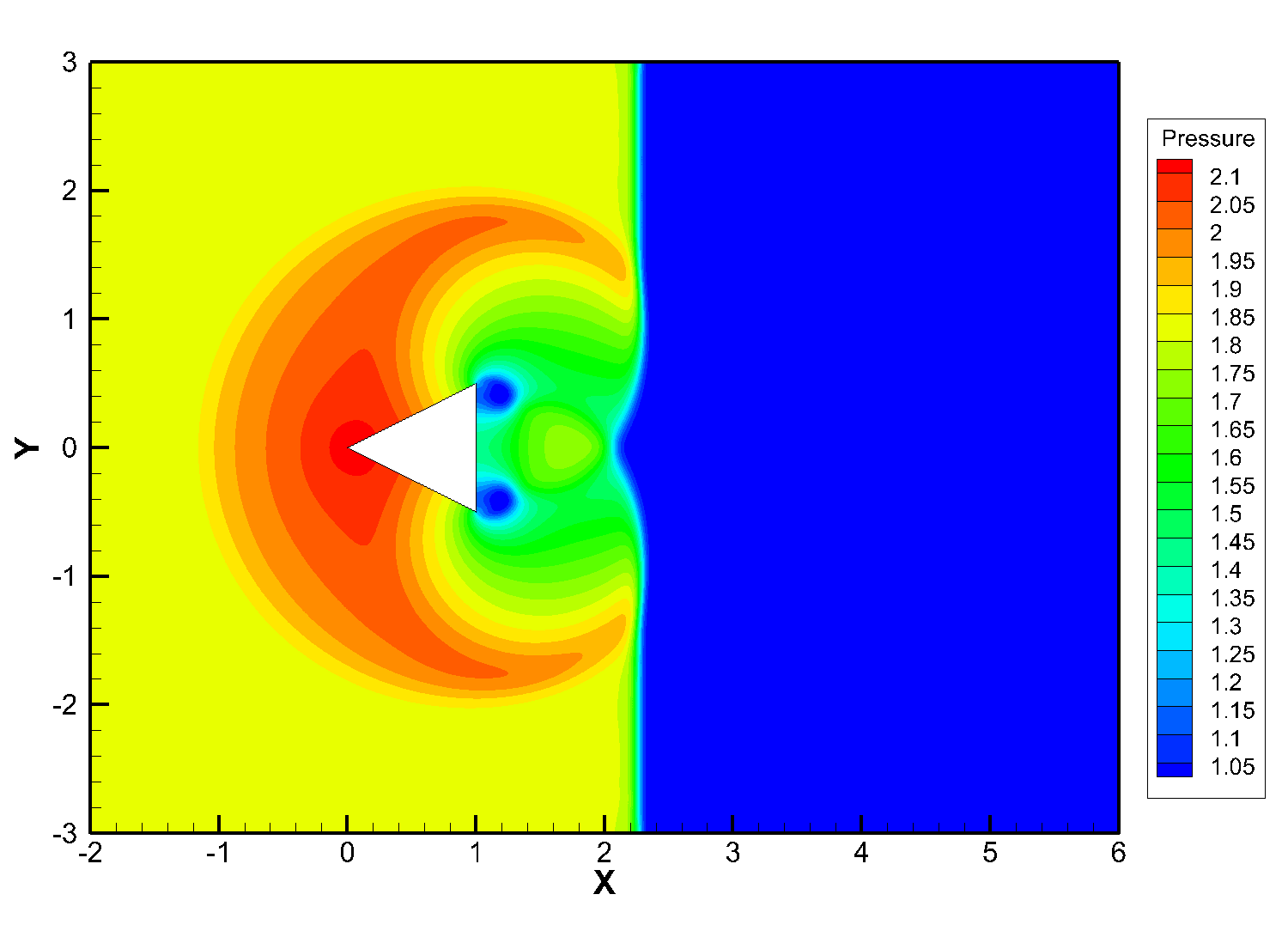}   
\end{center} 
	\caption{Shock-wedge interaction problem in 2D. Pressure contours 
		using the hybrid FV/FE method on an unstructured triangular 
		grid with mesh spacing $h=1/100$. Output times 
		from top to bottom: $t=1.5$, $t=2.0$, $t=2.5$.}
	\label{fig.shockwedge}
\end{figure}

In this section, we consider a flow that involves the interaction of a mild shock wave 
with a two-dimensional wedge, see also \cite{DumbserKaeser07,SolidBodies}. 
Experimental reference data for this test are available in form of Schlieren photographs, see 
\cite{albumfluidmotion,schardin}. 
The computational domain is given by $\Omega = [-2,6]\times[-3,3]$, excluding a wedge 
of length $L=1$ and height $H=1$ with its tip located in the origin. On all three edges 
of the wedge, we impose inviscid wall boundary conditions, while the upper and lower boundaries 
are periodic. On the left and on the right boundary, we impose the initial condition 
as Dirichlet boundary condition.  
The initial condition for a right-moving shock wave with shock Mach number $M_s=1.3$, initially located in $x=-1$, 
is setup according to the Rankine-Hugoniot relations, see \cite{DumbserKaeser07}. The pre-shock state (for $x>-1$) 
is given by $\rho_R = 1.4$, $\mathbf{u}_R=0$, and $p_R = 1.0$. 
A triangular mesh with a characteristic mesh spacing of $h=1/100$ is employed, leading to a total of 1080342 triangles. \textcolor{black}{For this test, we set $c_\alpha=0.5$.}  The pressure contours obtained with our hybrid FV/FE method are depicted in Figure \ref{fig.shockwedge} at 
several times. The location and shape of the shock and of the vortices shed behind the wedge compare qualitatively 
with those shown in \cite{DumbserKaeser07,SolidBodies,albumfluidmotion,schardin}. 

\subsection{Supersonic flow at $M=3$ over a circular blunt body} 

This last numerical test problem deals with the supersonic flow over a circular cylinder at Mach number $M=3$. The computational domain is the part for which $x\leq 0$ of a circle of radius $R=2$ centered in $\mathbf{x}_c=(0.5,0)$, from which a circular blunt body of radius $r_b=0.5$ centered in the origin is subtracted. The domain is discretized using 
a triangular mesh of characteristic mesh spacing $h=1/200$ composed of 348964 triangles. 
\textcolor{black}{For this test, we set $c_\alpha=3$.} 

\begin{figure}[h]
	\begin{center} 
		\begin{tabular}{c} 
			\includegraphics[width=0.5\textwidth]{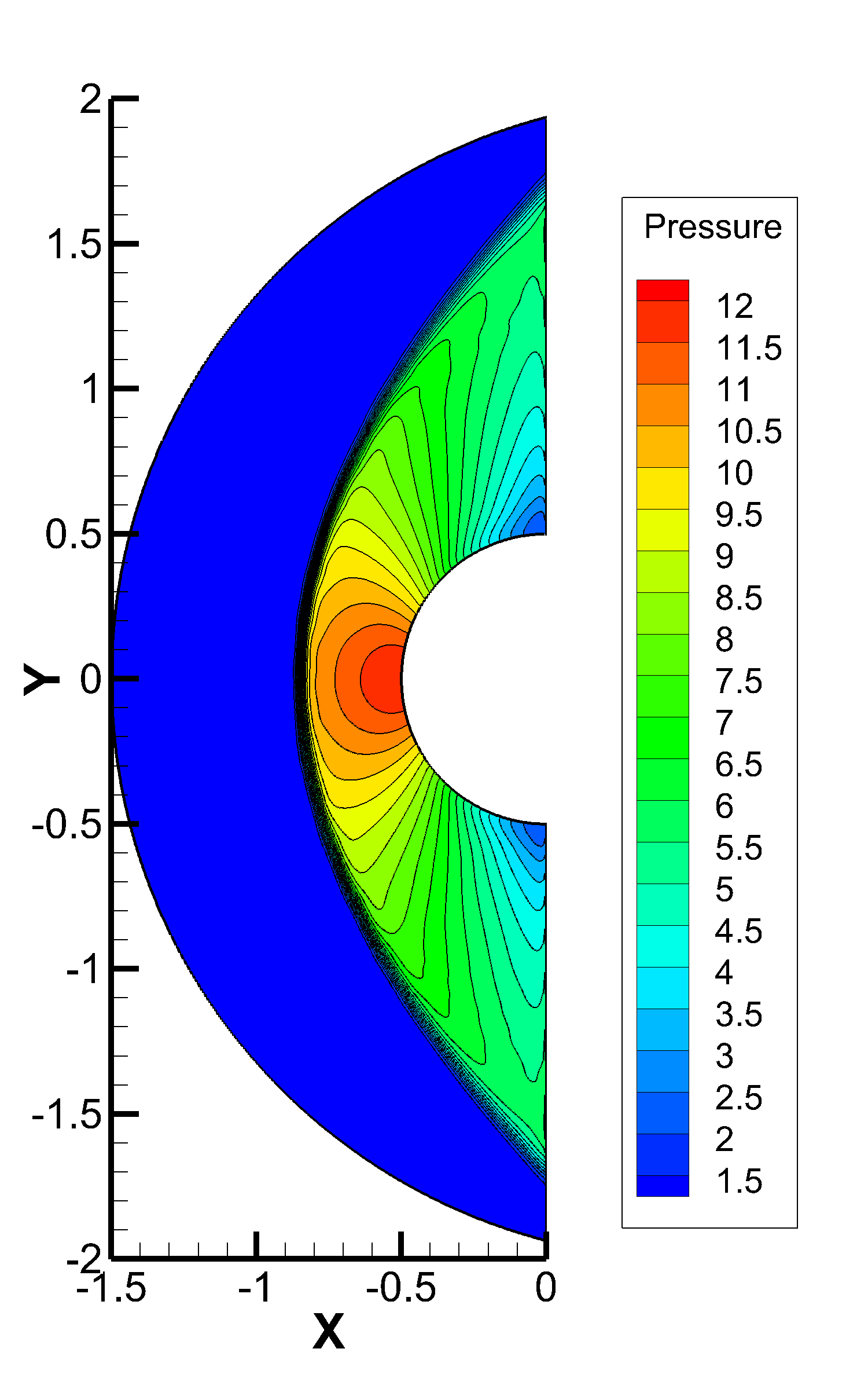}  
		\end{tabular}
	\end{center} 
	\caption{Pressure contours for the $M=3$ flow over a circular blunt body at time $t=5$. 
			}
	\label{fig.bluntbody}
\end{figure}

The initial condition is $\rho=1.4$, $\mathbf{u}=(3,0)$, $p=1$ in the entire computational domain. On the blunt body, inviscid wall boundary conditions are imposed. On the left inflow boundary, we impose the initial condition as Dirichlet boundary condition, while outflow is set on the right boundary. The computational results obtained with our hybrid FV/FE  scheme are depicted in Figure \ref{fig.bluntbody} at time $t=5$. The typical bow shock forms in front of the blunt body. 

\section{Conclusions}\label{sec:conclusions}
In this paper a novel asymptotic-preserving semi-implicit hybrid FV/FE algorithm has been proposed for the solution of all Mach number flows on staggered unstructured grids. 
The initial semi-discretization in time of Navier-Stokes equations allows a 
partial decoupling of the calculation of the linear momentum and of the density 
with respect to the solution of the pressure correction system. The first stage of the method 
involves the computation of the new density and an intermediate approximation 
of the linear momentum and total energy density, which account for the contribution 
of convection, diffusion, and gravity terms in the related conservative equations. 
Moreover, the pressure gradient at the previous time step is included so that 
the velocities would need only to be corrected with the pressure difference once 
it is computed. A local ADER methodology is employed to achieve a second order 
scheme in space and time, where we benefit from the dual mesh structure to 
approximate the gradients involved in the half in time reconstruction, hence reducing 
the stencil with respect to classical ADER methods. This procedure yields a 
good intermediate approximation of the linear momentum and total energy density 
to be provided for the computation of the pressure unknown on the projection stage. 
The splitting proposed following \cite{TV12,TD17} leads to an efficient numerical 
method in which the sound velocity is avoided on the eigenvalues computation of 
the transport-diffusion equations, approximated using an explicit scheme. Then, 
the pressure system is solved using classical implicit continuous finite element methods. 
Accordingly, the time step computed through the CFL condition is only limited by 
the flow velocity reducing the computational cost of the overall method. 
A key point of the proposed algorithm is the Picard iteration procedure 
that allows an iterative update of the linear momentum, enthalpy, and pressure 
variables avoiding the solution of a complex nonlinear system for the pressure. 
Once the pressure difference between two consecutive time steps is computed, 
the linear momentum and energy are corrected. The proposed methodology has been 
carefully validated by comparison of the obtained results with available 
analytical and numerical solutions. The numerical tests ranging from the 
incompressible limit to supersonic flows show the capability of 
the method to address complex flow phenomena. 

In the future, we plan to extend the hybrid FV-FE methodology in the context of shallow 
water equations making use of the seminal ideas presented in 
\cite{CasulliVOF,CasulliCheng1992, CasulliWalters2000,KramerStelling,TD14sw}. 
Moreover, more complex PDE systems, including natural involution constraints, 
like MHD equations, will be considered. To this end, the development of a 
structure-preserving scheme verifying the divergence-free condition will be essential, \cite{Pow97,MOSSV00,Bal03,BDA14,DBTF19}.

%
\section*{Acknowledgements}
This work was financially supported by the Italian Ministry of Education, University 
and Research (MIUR) in the framework of the PRIN 2017 project \textit{Innovative numerical methods for evolutionary partial differential equations and  applications} and via the  Departments of Excellence  Initiative 2018--2022 attributed to DICAM of the University of Trento (grant L. 232/2016). Furthermore, LR and MEV have received funding by Spanish MCIU under project MTM2017-86459-R and by FEDER and Xunta de Galicia funds under the ED431C 2017/60 project. SB was also funded by INdAM via a GNCS grant for young researchers and by a \textit{UniTN starting grant} of the University of Trento. SB, LR and MD are members of the GNCS group of INdAM.   

%

\bibliographystyle{plain}
\bibliography{./mibiblio}


\end{document}